\documentclass[11pt, a4paper, dvipsnames]{article}
\usepackage[margin=2cm]{geometry}
\usepackage[utf8]{inputenc}
\usepackage[T1]{fontenc}
\usepackage{amsmath,amssymb,amstext,amsthm}
\usepackage{booktabs}
\usepackage{graphicx}
\usepackage{hyperref}
\usepackage{subcaption}
\usepackage{array}
\usepackage{float}
\usepackage{longtable}
\usepackage[maxnames=199,backref=true]{biblatex}
\renewbibmacro{in:}{} 

\DeclareFieldFormat[article,inbook,incollection,inproceedings]
  {title:sentencecase}{\MakeSentenceCase*{#1}}
\DeclareFieldFormat{title:sentencecase}{#1} 
\DeclareFieldFormat{url}{%
  \mkbibacro{URL}\addcolon
  \ifkeyword{breakurl}{\newline}{\space}%
  \url{#1}} 

\renewbibmacro*{title}{%
  \ifboolexpr{
    test {\iffieldundef{title}}
    and
    test {\iffieldundef{subtitle}}
  }
    {}
    {\printtext[title]{%
       \printfield[title:sentencecase]{title}%
       \setunit{\subtitlepunct}%
       \printfield[title:sentencecase]{subtitle}}}%
  \newunit}

\DeclareCaptionLabelFormat{figsub}{Figure~\thefigure#2} 
\usepackage{pgfplots}
\pgfplotsset{compat=1.16}
\pgfplotstableread[col sep=comma]{data/length-fit-data.csv}{\datatable} 
\pgfplotstableread[col sep=comma]{data/angle-hist-100.csv}{\angletable}

\usepackage{tikz}
\usetikzlibrary{patterns}
\definecolor{colorF}{HTML}{F4C800}
\definecolor{colorI}{HTML}{803E75}
\definecolor{colorL}{HTML}{FF6800}
\definecolor{colorN}{HTML}{C10020}
\definecolor{colorP}{HTML}{007D34}
\definecolor{colorT}{HTML}{F6768E}
\definecolor{colorU}{HTML}{00A1C2}
\definecolor{colorV}{HTML}{93AA00}
\definecolor{colorW}{HTML}{FF1493}
\definecolor{colorX}{HTML}{4363D8}
\definecolor{colorY}{HTML}{8F3931}
\definecolor{colorZ}{HTML}{A0A0A0}

\usepackage{algorithm}
\usepackage{algpseudocode}

\newtheorem{definition}{Definition}
\newtheorem{theorem}{Theorem}
\newtheorem{corollary}{Corollary}

\newtheorem{lemma}{Lemma}
\newcommand{\sbinom}[2]{\genfrac{(}{)}{0pt}{2}{#1}{#2}}

\newcommand\maxpolygon{\textsc{MaxPolygon}}

\makeatletter
\renewcommand\@maketitle{%
  \newpage
  \null
  \vfill                
  \begin{center}%
  \let \footnote \thanks
    {\LARGE \@title \par}%
    \vskip 1.5em%
    {\large
      \lineskip .5em%
      \begin{tabular}[t]{c}%
        \@author
      \end{tabular}\par}%
    \vskip 1em%
    {\large \@date}%
  \end{center}%
  \par
  \vskip 1.5em}
\makeatother

\begin{document}

\title{Optimal Polyomino Fences}
\author{%
  \parbox{\dimexpr\textwidth-2\tabcolsep\relax}{\centering
    Alexis Langlois-Rémillard, Mia N. Mü\ss{}ig,
    Helmut Podhaisky, Erika Roldan, Günter Rote
\\[1.5ex]
 Maximilian Brömme\thanks{Part of the student team at Freie
   Universität Berlin that worked on the 3D walls},
 Matthes Norden\footnotemark[\value{footnote}],
 Maximilian Prietzel\footnotemark[\value{footnote}],
 Daniel Yu\footnotemark[\value{footnote}],
 Mykhailo Lyader\thanks{Worked with Helmut Podhaisky as part of a
   school project
 }}}

\date{November 12, 2026}

\AddToHook{enddocument}[addresses]
{
  \vskip1cm plus 2cm \goodbreak\bigskip{\footnotesize \parindent0pt
    
(A.~Langlois-Rémillard) \textsc{Aalto University, Otakaari 1 F
02150, Espoo, Finland}\par\nobreak
\textit{Email address:} \texttt{alexis.langlois-remillard@aalto.fi} \medskip

(M.~N.~Mü\ss{}ig) \textsc{Ludwig-Maximilians-Universität München, Germany}\par
\textit{Email address:} \texttt{nienna@miamuessig.de}\medskip

(H.~Podhaisky) \textsc{Martin-Luther-Universität Halle-Wittenberg, Germany}\par
\textit{Email address:} \texttt{helmut.podhaisky@mathematik.uni-halle.de}\medskip

(G.~Rote) \textsc{Freie Universität Berlin, Germany}\par
\textit{Email address:} \texttt{rote@inf.fu-berlin.de}\medskip

(E.~Roldan) \textsc{Max Planck Institute for Mathematics in the Sciences, Leipzig, Germany; ScaDS.AI, Leipzig University, Germany; Université de Genève, Switzerland}\par
\textit{Email addresses:} \texttt{roldan@mis.mpg.de}, \texttt{erika.roldan@unige.ch}\medskip

(M.~Norden) \textsc{Freie Universität Berlin, Germany}\par
\textit{Email address:} \texttt{nordstit@gmail.com}
\medskip

(M.~Brömme) \textsc{Freie Universität Berlin, Germany}\par
\textit{Email address:} \texttt{maximilian.broemme@aviomatik.de}\medskip

(M.~Prietzel) \textsc{Freie Universität Berlin, Germany}\par
\textit{Email address:} \texttt{maximilian.prietzel@proton.me}\medskip

(D.~Yu) \textsc{Freie Universität Berlin, Germany}\par
\textit{Email address:} \texttt{barzindyne@gmail.com}\medskip

(M.~Lyader) \textsc{Georg-Cantor-Gymnasium Halle, Germany}\par
\textit{Email address:} \texttt{mykhailolyader@gmail.com}\medskip

}}
\DeclareHookRule{enddocument}{addresses}{before}{*}

\maketitle

\begin{abstract}
    The $n$-omino fence problem asks how many cells can be surrounded by placing every $n$-omino exactly once in the grid. We present a new dynamic programming approach for this problem, allowing us to prove optimality for the pentomino and hexomino fence and to classify all possible pentomino orders by optimal area. For the first time we also provide strong heuristic solutions and upper bounds for all fences up to $n = 15$, after which we even approximate the $500$-omino fence with computational experiments. Finally, we generalize the problem to 3D, where we prove optimality for pentacubes and derive heuristic solutions for hexacubes and heptacubes.
\end{abstract}
\vfill\vfill 
\clearpage  

\tableofcontents

\section{Introduction}\label{sec:introduction}

We consider the square grid $\mathbb{Z}^2$, whose cells are unit squares. Two cells are \emph{4-adjacent} if they share an edge, and \emph{8-adjacent} if they share an edge or a vertex. A set of cells is \emph{4-connected} (resp.\ \emph{8-connected}) if its cells form a connected graph under 4-adjacency (resp.\ 8-adjacency). The \emph{free $n$-ominoes} are the 4-connected sets of $n$ cells, where two sets are considered the same if one can be obtained from the other by translation, rotation, and reflection. There are $12$, $35$, $108$, and $369$ free $n$-ominoes for $n = 5, 6, 7, 8$, respectively.

The \emph{$n$-omino fence problem} asks: what is the largest area (number of cells) that can be enclosed by placing every $n$-omino exactly once in the grid? The pieces must not overlap, and they must touch along full edges, not only at corners. We give a precise definition of the problem in Section~\ref{sec:precise-def}.

Victor Feser first presented such fence-like puzzles with pentominoes
in 1968~\cite{Fe68}. One of them, which asks for the largest area
enclosed by a fence of the 12 pentominoes, became widely known when it
appeared in Martin Gardner's \emph{Mathematical Games} column
in the \emph{Scientific American}
in
1973~\cite{Ga73}. Gardner presented his best attempt, enclosing 127
squares, and
reported that Donald E. Knuth had found a
fence of area 128, 
together with an informal proof that 128 cannot be exceeded.
In 1978, Takakazu Shimauchi published a more rigorous argument for 128
being optimal~\cite{Ta78}. In a 2025 paper by a subset of the
authors~\cite{previous-german} this result was reconfirmed with an
integer linear-programming (ILP) model, which was solved by the solver
Gurobi.
They also found an improved solution of area $1586$ (reported originally as $1588$ because of a counting error) for the hexomino fence.

In this paper, we continue our work on the problem. After the problem definition in Section~\ref{sec:precise-def}, we start with the pentomino fence in Section~\ref{sec:pentomino}. We have taken great care to make the presented dynamic programming (DP) approach accessible to a general mathematically interested audience, with detailed explanations and extensive visualizations. Additionally, the code verifying pentomino optimality does not rely on external solvers and runs in a few minutes even on low-budget hardware.

In Section~\ref{sec:hexomino}, we look at hexomino
fences, whose maximum area was unknown until now. We use an optimized
version of the dynamic programming approach to prove that $1597$ is
the maximum area, and we identify all optimal hexomino
fences. We proceed to heptominoes and octominoes,
without reaching optimality (Section~\ref{sec:hepto-okto}), and then far beyond (Section~\ref{sec:larger-fences}). Finally, in Section~\ref{sec:3d}, we turn to three dimensions: we determine the maximum enclosable volume with pentacubes and find heuristic solutions for hexacubes and heptacubes.

\paragraph{Repository with programs and results.} The programs and result files of our calculations are published on GitHub~\cite{github}. We also provide an interactive exploration of the main results~\cite{FenceChallenge}. The workflow of the computations, including the programs and the files that they exchange, is described in the README file of the repository.

\section{Precise definition of the problem}
\label{sec:precise-def}

\begin{figure}[htb]
  \centering
  \begin{tikzpicture}[scale=0.45, y=-1cm]
    \definecolor{pieceA}{HTML}{F4C800}
    \definecolor{pieceB}{HTML}{B785AD}
    \definecolor{pieceC}{HTML}{FFBC3A}
    \definecolor{pieceD}{HTML}{E86A7A}
    \definecolor{pieceE}{HTML}{5DBB7F}
    \definecolor{pieceF}{HTML}{F6768E}
    \definecolor{pieceG}{HTML}{5CC6DC}
    \definecolor{pieceH}{HTML}{B5CC4D}
    \definecolor{pieceI}{HTML}{FF74BF}
    \definecolor{pieceJ}{HTML}{8FA3EA}
    \definecolor{pieceK}{HTML}{C98A84}
    \definecolor{pieceL}{HTML}{B8B8B8}
    \definecolor{pieceM}{HTML}{9C86C2}
    \definecolor{pieceN}{HTML}{CEA262}
    \definecolor{pieceO}{HTML}{B5A89F}
    \definecolor{pieceP}{HTML}{F7866B}
    \definecolor{pieceQ}{HTML}{FFD08A}
    \definecolor{pieceR}{HTML}{A8E6A1}
    \definecolor{pieceS}{HTML}{C9B6F2}
    \definecolor{pieceT}{HTML}{FFB38A}
    \definecolor{pieceU}{HTML}{9ED9D0}
    \foreach \tx/\ty in {
      2/7, 3/7, 3/9, 3/10, 4/7, 4/8, 4/9, 4/10, 5/9, 5/10, 6/9, 6/10,
      6/11, 6/12, 7/9, 8/2, 8/3, 8/4, 8/5, 8/6, 8/7, 8/8, 8/11, 8/12,
      9/3, 9/4, 9/6, 9/7, 9/8, 9/9, 9/10, 9/11, 9/12, 10/3, 10/6, 10/7,
      10/8, 10/9, 10/10, 10/12, 11/3, 11/4, 11/5, 11/6, 11/7, 11/8, 11/9, 11/10,
      12/7, 12/9, 12/10
    } {
      \fill [pattern=north east lines, pattern color=black!70] (\tx,\ty) rectangle ++(1,1);
    }
    \foreach \tx/\ty in {
      1/4, 2/2, 2/5, 3/1, 3/2, 3/3, 3/4, 3/5, 3/12, 3/13, 3/14, 4/1,
      4/2, 4/3, 4/4, 4/5, 4/12, 4/14, 5/2, 5/3, 5/4, 5/5, 6/2, 6/3,
      6/4, 6/5, 6/7, 13/4, 14/4, 14/5, 14/7, 14/8, 14/9, 15/4, 15/5, 15/9,
      16/5, 16/6, 16/9, 17/4, 17/5, 17/6, 17/7, 17/8, 17/9, 17/10, 18/4, 18/5,
      19/10, 20/10
    } {
      \fill [pattern=crosshatch, pattern color=gray!50] (\tx,\ty) rectangle ++(1,1);
    }
    \foreach \x/\y in {1/0, 2/0, 3/0, 1/1, 2/1, 1/2}
      \draw[fill=pieceA] (\x,\y) rectangle ++(1,1);
    \foreach \x/\y in {4/0, 5/0, 5/1, 6/1}
      \draw[fill=pieceB] (\x,\y) rectangle ++(1,1);
    \foreach \x/\y in {7/1, 7/2, 7/3, 7/4, 7/5, 1/6, 2/6, 3/6, 4/6, 5/6, 6/6, 7/6, 5/7, 7/7, 5/8, 6/8, 7/8}
      \draw[fill=pieceC] (\x,\y) rectangle ++(1,1);
    \foreach \x/\y in {8/1, 9/1, 10/1, 9/2, 10/2, 11/2}
      \draw[fill=pieceD] (\x,\y) rectangle ++(1,1);
    \foreach \x/\y in {12/2, 12/3, 12/4, 12/5, 12/6, 13/5, 13/6, 13/7, 12/8, 13/8, 13/9, 13/10}
      \draw[fill=pieceE] (\x,\y) rectangle ++(1,1);
    \foreach \x/\y in {13/3, 14/3, 15/3, 16/3, 16/4}
      \draw[fill=pieceF] (\x,\y) rectangle ++(1,1);
    \foreach \x/\y in {17/3, 18/3, 19/3, 19/4, 19/5}
      \draw[fill=pieceG] (\x,\y) rectangle ++(1,1);
    \foreach \x/\y in {18/6, 19/6, 20/6, 18/7, 20/7, 18/8, 21/8, 18/9, 19/9, 20/9, 21/9, 18/10, 21/10, 18/11, 19/11, 20/11, 21/11}
      \draw[fill=pieceH] (\x,\y) rectangle ++(1,1);
    \foreach \x/\y in {14/10, 15/10, 16/10, 16/11, 17/11}
      \draw[fill=pieceI] (\x,\y) rectangle ++(1,1);
    \foreach \x/\y in {10/11, 11/11, 12/11, 13/11, 11/12}
      \draw[fill=pieceJ] (\x,\y) rectangle ++(1,1);
    \foreach \x/\y in {11/13, 12/13, 11/14, 12/14, 12/15}
      \draw[fill=pieceK] (\x,\y) rectangle ++(1,1);
    \foreach \x/\y in {6/13, 7/13, 8/13, 9/13, 10/13, 7/12,
      7/11, 7/10, 8/10, 8/9}
      \draw[fill=pieceL] (\x,\y) rectangle ++(1,1);
      \foreach \x/\y in {3/11,
        4/11, 5/11, 5/12, 4/13, 5/13}
      \draw[fill=pieceM] (\x,\y) rectangle ++(1,1);
    \foreach \x/\y in {2/10, 2/11, 2/12, 2/13, 2/14, 2/15, 3/15, 4/15, 5/15, 5/14}
      \draw[fill=pieceN] (\x,\y) rectangle ++(1,1);
    \foreach \x/\y in {1/7, 1/8, 2/8, 3/8, 2/9}
      \draw[fill=pieceO] (\x,\y) rectangle ++(1,1);
    \foreach \x/\y in {0/3, 1/3, 2/3, 0/4, 2/4, 0/5, 1/5}
      \draw[fill=pieceP] (\x,\y) rectangle ++(1,1);
    \foreach \x/\y in {14/0, 15/0, 16/0, 17/0, 15/1, 15/2, 17/1}
      \draw[fill=pieceQ] (\x,\y) rectangle ++(1,1);
    \foreach \x/\y in {18/0, 19/0, 20/0, 21/0, 21/1, 20/2, 21/2, 21/3}
      \draw[fill=pieceR] (\x,\y) rectangle ++(1,1);
    \foreach \x/\y in {14/6, 15/6, 15/7, 16/7, 15/8, 16/8}
      \draw[fill=pieceT] (\x,\y) rectangle ++(1,1);
    \foreach \x/\y in {10/4, 9/5, 10/5}
      \draw[fill=pieceU] (\x,\y) rectangle ++(1,1);
    \node[font=\scriptsize] at (2.5,1.5) {\textit{A}};
    \node[font=\scriptsize] at (5.5,0.5) {\textit{B}};
    \node[font=\scriptsize] at (5.5,6.5) {\textit{C}};
    \node[font=\scriptsize] at (9.5,1.5) {\textit{D}};
    \node[font=\scriptsize] at (13.5,6.5) {\textit{E}};
    \node[font=\scriptsize] at (15.5,3.5) {\textit{F}};
    \node[font=\scriptsize] at (18.5,3.5) {\textit{G}};
    \node[font=\scriptsize] at (19.5,9.5) {\textit{H}};
    \node[font=\scriptsize] at (16.5,10.5) {\textit{I}};
    \node[font=\scriptsize] at (11.5,11.5) {\textit{J}};
    \node[font=\scriptsize] at (12.5,14.5) {\textit{K}};
    \node[font=\scriptsize] at (8.5,13.5) {\textit{L}};
    \node[font=\scriptsize] at (5.5,12.5) {\textit{M}};
    \node[font=\scriptsize] at (2.5,13.5) {\textit{N}};
    \node[font=\scriptsize] at (2.5,8.5) {\textit{O}};
    \node[font=\scriptsize] at (0.5,4.5) {\textit{P}};
    \node[font=\scriptsize] at (15.5,0.5) {\textit{Q}};
    \node[font=\scriptsize] at (21.5,1.5) {\textit{R}};
    \node[font=\scriptsize] at (15.5,7.5) {\textit{S}};
    \node[font=\scriptsize] at (10.5,5.5) {\textit{T}};
  \end{tikzpicture}
  \caption{An arrangement of $20$ polyomino pieces of different sizes,
    labeled \textit{A} to~\textit{T}, which illustrates the definitions of
    Section~\ref{sec:precise-def}.
    The largest single region, which determines the area ($51$~cells), is hatched.
    The other five regions ($21$, $21$, $5$, $2$, and $1$~cells) are
    cross-hatched; they do not count.
    The white cells surrounded by \textit{F}, \textit{Q}, \textit{R}, and \textit{G},
    and the three white cells in the upper part of~\textit{H}
    are not enclosed, since they are 8-connected to the outside.
    Because of the free-floating triomino~\textit{T},
    the arrangement is not a valid fence.}
  \label{fig:weird}
\end{figure}
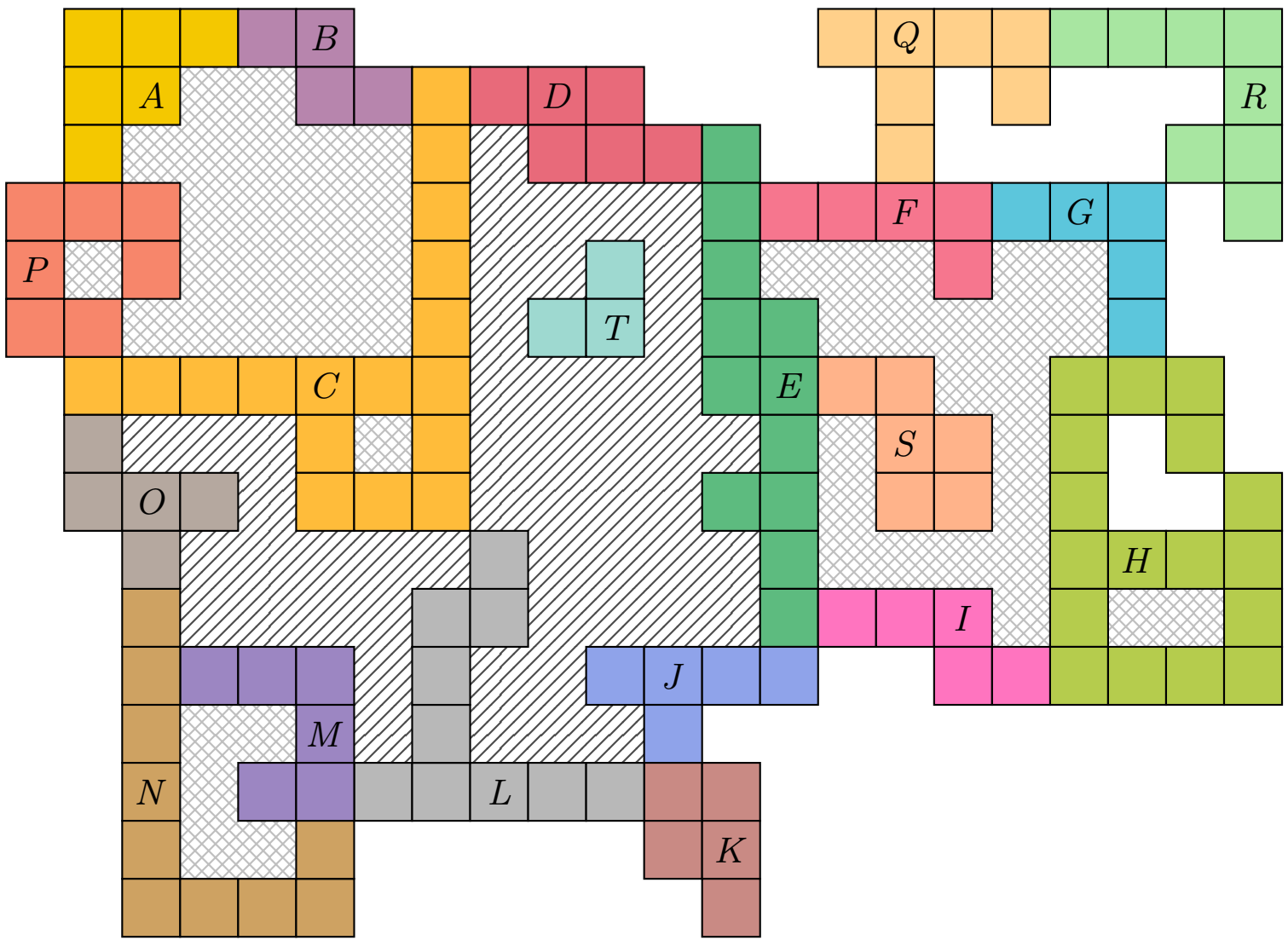

The informal description in the introduction leaves some questions
open that must be settled before a computer search can claim to have
found an optimal fence: Which cells count as enclosed? What if the
enclosed cells fall apart into several parts? What happens with pieces
that lie inside the enclosed area? In this section, we give precise
definitions. \autoref{fig:weird} shows an arrangement of pieces of
various sizes that exhibits many of the fine points. It is far from
optimal, and because of the free-floating triomino~\textit{T}, it is
not even a valid fence. Nevertheless, we will refer to it throughout
this section.

\subsection{Fences and their area}
\label{sec:def-area}

Recall from the introduction that two cells are \emph{4-adjacent} if
they share an edge, and \emph{8-adjacent} if they share an edge or a
vertex. We are given a finite collection of polyominoes, the
\emph{pieces}. In the $n$-omino fence problem, these are the free
$n$-ominoes, each of them exactly once, but the definitions make sense
for any collection of polyominoes, as in \autoref{fig:weird}. A
\emph{placement} puts every piece into the grid, after an arbitrary
rotation, reflection, and translation.

\begin{definition}\label{def:fence}
  A placement is a \emph{fence} if the placed pieces are pairwise
  disjoint and their union~$\Phi$ is 4-connected.
  The cells of~$\Phi$ are the \emph{fence cells}, and all other cells
  are \emph{free}. Since $\Phi$ is finite, exactly one 8-connected
  component of free cells is infinite; we call it the \emph{outside}.
  The free cells that do not belong to the outside are
  \emph{enclosed}, and the 8-connected components of the set of
  enclosed cells are the \emph{regions} of the fence.

  We define
  the \emph{area} of the fence as the number of cells of a largest
  region, or $0$ if there is no region.
\end{definition}

The \emph{$n$-omino fence problem} asks for the largest area of a
fence built by using each free $n$-omino exactly once.

\paragraph{Notions of connectivity.}
Free cells are regarded as connected even if they touch only at a
corner, whereas pieces must share an edge in order to be connected.
These two conventions fit together, by a basic fact of digital  
topology~\cite{Rosenfeld79,KongRosenfeld89}.
Enclosing an
8-connected set~$X$ of cells and cutting it off from other 8-connected regions
requires a 4-connected barrier (and vice versa).
Thus, the pieces of such a
wall must be joined along edges. \autoref{fig:weird} shows the
consequences:
\begin{itemize}
\item The eight white cells surrounded by \textit{F}, \textit{Q},
  \textit{R}, and \textit{G} look as if they were enclosed, but they
  are 8-connected to the outside through the diagonal gap where
  \textit{G} and \textit{R} touch at a corner. For the same reason,
  the three white cells in the upper part of~\textit{H} are not
  enclosed: They reach the outside through a diagonal gap between two
  cells of the piece~\textit{H} itself.
\item Conversely, enclosed cells that touch only at a corner belong to
  the same region. The hatched region consists of a part of $15$ cells
  on the left and a part of $36$ cells in the middle, which are
  connected only through the diagonal gap between \textit{C}
  and~\textit{L}. Together, they form a single region of $51$ cells.
  Similarly, the cell that is surrounded by~\textit{P} on all four
  sides belongs to the region at the upper
  left, since it touches this region at a corner.
\end{itemize}

\paragraph{Connected area.}
In our definition we only count the largest region as the area of the fence. The arrangement in \autoref{fig:weird} has six regions, with $51$,
$21$, $21$, $5$, $2$, and $1$ cells, respectively, and its area
according to our definition
is~$51$. The smaller regions are
cross-hatched: the regions of $21$
cells at the upper left and at the right, the region of $5$ cells
between \textit{M} and~\textit{N}, and the holes inside the pieces
\textit{C} and~\textit{H}. Counting only the largest region captures the
intention of the original puzzle, which asks for a single field of a
``pentomino farm''~\cite{Fe68}.
In particular, a hole inside a single piece
will hardly ever determine the area.
This becomes relevant
when we consider the octominoes in \autoref{sec:hepto-okto-lower}.
We discuss the alternative of counting the total enclosed
area in
\autoref{sec:total-area}.

\paragraph{Connectivity of the fence, and unused pieces.}
The free-floating triomino~\textit{T} of \autoref{fig:weird} violates the requirement that
$\Phi$ must be 4-connected. It lies inside the hatched region, but its
cells are not free and do not count: The hatched region has $51$
cells, although its outer boundary encloses $54$ cells. Moving
\textit{T} to the outside would enlarge the region to $54$~cells.

The requirement of 4-connectivity does not change the maximum area.
To see this, consider any placement of pairwise disjoint pieces, not
necessarily 4-connected, with a largest region~$\mathcal R$. By the
fact quoted above, the fence cells that are 8-adjacent to~$\mathcal R$
and lie in the unbounded 4-connected component of the complement
of~$\mathcal R$ form a 4-connected set. Hence they belong to a single
4-connected component~$K$ of~$\Phi$, which is a union of pieces. We
remove all pieces that do not belong to~$K$ and attach them, one after
the other, to the outside of the fence, for example directly on top of
the currently highest fence cell. An attached piece lies completely in
the outside, and therefore it cannot destroy or shrink a region.
The cells of the removed pieces that lay inside~$\mathcal R$, like
\textit{T}, become free and join the region, because they are
surrounded by~$\mathcal R$ and~$K$. The result is a fence whose largest region
contains~$\mathcal R$.

By the same procedure, pieces that play no role in enclosing the
largest region can always be moved to the outside. In the $n$-omino
fence problem, where every piece must be used, such pieces are thus
harmless: effectively, the problem allows using an arbitrary subset
of the pieces.

\subsection{The boundary around the largest region}
\label{sec:pieces-around}

\begin{figure}[ht]
    \centering
    \begin{tikzpicture}[scale=0.35, y=-1cm]
      \foreach \tx/\ty in {
        2/2, 3/2, 4/2, 6/2, 7/2,
        1/3, 2/3, 3/3, 4/3, 5/3, 6/3, 7/3, 8/3, 9/3, 10/3,
        1/4, 2/4, 3/4, 4/4, 5/4, 6/4, 7/4, 8/4, 9/4, 10/4,
        1/5, 2/5, 3/5, 4/5, 5/5, 6/5, 7/5, 8/5, 9/5, 10/5, 11/5,
        1/6, 2/6, 3/6, 4/6, 5/6, 6/6, 7/6, 8/6, 9/6, 10/6, 11/6, 12/6,
        1/7, 2/7, 3/7, 4/7, 5/7, 6/7, 7/7, 8/7, 9/7, 10/7, 11/7, 12/7,
        1/8, 2/8, 3/8, 4/8, 5/8, 6/8, 7/8, 8/8, 9/8, 10/8, 11/8, 12/8,
        2/9, 3/9, 4/9, 5/9, 6/9, 7/9, 8/9, 9/9, 10/9, 11/9, 12/9,
        2/10, 3/10, 4/10, 5/10, 6/10, 7/10, 8/10, 9/10, 10/10, 11/10, 12/10,
        3/11, 4/11, 5/11, 6/11, 7/11, 8/11, 9/11, 10/11, 11/11,
        2/12, 4/12, 5/12, 6/12, 7/12, 8/12, 9/12, 10/12, 11/12,
        4/13, 5/13, 6/13, 7/13, 8/13, 9/13, 10/13
      } {
        \fill [pattern=north east lines, pattern color=gray!30] (\tx,\ty) rectangle ++(1,1);
      }

      \begin{scope}[shift={(1,0)}]
        \foreach \x/\y in {0/0, 0/1, 0/2, 1/1, 2/1}
          \draw[fill=colorT!40] (\x,\y) rectangle ++(1,1);
      \end{scope}

      \begin{scope}[shift={(4,0)}]
        \foreach \x/\y in {0/1, 1/0, 1/1, 1/2, 2/1}
          \draw[fill=colorX!40] (\x,\y) rectangle ++(1,1);
      \end{scope}

      \begin{scope}[shift={(7,1)}]
        \foreach \x/\y in {0/0, 1/0, 1/1, 2/0, 2/1}
          \draw[fill=colorP!40] (\x,\y) rectangle ++(1,1);
      \end{scope}

      \begin{scope}[shift={(0,2)}]
        \foreach \x/\y in {0/0, 0/1, 0/2, 0/3, 0/4}
          \draw[fill=colorI!40] (\x,\y) rectangle ++(1,1);
      \end{scope}

      \begin{scope}[shift={(10,2)}]
        \foreach \x/\y in {0/0, 1/0, 1/1, 1/2, 2/2}
          \draw[fill=colorZ!40] (\x,\y) rectangle ++(1,1);
      \end{scope}

      \begin{scope}[shift={(12,5)}]
        \foreach \x/\y in {0/0, 1/0, 1/1, 1/2, 1/3}
          \draw[fill=colorL!40] (\x,\y) rectangle ++(1,1);
      \end{scope}

      \begin{scope}[shift={(0,7)}]
        \foreach \x/\y in {0/0, 0/1, 0/2, 1/2, 1/3}
          \draw[fill=colorN!40] (\x,\y) rectangle ++(1,1);
      \end{scope}

      \begin{scope}[shift={(12,9)}]
        \foreach \x/\y in {0/2, 1/0, 1/1, 1/2, 2/1}
          \draw[fill=colorF!40] (\x,\y) rectangle ++(1,1);
      \end{scope}

      \begin{scope}[shift={(1,11)}]
        \foreach \x/\y in {0/0, 0/1, 0/2, 1/0, 1/2}
          \draw[fill=colorU!40] (\x,\y) rectangle ++(1,1);
      \end{scope}

      \begin{scope}[shift={(3,12)}]
        \foreach \x/\y in {0/0, 0/1, 0/2, 1/2, 2/2}
          \draw[fill=colorV!40] (\x,\y) rectangle ++(1,1);
      \end{scope}

      \begin{scope}[shift={(10,12)}]
        \foreach \x/\y in {0/2, 1/1, 1/2, 2/0, 2/1}
          \draw[fill=colorW!40] (\x,\y) rectangle ++(1,1);
      \end{scope}

      \begin{scope}[shift={(6,14)}]
        \foreach \x/\y in {0/0, 1/0, 2/0, 2/1, 3/0}
          \draw[fill=colorY!40] (\x,\y) rectangle ++(1,1);
      \end{scope}

      \foreach \sx/\sy/\ex/\ey in {
        11/5/12/5, 12/5/12/6, 2/2/3/2, 2/3/2/2, 6/2/7/2, 6/3/6/2,
        4/2/5/2, 5/2/5/3, 5/3/6/3, 2/9/1/9, 1/9/1/8, 7/14/6/14,
        9/3/10/3, 1/3/2/3, 1/4/1/3, 12/6/13/6, 13/6/13/7, 1/6/1/5,
        13/8/13/9, 5/14/4/14, 4/14/4/13, 9/14/8/14, 1/8/1/7, 7/2/8/2,
        8/2/8/3, 13/11/12/11, 13/10/13/11, 6/14/5/14, 10/3/11/3, 11/3/11/4,
        12/11/12/12, 3/11/2/11, 2/11/2/10, 3/2/4/2, 4/12/3/12, 3/12/3/11,
        1/5/1/4, 2/12/3/12, 3/13/2/13, 2/13/2/12, 3/12/3/13, 13/7/13/8,
        4/13/4/12, 8/3/9/3, 1/7/1/6, 8/14/7/14, 13/9/13/10, 10/14/9/14,
        11/4/11/5, 11/14/10/14, 11/13/11/14, 12/13/11/13, 12/12/12/13, 2/10/2/9
      } {
        \draw[black, line width=0.8pt, ->, shorten >=0.5pt, shorten <=0.5pt] (\ex, \ey) -- (\sx, \sy);
      }
    \end{tikzpicture}
    \caption{The counterclockwise boundary around the largest region is indicated by the arrows. Note that the pentomino fence has area $119$, since the single extra cell in the bottom-left corner is 8-connected to the other cells.}
    \label{fig:edge-boundary}
\end{figure}
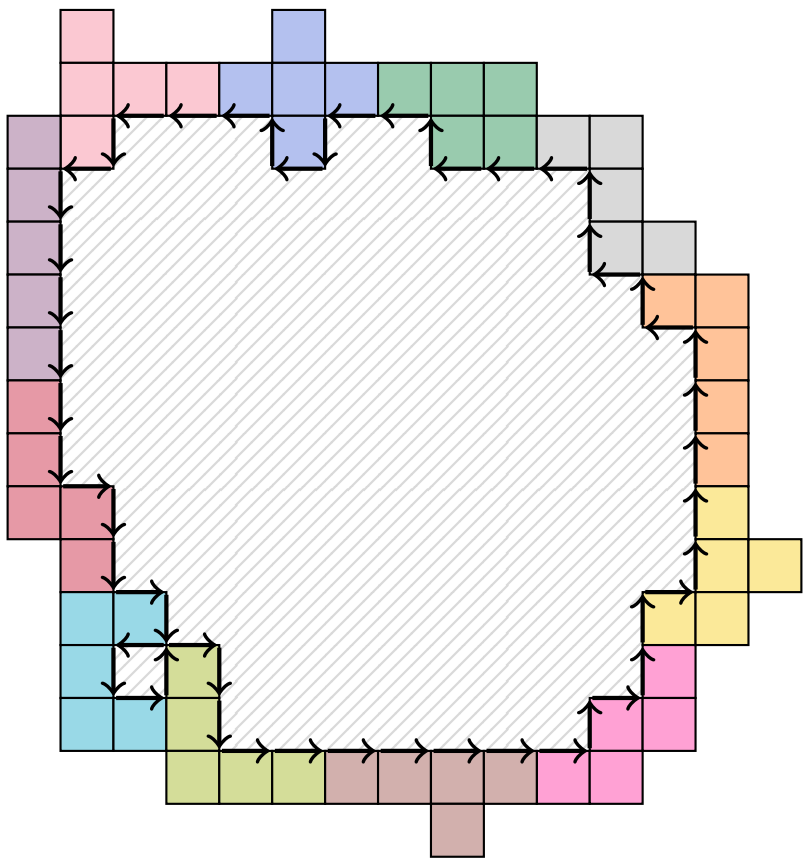

Let $\mathcal R$ be a region of a fence. Its \emph{boundary} consists of
the unit edges between a cell of~$\mathcal R$ and a cell not
in~$\mathcal R$. The latter is always a fence cell, because a free cell
that is adjacent to~$\mathcal R$ would belong to~$\mathcal R$. Every
4-connected component of the complement of~$\mathcal R$ contains such
fence cells, and since $\Phi$ is 4-connected and disjoint
from~$\mathcal R$, it lies in a single 4-connected component. Hence the
complement of~$\mathcal R$ is 4-connected; in other words, no piece can
lie inside~$\mathcal R$. It follows that the boundary of~$\mathcal R$
forms a single closed walk, which we orient counterclockwise
around~$\mathcal R$, so that $\mathcal R$ lies on the left side and the
fence on the right side. At a vertex where two cells of~$\mathcal R$
touch diagonally, the walk passes twice. There, we connect the edges
such that the walk turns right and continues along the same fence cell. Examples are the
corner between \textit{C} and~\textit{L} in \autoref{fig:weird} and the
bottom-left corner in \autoref{fig:edge-boundary}. The area enclosed by
this closed walk, as computed by the Shoelace formula
in \autoref{sec:boundary}, is exactly~$|\mathcal R|$.

Walking along the boundary of~$\mathcal R$, we record the fence cells
that we pass: for each edge, the fence cell on its right side, and at
each left turn, where the walk goes around a convex corner
of~$\mathcal R$, additionally the cell that is diagonally opposite to
the corner cell of~$\mathcal R$. This cell is a fence cell as well,
since it is 8-adjacent to~$\mathcal R$. In this way, we record exactly
the fence cells that are 8-adjacent to~$\mathcal R$. Consecutive recorded
cells are 4-adjacent: After a straight step, the two cells on the right
side are 4-adjacent. After a right turn, the cell on the right side is
the same. At a left turn, the diagonally opposite cell is 4-adjacent to
the cells on the right side before and after the turn. Replacing every
recorded cell by the piece containing it and merging consecutive
repetitions, we obtain a cyclic sequence of pieces, the \emph{pieces
  around~$\mathcal R$}. By construction, consecutive pieces share an
edge.

For the largest region in \autoref{fig:weird}, the pieces around the
region, in counterclockwise order starting at the bottom, are
\textit{L}, \textit{K}, \textit{J}, \textit{E}, \textit{D},
\textit{C}, \textit{O}, \textit{N}, \textit{M}. The piece~\textit{K}
touches the region only at a vertex, and it is recorded only at a left
turn. Nevertheless, it is indispensable: \textit{J} and \textit{L}
touch only at a corner, and without~\textit{K}, the region would leak
to the outside through this diagonal gap. For the region at the right,
the sequence is \textit{I}, \textit{H}, \textit{G}, \textit{F},
\textit{E}, \textit{S}, \textit{E}. The piece~\textit{E} occurs twice,
because \textit{S} is attached to~\textit{E} and protrudes into the
region. Removing~\textit{S} would enlarge this region. The following
lemma shows that this cannot happen for a largest region of an optimal
fence.

\begin{lemma}\label{lem:interval}
  Let $\mathcal R$ be a largest region of a fence of maximum area.
  Then every piece occurs at most once in the cyclic sequence of pieces
  around~$\mathcal R$. Consequently, for every piece~$P$, the edges of
  the boundary of~$\mathcal R$ that belong to~$P$ form a single
  connected (possibly empty) interval along the boundary of~$P$.
\end{lemma}

\begin{proof}
  Suppose that some piece~$P$ occurs twice, so that the cyclic sequence
  has the form $P,\sigma_1,P,\sigma_2$ with nonempty sequences
  $\sigma_1$ and~$\sigma_2$. Let $a$ and~$c$ be cells of~$P$ that were
  recorded at the two occurrences, and let $\eta$ be a 4-connected path
  from $a$ to~$c$ within~$P$. The path~$\eta$ runs outside
  of~$\mathcal R$ and connects two points on the boundary
  of~$\mathcal R$. By the Jordan curve theorem, it cuts off a bounded
  \emph{pocket} from the complement of~$\mathcal R$. The pocket is
  bounded by~$\eta$ and by the part of the boundary walk
  between the two occurrences of~$P$ on one side, say, on the side
  of~$\sigma_1$. Every piece of~$\sigma_1$ lies in the pocket, because
  it is 4-connected and can neither cross~$P$ nor~$\mathcal R$.

  We remove all pieces that lie in the pocket and attach them to the
  outside of the fence, as described in \autoref{sec:def-area}. The
  remaining pieces are still 4-connected, because every 4-connected
  path in~$\Phi$ that passes through the pocket must enter and leave it
  through~$P$. The cells of the pocket become free. They are surrounded
  by~$P$ and~$\mathcal R$, and $\mathcal R$ is still surrounded by~$P$
  and the pieces of~$\sigma_2$. Hence the freed cells are enclosed. Since the recorded cells of~$\sigma_1$ are
  8-adjacent to~$\mathcal R$, the freed cells together with~$\mathcal R$
  belong to a region that is strictly larger than~$\mathcal R$, contradicting
  the optimality of the fence.

  The idea is illustrated in \autoref{fig:edge-boundary-double}: there,
  $P$ is the I-pentomino, $\sigma_1$ consists of the X-pentomino, and the
  X-pentomino lies in the pocket cut off by the two occurrences of~$P$. A piece of the
  pocket may of course also be adjacent to a different enclosed region, as
  in \autoref{fig:edge-boundary-double-special}; but then it occurs in the
  boundary walk of that region and not in that of~$\mathcal R$, so each
  boundary walk still contains it only once.

  For the second statement, consider two consecutive edges of the
  boundary walk that both have a cell of~$P$ on their right side.
  Between them, the walk either goes straight, or it turns right around
  a corner of the same cell, or it turns left around a vertex whose
  diagonally opposite cell belongs to~$P$ as well (since otherwise
  another piece would be recorded between them). In each case, the two
  edges are consecutive on the boundary of~$P$, traversed clockwise
  around~$P$.
\end{proof}

\begin{figure}[h]
    \centering
    \begin{minipage}{0.45\textwidth}
        \centering
        \begin{tikzpicture}[scale=0.35, y=-1cm]

\fill [pattern=north east lines, pattern color=gray!30] (1,2) rectangle ++(4,7);
            
            \begin{scope}[shift={(1,0)}]
                \foreach \x/\y in {0/0, 0/1, 0/2, 1/1, 2/1}
                    \draw[fill=colorT!40] (\x,\y) rectangle ++(1,1);
            \end{scope}

            \begin{scope}[shift={(0,2)}]
                \foreach \x/\y in {0/0, 0/1, 0/2, 0/3, 0/4}
                    \draw[fill=colorI!40] (\x,\y) rectangle ++(1,1);
            \end{scope}

            \begin{scope}[shift={(1,4)}]
                \foreach \x/\y in {0/1, 1/0, 1/1, 1/2, 2/1}
                    \draw[fill=colorX!40] (\x,\y) rectangle ++(1,1);
            \end{scope}

            \begin{scope}[shift={(0,7)}]
                \foreach \x/\y in {0/0, 0/1}
                    \draw[fill=colorL!40] (\x,\y) rectangle ++(1,1);
            \end{scope}

            \foreach \sx/\sy/\ex/\ey in {
                3/2/4/2, 2/2/3/2, 2/3/2/2, 1/3/2/3, 1/4/1/3, 1/5/1/4,
                2/5/1/5, 2/4/2/5, 3/4/2/4, 3/5/3/4, 4/5/3/5, 4/6/4/5,
                3/6/4/6, 3/7/3/6, 2/7/3/7, 2/6/2/7, 1/6/2/6, 1/7/1/6,
                1/8/1/7, 1/9/1/8
            } {
                \draw[black, line width=0.8pt, ->, shorten >=0.5pt, shorten <=0.5pt] (\ex, \ey) -- (\sx, \sy);
            }
        \end{tikzpicture}
        \subcaption{Here the I-pentomino occurs twice in the boundary, but we can simply remove the \textit{pocket} formed by the X-pentomino, which only increases the area.}
        \label{fig:edge-boundary-double}
    \end{minipage}
    \hfill
    \begin{minipage}{0.45\textwidth}
        \centering
        \begin{tikzpicture}[scale=0.35, y=-1cm]
            \foreach \tx/\ty in {
                2/2,
                1/3, 2/3,
                1/4,
                1/6,
                1/7, 2/7
            } {
                \fill [pattern=north east lines, pattern color=gray!30] (\tx,\ty) rectangle ++(1,1);
            }

            \begin{scope}[shift={(1,0)}]
                \foreach \x/\y in {0/0, 0/1, 0/2, 1/1, 2/1}
                    \draw[fill=colorT!40] (\x,\y) rectangle ++(1,1);
            \end{scope}

            \begin{scope}[shift={(0,2)}]
                \foreach \x/\y in {0/0, 0/1, 0/2, 0/3, 0/4}
                    \draw[fill=colorI!40] (\x,\y) rectangle ++(1,1);
            \end{scope}

            \begin{scope}[shift={(1,4)}]
                \foreach \x/\y in {0/1, 1/0, 1/1, 1/2, 2/1}
                    \draw[fill=colorX!40] (\x,\y) rectangle ++(1,1);
            \end{scope}

            \begin{scope}[shift={(0,7)}]
                \foreach \x/\y in {0/0, 0/1}
                    \draw[fill=colorL!40] (\x,\y) rectangle ++(1,1);
            \end{scope}

            \foreach \tx/\ty in {4/2, 4/1, 3/2, 3/3, 3/4} {
                \fill[colorN!40] (\tx,\ty) rectangle ++(1,1);
                \draw (\tx,\ty) rectangle ++(1,1);
            }

            \foreach \tx/\ty in {3/6, 3/7, 3/8, 2/8, 1/8} {
                \fill[colorV!40] (\tx,\ty) rectangle ++(1,1);
                \draw (\tx,\ty) rectangle ++(1,1);
            }

            \foreach \sx/\sy/\ex/\ey in {
                2/2/3/2,
                3/2/3/3,
                3/3/3/4,
                3/4/2/4,
                2/4/2/5,
                2/5/1/5,
                1/5/1/4,
                1/4/1/3,
                1/3/2/3,
                2/3/2/2
            } {
                \draw[black, line width=0.8pt, ->, shorten >=0.5pt, shorten <=0.5pt] (\ex, \ey) -- (\sx, \sy);
            }

            \foreach \sx/\sy/\ex/\ey in {
                1/6/2/6,
                2/6/2/7,
                2/7/3/7,
                3/7/3/8,
                3/8/2/8,
                2/8/1/8,
                1/8/1/7,
                1/7/1/6
            } {
                \draw[black, line width=0.8pt, ->, shorten >=0.5pt, shorten <=0.5pt] (\ex, \ey) -- (\sx, \sy);
            }

        \end{tikzpicture}
        \subcaption{If a piece in the \textit{pocket} is used elsewhere in the fence, we must have two separate regions of enclosed cells and in their separate boundaries each pentomino occurs only once.}
        \label{fig:edge-boundary-double-special}
    \end{minipage}
\end{figure}

Lemma~\ref{lem:interval} justifies the model that we use in the rest
of the paper. The pieces around the largest region~$\mathcal R$ of an
optimal fence form a cyclic sequence $\pi_1,\ldots,\pi_k$ of distinct
pieces, in which consecutive pieces share an edge. Every piece~$\pi_i$
contributes a single consecutive section of its boundary to the
boundary walk of~$\mathcal R$. By the same local analysis as in the
proof, the edge of~$\pi_i$ immediately before this section is shared
with~$\pi_{i-1}$, and the edge immediately after it is shared
with~$\pi_{i+1}$; these will be the \emph{incoming} and \emph{outgoing
  attachment edges} in Section~\ref{sec:pentomino}. The section can be
empty if the piece touches~$\mathcal R$ only at a vertex, like
\textit{K} in \autoref{fig:weird}. The area~$|\mathcal R|$ is obtained by
applying the Shoelace formula to the concatenation of these sections.
Pieces that do not occur around~$\mathcal R$ play no role. A piece may
border several regions, like~\textit{C} in \autoref{fig:weird}, but
only its contact with the largest region matters.

\section{Pentomino fences}\label{sec:pentomino}

\subsection{Edge boundary}\label{sec:boundary}

As already mentioned in Section~\ref{sec:pieces-around}, instead of counting the number of cells, the algorithm will focus on the counterclockwise boundary shown in Figure~\ref{fig:edge-boundary}. So we are given a closed edge boundary with counterclockwise edges
$f_1, \ldots, f_t$, where edge $f_{i}$ starts from coordinates $(x(s_i), y(s_i))$ and ends at $(x(e_i), y(e_i))$. Then the Shoelace formula, first described by Albrecht Meister in 1769~\cite{meister1769}, gives us the area as:
\begin{equation}\label{eq:shoelace}
    A = \frac{1}{2} \left[ \sum_{i = 1}^t x(s_i) y(e_i) - x(e_i) y(s_i) \right].
\end{equation}
Note that since the edges are connected we have $x(e_i) = x(s_{i+1})$ and $y(e_i) = y(s_{i+1})$, with the wrap-around $e_t = s_1$.
  
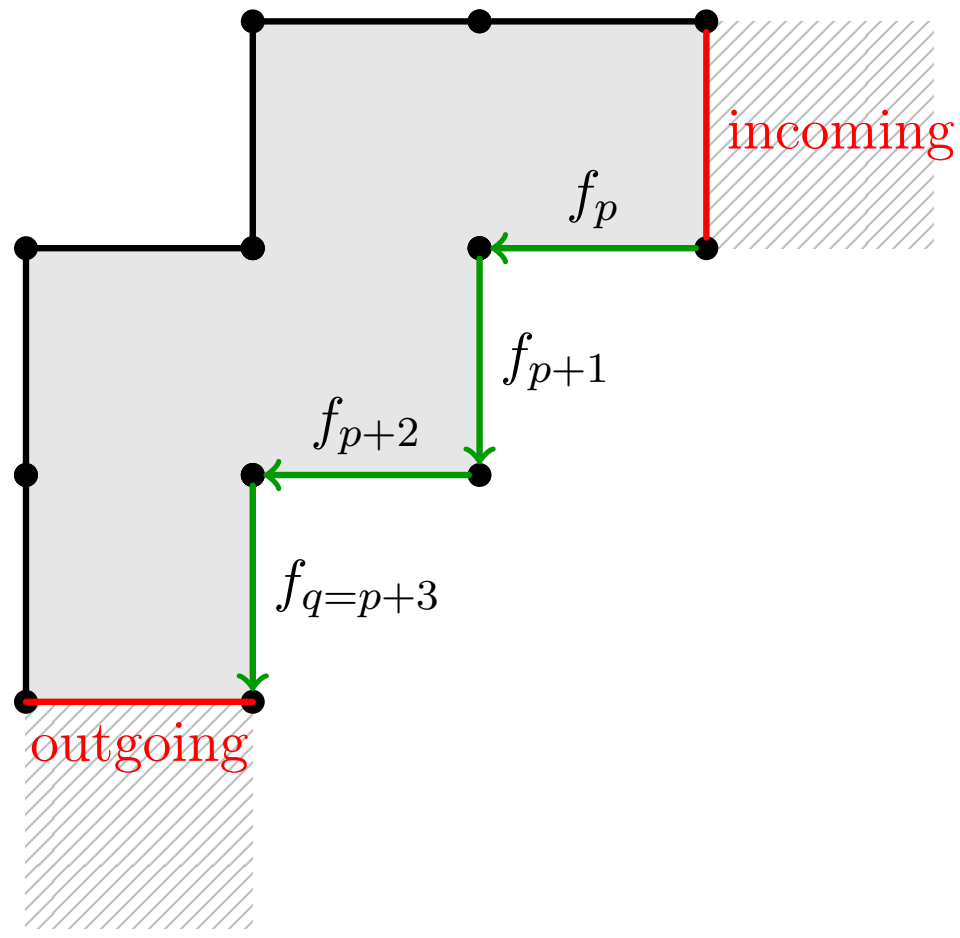
\begin{figure}[h]
    \centering
    \begin{tikzpicture}[scale=1.5, line cap=round, line join=round]
        \foreach \tx/\ty in {0/0, 0/1, 1/1, 1/2, 2/2} {
            \fill[gray!20] (\tx-0.5, \ty+0.5) rectangle ++(1,-1);
        }

        \foreach \tx/\ty in {0/-1, 3/2} {
            \fill[pattern=north east lines, pattern color=gray!50] (\tx-0.5, \ty+0.5) rectangle ++(1,-1);
        }

        \foreach \tx/\ty in {0/0, 0/1, 1/1, 1/2, 2/2} {
            \foreach \dx in {-0.5, 0.5} {
                \foreach \dy in {-0.5, 0.5} {
                    \fill (\tx+\dx, \ty+\dy) circle (1.5pt);
                }
            }
        }

        \draw[
            line width=1.2pt,
            shorten >=2pt,
            shorten <=2pt
        ]
            (-0.5, 0.5) -- (-0.5, 1.5)
            (1.5, 2.5) -- (2.5, 2.5)
            (-0.5, 0.5) -- (-0.5, -0.5)
            (-0.5, 1.5) -- (0.5, 1.5)
            (0.5, 1.5) -- (0.5, 2.5)
            (0.5, 2.5) -- (1.5, 2.5);

        \draw[
            red,
            line width=1.2pt,
            shorten >=2pt,
            shorten <=2pt
        ]
            (0.5, -0.5) -- (-0.5, -0.5)
            (2.5, 2.5) -- (2.5, 1.5);

        \foreach \sx/\sy/\ex/\ey in {2.5/1.5/1.5/1.5, 1.5/1.5/1.5/0.5, 1.5/0.5/0.5/0.5, 0.5/0.5/0.5/-0.5} {
            \draw[green!60!black, line width=1.2pt, ->, shorten >=2pt, shorten <=2pt] (\sx, \sy) -- (\ex, \ey);
        }

        \node[above] at (2.0, 1.5) {$f_p$};
        \node[right] at (1.5, 1.0) {$f_{p+1}$};
        \node[above] at (1.0, 0.5) {$f_{p+2}$};
        \node[right] at (0.5, 0.0) {$f_{q=p+3}$};

        \node[below, red] at (0.0, -0.5) {outgoing};
        \node[right, red] at (2.5, 2.0) {incoming};

    \end{tikzpicture}
    \caption{One possibility of using the W-pentomino. The green
      directed edges form part of the
      boundary of the enclosed region. The red edges are the attachment edges, i.e., the edges
      shared with the neighboring pentominoes, whose cells are
      shaded. They will become important later for ensuring
      edge-to-edge connectivity.}
    \label{fig:pento-possibility}
\end{figure}

To compute the Shoelace formula more efficiently during the search
algorithm, we use a lookup table and precompute all possibilities for how a pentomino can be part of the edge boundary. This means that, for each pentomino, we go over its up to $8$ distinct orientations, obtained by rotation and reflection, and then look at all consecutive edge sections $\{f_{p}, \ldots, f_{q}\}$ of the boundary of the pentomino itself. As the example in Figure~\ref{fig:pento-possibility} shows, these edges are oriented counterclockwise around the edge boundary, but clockwise with respect to the pentomino's own boundary. Then we calculate
\begin{equation*}
  \Delta A_{p, q} := \sum_{i = p}^q x(s_i) y(e_i) - x(e_i) y(s_i),
\end{equation*}
while assuming that the start of the first boundary edge, i.e., $(x(s_{p}), y(s_{p}))$, is shifted to $(0, 0)$. Note that, as seen with piece~\textit{K} in Figure~\ref{fig:weird}, the incoming and outgoing attachment edges (shown in red in Figure~\ref{fig:pento-possibility}) can also be directly adjacent. In that case, we just have $\{f_{p}, \ldots, f_{q}\} = \emptyset$ and $\Delta A_{p, q} = 0$. During the search, we use the following lemma.

\begin{lemma}\label{lem:shoelace-gain}
    When we place the pentomino with the start of the first boundary edge at $(\Delta x, \Delta y)$, the contribution of its edges to the Shoelace formula is
    \begin{equation} \label{eq:shoelace-gain}
      \sum_{i=p}^q \big[x(s_i) y(e_i) - x(e_i) y(s_i)\big] = \Delta A_{p, q} + \Delta x d_y - \Delta y d_x,
    \end{equation}
    where we define $(d_x, d_y)$ as the difference between the end $e_q$ and start vertex $s_p$ of the edge section $\{f_{p}, \ldots, f_{q}\}$.
\end{lemma}

\begin{proof}
  Write $x_0(\cdot), y_0(\cdot)$ for the precomputed coordinates used to obtain $\Delta A_{p,q}$, i.e., those satisfying $(x_0(s_p), y_0(s_p)) = (0,0)$ and $(x_0(e_q), y_0(e_q)) = (d_x, d_y)$. This means

  \begin{equation*}
    x(s_i) = x_0(s_i) + \Delta x, \quad y(s_i) = y_0(s_i) + \Delta y,
  \end{equation*}

  and so for a single Shoelace term we have
  \begin{align*}
    x(s_i) y(e_i) - x(e_i) y(s_i)
    &= \big(x_0(s_i)+\Delta x\big)\big(y_0(e_i)+\Delta y\big) - \big(x_0(e_i)+\Delta x\big)\big(y_0(s_i)+\Delta y\big) \\
    &= \big[x_0(s_i) y_0(e_i) - x_0(e_i) y_0(s_i)\big] + \Delta x \big[y_0(e_i) - y_0(s_i)\big] + \Delta y \big[x_0(s_i) - x_0(e_i)\big].
  \end{align*}

  Summing over $i = p, \ldots, q$, the first bracket sums to $\Delta A_{p,q}$ by definition. For the other two sums, connectivity gives $x_0(e_i) = x_0(s_{i+1})$ and $y_0(e_i) = y_0(s_{i+1})$, so both telescope to
  \begin{equation*}
    \sum_{i=p}^q \big[y_0(e_i) - y_0(s_i)\big] = y_0(e_q) - y_0(s_p), \qquad
    \sum_{i=p}^q \big[x_0(s_i) - x_0(e_i)\big] = x_0(s_p) - x_0(e_q).
  \end{equation*}
  Since $(x_0(s_p), y_0(s_p)) = (0,0)$ and $(x_0(e_q), y_0(e_q)) = (d_x, d_y)$, these reduce to $d_y$ and $-d_x$, respectively, and the formula is proven.
\end{proof}

Before proceeding, we want to discuss a nice geometric interpretation of $\Delta A_{p, q}$. All steps of the dynamic programming search for the pentominoes can be understood without it, but it will be very useful for the hexomino fence.

\begin{figure}[h]
    \centering
\noindent
        \begin{tikzpicture}[scale=1.3, line cap=round, line join=round]
            \foreach \tx/\ty in {0/0, 0/1, 1/1, 1/2, 2/2} {
                \fill[gray!20] (\tx-0.5, \ty+0.5) rectangle ++(1,-1);
            }

            \fill[pattern=dots, pattern color=blue] (3.3,-1.3) -- (2.5,1.5) -- (0.5,-0.5) -- cycle;
            \fill[pattern=dots, pattern color=green!60!black] (2.5,1.5) -- (1.5,1.5) -- (1.5,0.5) -- cycle;
            \fill[pattern=dots, pattern color=green!60!black] (1.5,0.5) -- (0.5,0.5) -- (0.5,-0.5) -- cycle;

            \foreach \tx/\ty in {0/0, 0/1, 1/1, 1/2, 2/2} {
                \foreach \dx in {-0.5, 0.5} {
                    \foreach \dy in {-0.5, 0.5} {
                        \fill (\tx+\dx, \ty+\dy) circle (1.5pt);
                    }
                }
            }

            \draw[line width=1.2pt, shorten >=2pt, shorten <=2pt]
                (-0.5, 0.5) -- (-0.5, 1.5)
                (1.5, 2.5) -- (2.5, 2.5)
                (-0.5, 0.5) -- (-0.5, -0.5)
                (-0.5, 1.5) -- (0.5, 1.5)
                (0.5, 1.5) -- (0.5, 2.5)
                (0.5, 2.5) -- (1.5, 2.5);

            \draw[red, line width=1.2pt, shorten >=2pt, shorten <=2pt]
                (0.5, -0.5) -- (-0.5, -0.5)
                (2.5, 2.5) -- (2.5, 1.5);

            \foreach \sx/\sy/\ex/\ey in {2.5/1.5/1.5/1.5, 1.5/1.5/1.5/0.5, 1.5/0.5/0.5/0.5, 0.5/0.5/0.5/-0.5} {
                \draw[green!60!black, line width=1.2pt, ->, shorten >=2pt, shorten <=2pt] (\sx, \sy) -- (\ex, \ey);
            }

            \node[above] at (2.0, 1.5) {$f_p$};
            \node[left] at (0.5, 0.0) {$f_{q}$};

            \node[right] at (2.6, 1.5) {$(\Delta x, \Delta y)$};
            \node[below left] at (0.5, -0.5) {$(\Delta x + d_x,\ \Delta y + d_y)$};

            \fill (3.3, -1.3) circle (1.5pt);
            \node[below right] at (3.3, -1.1) {$(0,0)$};

            \draw[blue, line width=1.2pt, ->, shorten >=2pt, shorten <=2pt] (3.3, -1.3) -- (2.5, 1.5);
            \draw[blue, line width=1.2pt, ->, shorten >=2pt, shorten <=2pt] (2.5, 1.5) -- (0.5, -0.5);
            \draw[blue, line width=1.2pt, ->, shorten >=2pt, shorten <=2pt] (0.5, -0.5) -- (3.3, -1.3);

            \node[right, blue] at (2.9, 0.1) {$b$};
            \node[right, blue] at (1.0, -0.1) {$c$};
            \node[below, blue] at (1.7, -0.9) {$a$};

            \node at (1.2, -1.5) {(a)};
        \end{tikzpicture}
        \hfill
        \begin{tikzpicture}[scale=1.3, line cap=round, line join=round]
            \foreach \tx/\ty in {0/0, 0/1, 1/1, 0/2, 1/2} {
                \fill[gray!20] (\tx-0.5, \ty+0.5) rectangle ++(1,-1);
            }

            \fill[pattern=dots, pattern color=blue] (3.3,-1.3) -- (1.5,1.5) -- (1.5,0.5) -- (1.0,0.5) -- (0.5,-0.5) -- cycle;
            \fill[pattern=dots, pattern color=red] (1.5,1.5) -- (1.5,0.5) -- (1.0,0.5) -- cycle;
            \fill[pattern=dots, pattern color=green!60!black] (1.0,0.5) -- (0.5,0.5) -- (0.5,-0.5) -- cycle;

            \foreach \vx/\vy in {-0.5/-0.5, -0.5/0.5, -0.5/1.5, -0.5/2.5, 0.5/-0.5, 0.5/0.5, 0.5/2.5, 1.5/0.5, 1.5/1.5, 1.5/2.5} {
                \fill (\vx, \vy) circle (1.5pt);
            }

            \draw[line width=1.2pt, shorten >=2pt, shorten <=2pt]
                (-0.5, -0.5) -- (-0.5, 0.5)
                (-0.5, 0.5) -- (-0.5, 1.5)
                (-0.5, 1.5) -- (-0.5, 2.5)
                (-0.5, 2.5) -- (0.5, 2.5)
                (0.5, 2.5) -- (1.5, 2.5);

            \draw[red, line width=1.2pt, shorten >=2pt, shorten <=2pt]
                (1.5, 2.5) -- (1.5, 1.5)
                (0.5, -0.5) -- (-0.5, -0.5);

            \foreach \sx/\sy/\ex/\ey in {1.5/1.5/1.5/0.5, 1.5/0.5/0.5/0.5, 0.5/0.5/0.5/-0.5} {
                \draw[green!60!black, line width=1.2pt, ->, shorten >=2pt, shorten <=2pt] (\sx, \sy) -- (\ex, \ey);
            }

            \node[right] at (1.6, 1.5) {$(\Delta x, \Delta y)$};
            \node[below left] at (0.5, -0.5) {$(\Delta x + d_x,\ \Delta y + d_y)$};

            \fill (3.3, -1.3) circle (1.5pt);
            \node[below right] at (3.3, -1.1) {$(0,0)$};

            \draw[blue, line width=1.2pt, ->, shorten >=2pt, shorten <=2pt] (3.3, -1.3) -- (1.5, 1.5);
            \draw[blue, line width=1.2pt, ->, shorten >=2pt, shorten <=2pt] (1.5, 1.5) -- (0.5, -0.5);
            \draw[blue, line width=1.2pt, ->, shorten >=2pt, shorten <=2pt] (0.5, -0.5) -- (3.3, -1.3);

            \node[blue] at (2.54, 0.19) {$b$};
            \node[blue] at (0.91, 0.69) {$c$};
            \node[blue] at (1.65, -1.00) {$a$};
            \node at (1.3, -1.5) {(b)};
        \end{tikzpicture}

        \caption{(a)
The difference $\Delta A_{p, q}$ is exactly twice the difference
between the area of the blue triangle and the polygon formed from
edges $a$, $b$ and $\{f_{p}, \ldots, f_{q}\}$. Here, this difference
is positive and drawn as two green triangles.
(b)~The difference between the area of the blue triangle and the
polygon formed from edges $a$, $b$ and $\{f_{p}, \ldots, f_{q}\}$
consists of gained area (green triangle) and lost area (red dotted
triangle).}
        \label{fig:pento-possibility-geometric-P}

\end{figure}
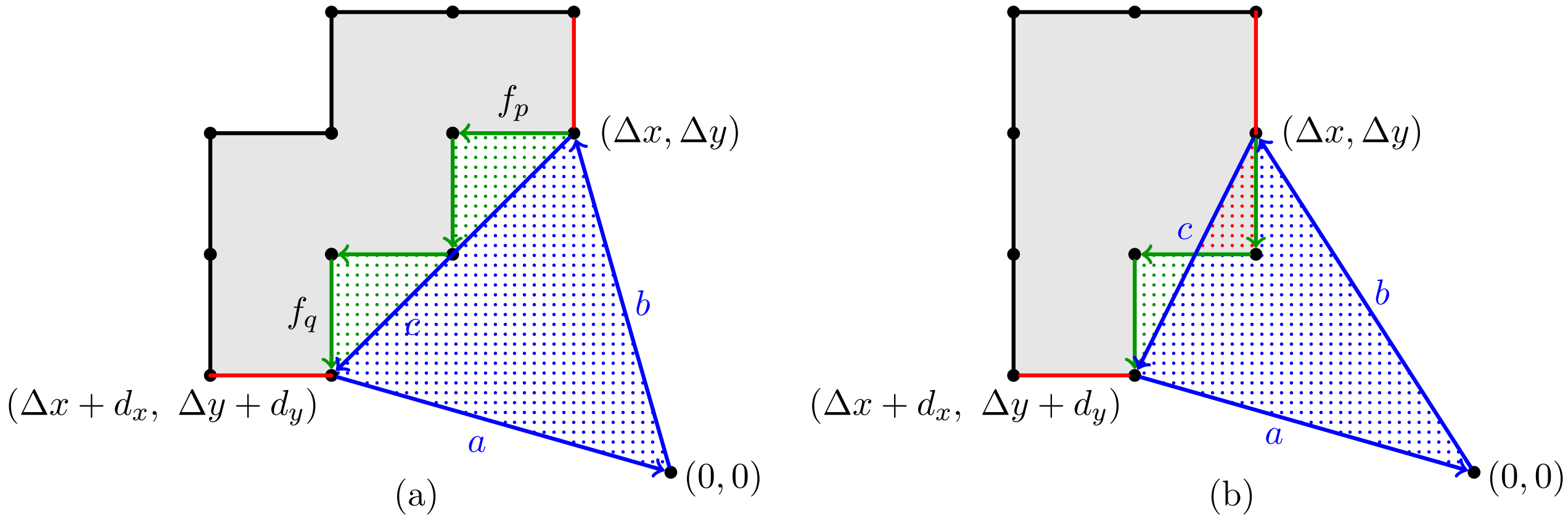

For the consecutive edge section $\{f_{p}, \ldots, f_{q}\}$, we connect its starting vertex $(\Delta x, \Delta y)$ with its ending vertex $(\Delta x + d_x,\ \Delta y + d_y)$ and the origin $(0, 0)$ to form a triangle, as seen in Figure~\ref{fig:pento-possibility-geometric-P}a. By the Shoelace formula, its area is $A = \frac{1}{2} \left[ \sum_{i \in \{a, b, c\}} x(s_i) y(e_i) - x(e_i) y(s_i) \right]$, and the summand corresponding to edge $c$ is:
\begin{equation*}
  x(s_c) y(e_c) - x(e_c) y(s_c) = \Delta x \cdot (\Delta y + d_{y}) - \Delta y \cdot (\Delta x + d_{x}) = \Delta x d_{y} - \Delta y d_{x}.
\end{equation*}

If we instead look at the polygon formed by edges $a$, $b$ and $\{f_{p}, \ldots, f_{q}\}$, the summands for $a$ and $b$ stay the same and by Equation~\eqref{eq:shoelace-gain} the remaining edges contribute $\Delta A_{p, q} + \Delta x d_y - \Delta y d_x$. Thus, the only differing term $\Delta A_{p, q}$ describes exactly what is added when we replace the direct vector $c$ with the polyomino boundary $\{f_{p}, \ldots, f_{q}\}$.

In the example shown in Figure~\ref{fig:pento-possibility-geometric-P}a, this corresponds to the area of the two green triangles. They have total area $1$, and since the Shoelace formula multiplies every contribution by $\frac{1}{2}$, the value of $\Delta A_{p, q}$ is exactly twice that, as a simple calculation confirms.

In some cases, like the P-pentomino in Figure~\ref{fig:pento-possibility-geometric-P}b, we can also lose some area by switching to the polyomino boundary instead of the direct vector (shown with red dots). So we have to subtract the areas that are now lost from the areas gained, which in this specific case comes out to $0$, but can also sometimes be a negative value.

When we later work with such vectors directly in Section~\ref{sec:reference-vectors}, the analogous quantity will be called the \emph{correction} (Definition~\ref{def:reference-vector}); there, the vectors connect the midpoints of the attachment edges instead of the vertices $s_p$ and~$e_q$.

\subsection{Dynamic programming with fixed order}

\subsubsection{Algorithm}

We begin with the restricted case in which the pentominoes are placed in a fixed order $\pi = (\pi_1, \ldots, \pi_{12})$, and discuss the full search over all orders in Section~\ref{sec:dp-unrestricted}.

The key idea is to build the edge boundary one pentomino at a time and track the running Shoelace sum along the way. Crucially, since the Shoelace formula~\eqref{eq:shoelace} is a simple sum over individual edges, we can evaluate it incrementally: each time we append a new pentomino's edge segment to the boundary, we add the corresponding summands directly, without needing the boundary to be closed yet. This allows us to reason about partial boundaries and define an exact objective for each intermediate state of the construction.

We fix the starting grid vertex at the origin $(0,0)$ and place pentominoes one by one, each time extending the open boundary chain by one pentomino's edge segment. After placing the first $i$ pentominoes, the partial boundary is a connected open chain running from $(0,0)$ to some endpoint $(x, y) \in \mathbb{Z}^2$.

\begin{definition}
  $T_i(x, y)$ is the maximum accumulated Shoelace sum over all ways to place $\pi_1, \ldots, \pi_i$ such that the partial boundary ends at $(x,y)$.
\end{definition}

The answer to the fixed-order problem is then $\frac{1}{2} T_{12}(0,0)$: the boundary must return to the origin to close, giving the area of the fence with the $\frac{1}{2}$ coming from the Shoelace formula~\eqref{eq:shoelace}. The base case is $T_0(0,0) = 0$, since the empty boundary starts and ends at the origin with no Shoelace summands.

For the recurrence, recall from Section~\ref{sec:boundary} that we precompute a set of all possibilities of using pentomino $\pi_i$ in one of its $8$ orientations. Each possibility is a triple $(\Delta A, d_x, d_y)$, where $\Delta A$ is the Shoelace contribution of the pentomino's edge segment when its starting vertex is placed at the origin, and $(d_x, d_y)$ is the displacement from that starting vertex to the endpoint of the segment. We denote the set of all possibilities for using pentomino $\pi_i$, over all its orientations, by $\mathrm{Poss}_i$.

To compute $T_i(x, y)$, we consider all ways of reaching $(x, y)$ by appending $\pi_i$ to a partial boundary ending at the predecessor position $(x - d_x,\, y - d_y)$. Applying~\eqref{eq:shoelace-gain} with $(\Delta x, \Delta y) = (x - d_x,\, y - d_y)$ gives the Shoelace gain of this step:
\begin{equation*}
  \Delta A + (x - d_x)\, d_y - (y - d_y)\, d_x.
\end{equation*}
The recurrence is therefore:
\begin{equation} \label{eq:dp-recurrence}
  T_i(x, y) \;=\; \max_{(\Delta A,\, d_x,\, d_y)\,\in\,\mathrm{Poss}_i} \Delta A + (x - d_x)\,d_y - (y - d_y)\,d_x + T_{i-1}(x - d_x,\, y - d_y).
\end{equation}

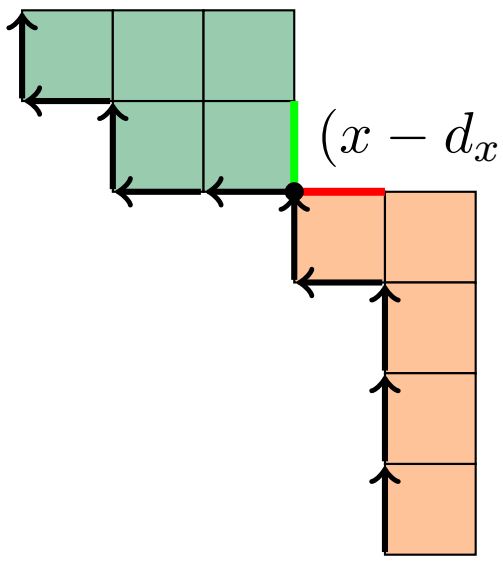
\begin{figure}[h]
    \centering
    \begin{tikzpicture}[scale=0.6, y=-1cm]
        \begin{scope}[shift={(-2,0)}]
            \foreach \x/\y in {0/0, 1/0, 2/0, 1/1, 2/1}
                \draw[fill=colorP!40] (\x,\y) rectangle ++(1,1);
            \foreach \xa/\ya/\xb/\yb in {0/1/0/0, 1/2/1/1, 1/1/0/1, 2/2/1/2, 3/2/2/2}
                \draw[black, line width=1.15pt, ->, shorten >=0.5pt, shorten <=0.5pt] (\xa,\ya) -- (\xb,\yb);
            \node[above left, inner sep=2pt] at (0,0) {$(x,y)$};
        \end{scope}
        \begin{scope}[shift={(0,0)}]
            \foreach \x/\y in {1/2, 2/2, 2/3, 2/4, 2/5}
                \draw[fill=colorL!40] (\x,\y) rectangle ++(1,1);
            \foreach \xa/\ya/\xb/\yb in {1/3/1/2, 2/4/2/3, 2/5/2/4, 2/6/2/5, 2/3/1/3}
                \draw[black, line width=1.15pt, ->, shorten >=0.5pt, shorten <=0.5pt] (\xa,\ya) -- (\xb,\yb);
        \end{scope}
        \draw[red, line width=1.5pt] (1,2) -- (2,2);
        \draw[green, line width=1.5pt] (1,2) -- (1,1);
        \fill[black] (1,2) circle (3pt);
        \node[above right, inner sep=4pt] at (1,2) {$(x - d_x, y - d_y)$};
    \end{tikzpicture}
    \caption{One example demonstrating that the simple recurrence equation does not enforce the edge-to-edge connectivity constraint. The green edge of the P-pentomino does not connect to the red edge of the L-pentomino.}
    \label{fig:orientation}
\end{figure}

\paragraph{Reintroducing the edge-to-edge connectivity constraint.} So far, the recurrence relation ignores the requirement that consecutive pentominoes must share a full grid edge and not merely touch at a corner: As the example in Figure~\ref{fig:orientation} shows, we only check that the boundary edges of the previous pentomino end where the boundary edges of the next one start.

For this we now look at the clockwise edge boundary of the pentomino itself, and define the \emph{incoming attachment edge} as the edge immediately before~$f_{p}$ (the first edge of the section on the inner boundary of the fence) and the \emph{outgoing attachment edge} as the edge immediately after~$f_{q}$ (the last edge of this section). We again refer to Figure~\ref{fig:pento-possibility} for a visualization; in Figure~\ref{fig:orientation}, the outgoing attachment edge of the L-pentomino is drawn in red and the incoming attachment edge of the P-pentomino in green.

To ensure that the pentominoes are connected edge to edge, we simply introduce the additional constraint that the new incoming attachment edge must lie in the same cardinal direction $\tau \in \{N, E, S, W\}$ with respect to the vertex $(x - d_{x}, y - d_{y})$ as the previous outgoing attachment edge.

Accordingly, we augment the DP state to $\psi = ((x, y), \tau)$, where $(x,y)$ is the current endpoint of the partial boundary and $\tau$ is the required incoming attachment edge orientation for the next pentomino. The precomputed possibilities are indexed by the incoming orientation, so $\mathrm{Poss}[\pi_i, \tau]$ contains the tuples $(\tau', \Delta A, d_x, d_y)$ valid when $\pi_i$ is attached from direction $\tau$, and $\tau'$ is the resulting outgoing attachment edge orientation. The initial state is $\psi_0 = ((0,0), \tau_0)$ for a fixed starting orientation $\tau_0$, and a completed fence corresponds to returning to $\psi_0$ after all 12 pentominoes have been placed.

Algorithm \hyperref[alg:optimal-ordered-fence]{\textsc{OptimalOrderedFence}} implements this procedure. The DP table $T_i$ is a dictionary mapping states $\psi = ((x,y),\tau)$ to the best Shoelace sum achieved so far by any partial boundary of length $i$ that reaches $\psi$. At each step $i$, we iterate over all states in $T_{i-1}$ and all valid possibilities for $\pi_i$ given the current orientation. The translation-corrected gain $A_{\mathrm{add}}$ is computed via~\eqref{eq:shoelace-gain} using the predecessor position (note that now $(x, y)$ refers to the predecessor coordinates, i.e., $(d_x, d_y)$ is already subtracted), and the new state $\psi'$ records the updated position and outgoing orientation. On the final step $i = 12$, only transitions that return to $\psi_0$ are permitted, enforcing closure of the boundary. The algorithm tracks the best closed-boundary Shoelace sum (twice the area) in $\hat A_{\max}$ and returns $A_{\max} := \frac{1}{2} \hat A_{\max}$ as the optimal area.

\begin{algorithm}[H]
  \caption{\textsc{OptimalOrderedFence}$(\pi)$}
  \label{alg:optimal-ordered-fence}
  \begin{algorithmic}[1]
    \State \textbf{Input:} order $\pi = (\pi_1, \ldots, \pi_{12})$, precomputed possibilities $\mathrm{Poss}$
    \State $\psi_0 \gets ((0,0), \tau_0)$ \Comment{fixed start position and orientation}
    \State $T_0[\psi_0] \gets 0$ \Comment{initialize the dynamic program}
    \State $\hat A_{\max} \gets 0$
    \For{$i = 1$ to $12$}
      \State $T_i \gets \emptyset$
      \For{each state $((x,y), \tau)$ with value $A$ in $T_{i-1}$}
        \For{each $(\tau', \Delta A, (d_x,d_y)) \in \mathrm{Poss}[\pi_i, \tau]$}
          \State $A_{\mathrm{add}} \gets \Delta A + x d_y - y d_x$ \Comment{translation correction in the Shoelace sum}
          \State $\psi' \gets ((x+d_x, y+d_y), \tau')$
          \If{$i = 12$ and $\psi' \neq \psi_0$}
            \State \textbf{continue} \Comment{enforce boundary closure on the final step}
          \EndIf
          \State $T_i[\psi'] \gets \max\{T_i[\psi'], A + A_{\mathrm{add}}\}$
          \If{$i = 12$ and $\psi' = \psi_0$} \label{algline:amax-check} \Comment{record the best closed boundary found so far}
            \State $\hat A_{\max} \gets \max\{\hat A_{\max}, A + A_{\mathrm{add}}\}$
          \EndIf
        \EndFor
      \EndFor \label{algline:finish-iter}
    \EndFor
    \State \Return $\frac{1}{2} \hat A_{\max}$
  \end{algorithmic}
\end{algorithm}

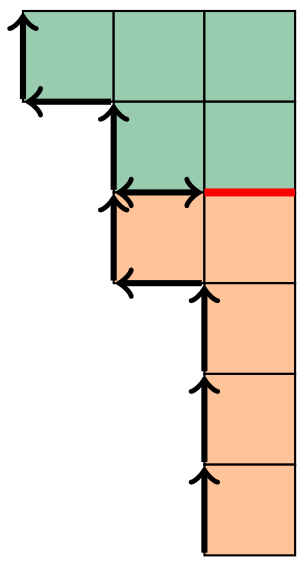
\begin{figure}[h]
    \centering
    \begin{tikzpicture}[scale=0.6, y=-1cm]
      \begin{scope}[shift={(0,0)}]
        \foreach \x/\y in {0/0, 1/0, 2/0, 1/1, 2/1}
          \draw[fill=colorP!40] (\x,\y) rectangle ++(1,1);
      \end{scope}

      \begin{scope}[shift={(0,0)}]
        \foreach \x/\y in {1/2, 2/2, 2/3, 2/4, 2/5}
          \draw[fill=colorL!40] (\x,\y) rectangle ++(1,1);
      \end{scope}

      \draw[red, line width=1.5pt] (2,2) -- (3,2);

      \foreach \xa/\ya/\xb/\yb in {0/1/0/0, 1/2/1/1, 1/3/1/2, 2/4/2/3, 2/5/2/4, 2/6/2/5, 2/3/1/3, 1/1/0/1}
        \draw[black, line width=1.15pt, ->, shorten >=0.5pt, shorten <=0.5pt] (\xa,\ya) -- (\xb,\yb);

      \draw[black, line width=1.15pt, <->, shorten >=0.5pt, shorten <=0.5pt] (1,2) -- (2,2);
    \end{tikzpicture}
    \caption{If consecutive pentominoes share several edges, then, depending on the choice of the attachment edge (drawn in red), the inner boundary can contain appendices (here the double-arrow edge). They cancel out in the Shoelace formula and contribute zero area, so they can be safely ignored.}
    \label{fig:shoelace-multiple}
\end{figure}

\paragraph{Multiple shared edges.} Note that our search method does not prevent consecutive pentominoes from sharing more than one edge, as shown in Figure~\ref{fig:shoelace-multiple}. Depending on which edge is chosen as the respective incoming/outgoing attachment edge, the remaining shared edges appear as back-and-forth appendices on the boundary. However, these appendices contribute zero to the Shoelace sum: each such edge is traversed once in each direction, so its two contributions cancel exactly.

\paragraph{Solution reconstruction and self-intersection.}
Each DP table entry stores only the current state $\psi = ((x,y), \tau)$, discarding all information about which orientations and attachment edges were used for the pentominoes placed so far. This compression is what makes the search tractable, but it means \hyperref[alg:optimal-ordered-fence]{\textsc{OptimalOrderedFence}} cannot directly output a concrete pentomino placement, nor can it detect whether any candidate placement is \emph{self-intersecting}, i.e.,\ whether two pentominoes overlap on the same grid cell.

\begin{figure}[h]
    \centering
    \begin{tikzpicture}[scale=0.4, y=-1cm]
      \foreach \tx/\ty in {
        5/2, 6/2, 7/2, 8/2, 9/2, 10/2,
        2/3, 3/3, 4/3, 5/3, 6/3, 7/3, 8/3, 9/3, 10/3, 11/3,
        2/4, 3/4, 4/4, 5/4, 6/4, 7/4, 8/4, 9/4, 10/4, 11/4,
        1/5, 2/5, 3/5, 4/5, 5/5, 6/5, 7/5, 8/5, 9/5, 10/5, 11/5, 12/5,
        1/6, 2/6, 3/6, 4/6, 5/6, 6/6, 7/6, 8/6, 9/6, 10/6, 11/6, 12/6,
        1/7, 2/7, 3/7, 4/7, 5/7, 6/7, 7/7, 8/7, 9/7, 10/7, 11/7, 12/7,
        1/8, 2/8, 3/8, 4/8, 5/8, 6/8, 7/8, 8/8, 9/8, 12/8,
        1/9, 2/9, 3/9, 4/9, 5/9, 6/9, 7/9, 8/9, 9/9, 10/9,
        1/10, 2/10, 3/10, 4/10, 5/10, 6/10, 7/10, 8/10, 9/10, 10/10, 12/10,
        1/11, 2/11, 3/11, 4/11, 5/11, 6/11, 7/11, 8/11, 9/11, 10/11, 11/11, 12/11,
        2/12, 3/12, 4/12, 5/12, 6/12, 7/12, 8/12, 9/12, 10/12, 11/12, 12/12,
        3/13, 4/13, 5/13, 6/13, 7/13, 8/13, 9/13, 6/14
      } {
        \fill [pattern=north east lines, pattern color=gray!30] (\tx,\ty) rectangle ++(1,1);
      }

      \foreach \tx/\ty in {10/8, 11/8, 11/9, 11/10} {
        \fill [pattern=dots, pattern color=gray!60] (\tx,\ty) rectangle ++(1,1);
      }

      \foreach \sx/\sy/\ex/\ey in {
        6/2/5/2, 7/2/6/2, 8/2/7/2, 9/2/8/2, 10/2/9/2, 11/2/10/2,
        5/2/5/3, 11/3/11/2, 3/3/2/3, 2/3/2/4, 4/3/3/3, 5/3/4/3,
        12/3/11/3, 12/4/12/3, 2/4/2/5, 12/5/12/4, 2/5/1/5, 1/5/1/6,
        13/5/12/5, 13/6/13/5, 1/6/1/7, 13/7/13/6, 1/7/1/8, 13/8/13/7,
        1/8/1/9, 12/9/13/9, 13/9/13/8, 1/9/1/10, 12/10/12/9,
        12/9/12/8, 12/8/11/8, 11/8/10/8, 10/8/10/9, 10/9/11/9,
        11/9/11/10, 11/10/11/11, 11/11/12/11, 12/11/12/10,
        1/10/1/11, 13/10/12/10, 13/11/13/10, 1/12/2/12, 1/11/1/12, 13/12/13/11,
        2/13/3/13, 2/12/2/13, 10/13/11/13, 11/13/12/13, 12/13/13/13,
        13/13/13/12, 3/14/4/14, 3/13/3/14, 4/14/5/14, 5/14/6/14,
        7/14/8/14, 8/14/9/14, 9/14/10/14, 10/14/10/13, 6/15/7/15,
        6/14/6/15, 7/15/7/14
      } {
        \draw[black, line width=0.8pt, ->, shorten >=0.5pt, shorten <=0.5pt] (\sx, \sy) -- (\ex, \ey);
      }

      \draw[black, line width=0.8pt, ->>, shorten >=0.5pt, shorten <=0.5pt] (12, 10) -- (12, 9);
    \end{tikzpicture}
    \caption{Example of a self-intersecting boundary that can arise when placed pentominoes overlap. The doubly-covered cells (dotted region) cause the boundary to cross itself; the Shoelace formula counts the area of this region twice via its winding number.}
    \label{fig:pento-self-inter}
\end{figure}
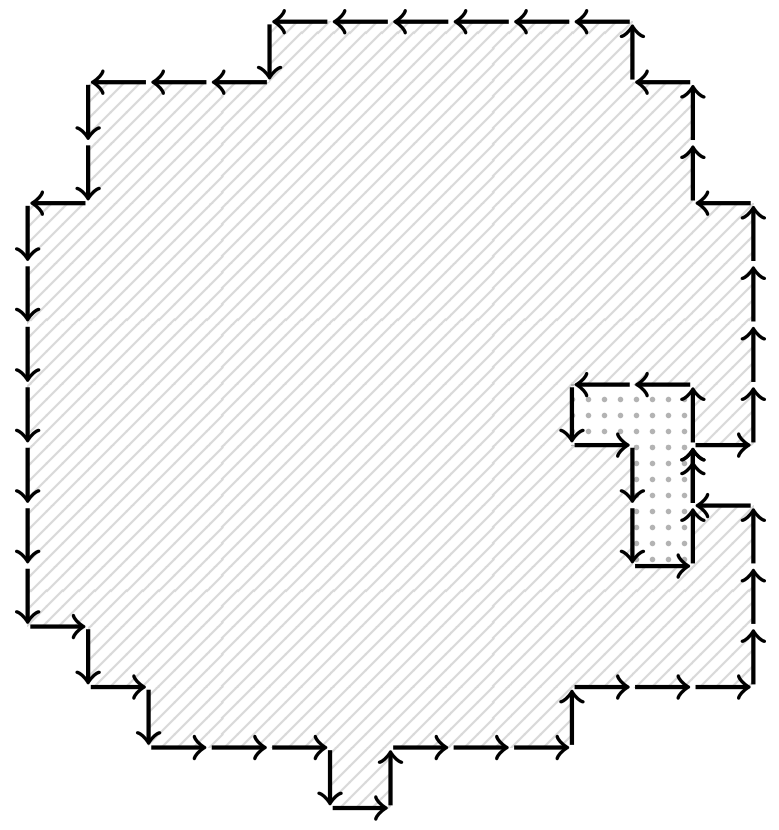

Self-intersection is a genuine concern because the Shoelace formula remains well-defined even for self-crossing boundaries: it computes area weighted by the winding number, so any doubly-covered region is counted twice. An example is shown in Figure~\ref{fig:pento-self-inter}, where overlapping pentominoes cause the boundary to cross itself and the Shoelace value exceeds the true enclosed area. This means the update
\[
  T_i[\psi'] \gets \max\{T_i[\psi'],\, A + A_{\mathrm{add}}\}
\]
could in principle promote a value originating from an invalid self-intersecting placement, making $\frac{1}{2}\hat A_{\max}$ technically an upper bound rather than a confirmed fence area. Additionally, the enclosed cells could consist of multiple regions and by our definition we would only count the largest of them as the fence's area.

The absence of a concrete solution and the self-intersection ambiguity can both be resolved by solution reconstruction. We maintain a second table $\mathrm{Pred}$ alongside $T$: whenever a transition with $A + A_{\mathrm{add}} > T_i[\psi']$ causes an update, we record $\mathrm{Pred}_i[\psi'] \gets (\psi, p)$, where $\psi$ is the predecessor state and $p$ is the possibility used for $\pi_i$. After the search, the optimal placement is recovered by tracing back through $\mathrm{Pred}$: starting from $\mathrm{Pred}_{12}[\psi_0]$ gives the state $\psi_{11}$ and how $\pi_{12}$ was placed; then $\mathrm{Pred}_{11}[\psi_{11}]$ gives $\psi_{10}$ and how $\pi_{11}$ was placed; and so on. With the explicit placement in hand, we can check directly whether any two pentominoes overlap. If no overlap is found, which was always the case in practice, the self-intersection concern is eliminated. After checking that all enclosed cells belong to the same region with a breadth-first search (BFS), we could always confirm $\frac{1}{2} \hat A_{\max}$ as the true optimal fence area for the order $\pi$.

To enumerate \emph{all} optimal solutions, one instead collects all predecessor pairs $(\psi, p)$ with $A + A_{\mathrm{add}} = T_i[\psi']$ and runs standard backtracking through $\mathrm{Pred}$, verifying each reconstructed solution individually. It should be noted that when adjacent pentominoes share multiple edges, two different paths could correspond to the same pentomino fence. Although this never occurred in practice for optimal fences, one has to verify, strictly speaking, that all the produced solutions are actually distinct.

\subsubsection{Solving the min-max problem}\label{sec:min-max}

The \emph{min-max problem} asks: over all possible orders $\pi$ of the 12 pentominoes, what is the minimum value of the maximum fence area achievable under that order? We solve this by running the algorithm \hyperref[alg:optimal-ordered-fence]{\textsc{OptimalOrderedFence}} exhaustively over all orders.

The number of distinct orders to consider is $\frac{11!}{2} = 19\,958\,400$, since by symmetry we can fix the first piece and the direction in which the pentominoes are placed. One nice optimization that speeds up the search is prefix reuse: if two consecutive orders in the enumeration share the same first $t$ pentominoes, the DP tables $T_1, \ldots, T_t$ and $\mathrm{Pred}_1, \ldots, \mathrm{Pred}_t$ are identical and can be carried over without recomputation. To maximize the applicability of this optimization, we process the orders in lexicographic order.

We ran the full search multi-threaded on a machine with an AMD Ryzen 7
7735HS CPU and 20\,GB of
DDR5 %
main
memory, which took around 2 days to complete.

\begin{figure}[h]
\centering
\begin{tikzpicture}
\begin{axis}[
            width=11cm,
            height=9cm,
            ybar,
            bar width=0.7,
            ymin=0,
            xtick=data,
            xticklabel style={font=\tiny, rotate=45, anchor=east},
            ylabel style={yshift=0.5em},
            scaled y ticks=base 10:-6,
            ytick scale label code/.code={},
            yticklabel={\pgfmathprintnumber{\tick}M},
            ymajorgrids=true,
            grid style=dashed,
            nodes near coords,
            nodes near coords style={font=\tiny, rotate=90, anchor=west, /pgf/number format/fixed, /pgf/number format/1000 sep={\,}},
            enlarge y limits={upper, value=0.25},
]
\addplot+[fill=blue, draw=none] coordinates {
              (111, 313)       (112, 16056)     (113, 163064)
              (114, 737289)    (115, 1693225)   (116, 2482741)
              (117, 2872532)   (118, 2795947)   (119, 2585031)
              (120, 2183052)   (121, 1677223)   (122, 1249939)
              (123, 803864)    (124, 399534)    (125, 193506)
              (126, 84228)     (127, 19416)     (128, 1440)
          };
\end{axis}
\end{tikzpicture}
\caption{Distribution of the maximum achievable fence area across all $\frac{11!}{2}$ pentomino orders. The minimum value of $111$ is achieved by $313$ orders.}
\label{fig:pent-ord}
\end{figure}
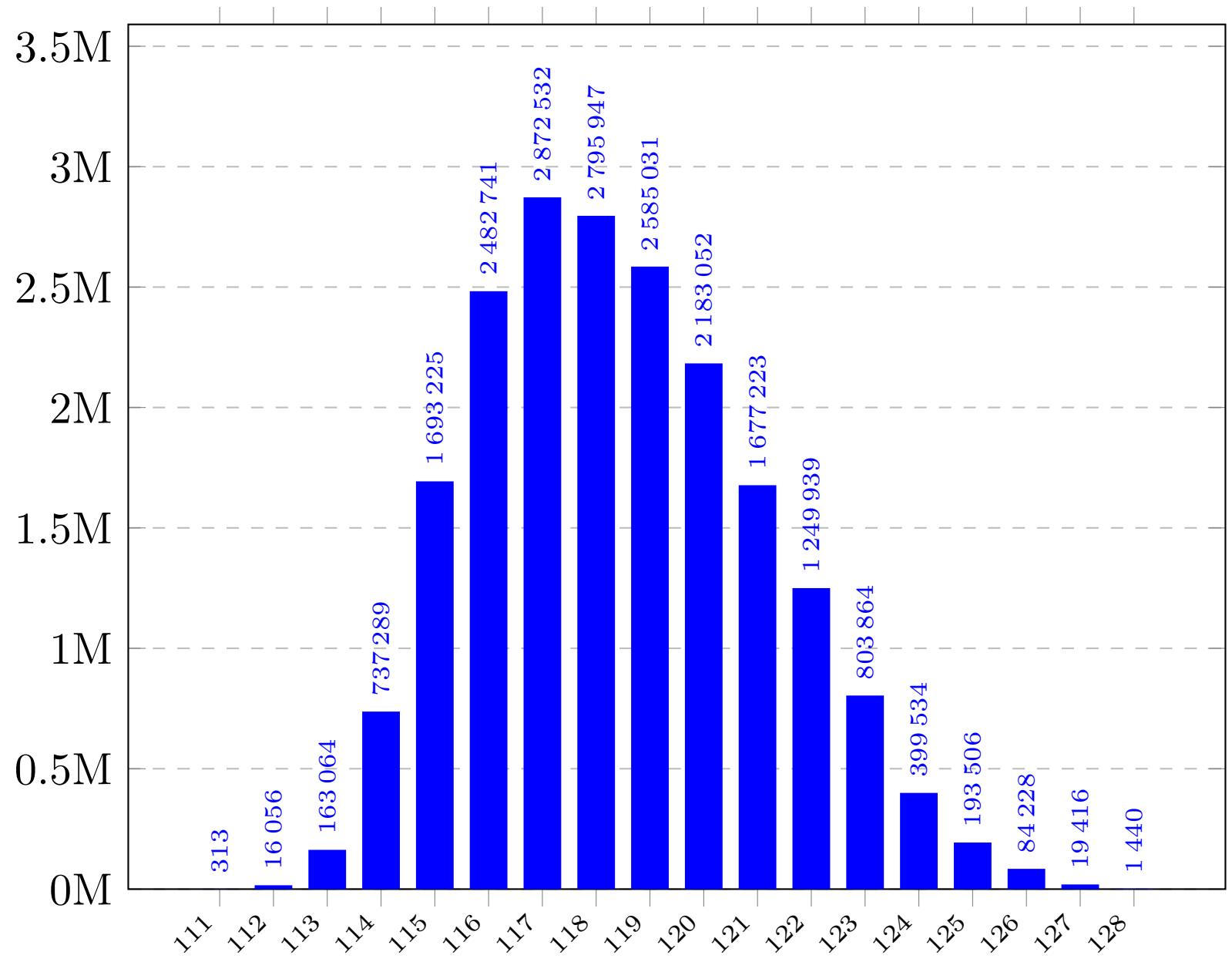

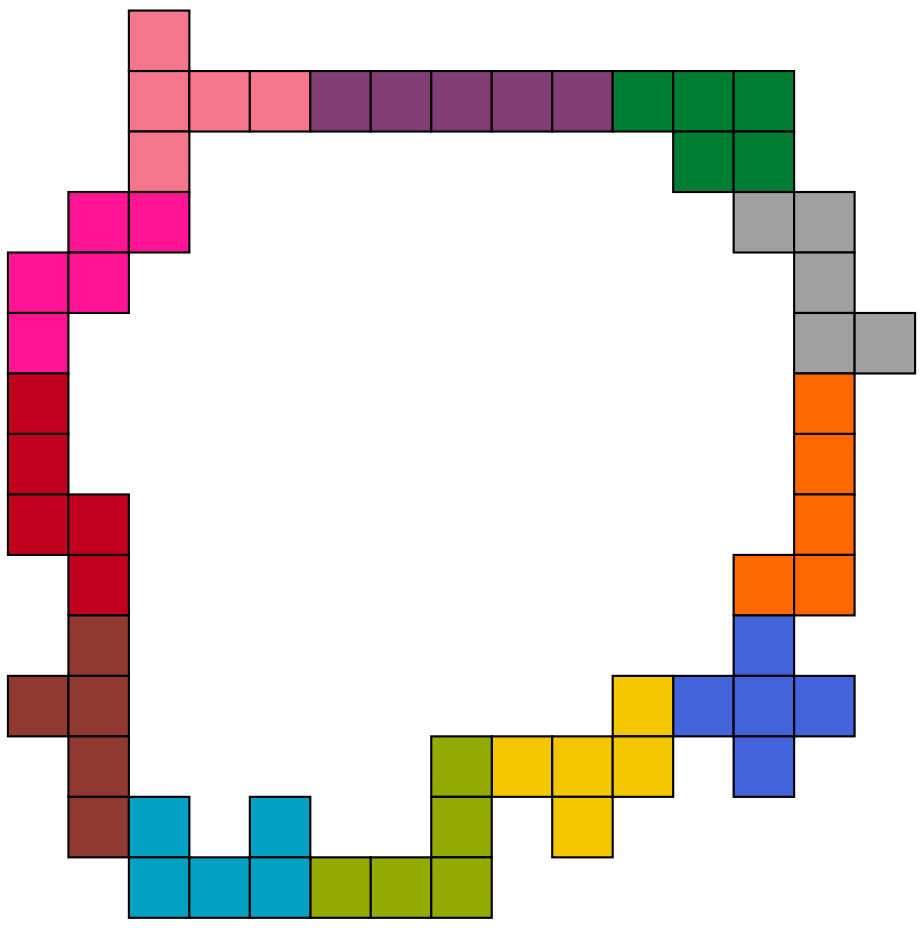
\begin{figure}[h]
  \centering
  \begin{tikzpicture}[scale=0.4, y=-1cm]
    \begin{scope}[shift={(2,0)}]
      \foreach \x/\y in {0/0, 0/1, 0/2, 1/1, 2/1}
        \draw[fill=colorT] (\x,\y) rectangle ++(1,1);
    \end{scope}
    \begin{scope}[shift={(5,1)}]
      \foreach \x/\y in {0/0, 1/0, 2/0, 3/0, 4/0}
        \draw[fill=colorI] (\x,\y) rectangle ++(1,1);
    \end{scope}
    \begin{scope}[shift={(10,1)}]
      \foreach \x/\y in {0/0, 1/0, 2/0, 1/1, 2/1}
        \draw[fill=colorP] (\x,\y) rectangle ++(1,1);
    \end{scope}
    \begin{scope}[shift={(0,3)}]
      \foreach \x/\y in {1/0, 2/0, 0/1, 1/1, 0/2}
        \draw[fill=colorW] (\x,\y) rectangle ++(1,1);
    \end{scope}
    \begin{scope}[shift={(12,3)}]
      \foreach \x/\y in {0/0, 1/0, 1/1, 1/2, 2/2}
        \draw[fill=colorZ] (\x,\y) rectangle ++(1,1);
    \end{scope}
    \begin{scope}[shift={(0,6)}]
      \foreach \x/\y in {0/0, 0/1, 0/2, 1/2, 1/3}
        \draw[fill=colorN] (\x,\y) rectangle ++(1,1);
    \end{scope}
    \begin{scope}[shift={(12,6)}]
      \foreach \x/\y in {1/0, 1/1, 1/2, 0/3, 1/3}
        \draw[fill=colorL] (\x,\y) rectangle ++(1,1);
    \end{scope}
    \begin{scope}[shift={(0,10)}]
      \foreach \x/\y in {1/0, 0/1, 1/1, 1/2, 1/3}
        \draw[fill=colorY] (\x,\y) rectangle ++(1,1);
    \end{scope}
    \begin{scope}[shift={(11,10)}]
      \foreach \x/\y in {1/0, 0/1, 1/1, 2/1, 1/2}
        \draw[fill=colorX] (\x,\y) rectangle ++(1,1);
    \end{scope}
    \begin{scope}[shift={(8,11)}]
      \foreach \x/\y in {2/0, 0/1, 1/1, 2/1, 1/2}
        \draw[fill=colorF] (\x,\y) rectangle ++(1,1);
    \end{scope}
    \begin{scope}[shift={(5,12)}]
      \foreach \x/\y in {2/0, 2/1, 0/2, 1/2, 2/2}
        \draw[fill=colorV] (\x,\y) rectangle ++(1,1);
    \end{scope}
    \begin{scope}[shift={(2,13)}]
      \foreach \x/\y in {0/0, 2/0, 0/1, 1/1, 2/1}
        \draw[fill=colorU] (\x,\y) rectangle ++(1,1);
    \end{scope}
  \end{tikzpicture}
  \caption{One of the $313$ min-max pentomino fence orders, enclosing a maximum area of $111$ grid cells.}
  \label{fig:minmax-fence}
\end{figure}

The distribution of optimal areas across all orders is shown in Figure~\ref{fig:pent-ord}. The minimum of the distribution, and so the answer to the min-max problem, is $111$, achieved by $313$ orders. One such min-max solution is shown in Figure~\ref{fig:minmax-fence}. For hexominoes, we have only sampled random orders (Appendix~\ref{sec:hex-dist}).

\paragraph{Interactive lab.} In the \emph{pentomino DP lab} (\url{https://erikaroldanroa.github.io/fence-challenge/pentomino-dp-lab/}) the reader can follow the dynamic program step by step for any order.

\subsection{Dynamic programming with unrestricted order} \label{sec:dp-unrestricted}

To solve the original pentomino fence problem, i.e., finding the largest area achievable by any valid placement, regardless of the order, we extend \hyperref[alg:optimal-ordered-fence]{\textsc{OptimalOrderedFence}} with one additional loop. Rather than advancing to a fixed next pentomino $\pi_i$ at each step, we iterate over all pentominoes not yet used in the current partial boundary and try each one as the next piece. The DP state is augmented to record the set of used pentominoes alongside the position and orientation, and the recurrence otherwise remains identical. Thus, we can keep everything else in Algorithm~\ref{alg:optimal-ordered-fence}.

Crucially, the state only needs to track \emph{which} pentominoes have been used, not the order in which they were placed. This is valid because the Shoelace sum accumulated so far is already fully captured by the value stored in the table: two partial boundaries that used the same set of pentominoes and arrive at the same state $\psi = ((x,y), \tau)$ with the same accumulated sum are entirely interchangeable from the perspective of any future extension, regardless of the order in which those pentominoes were placed.

\paragraph{State space size.}

We can easily bound the size of each table $T_{i}$ by the total number of possible state attribute combinations. There are always just $4$ possibilities for the orientation, and $\sbinom{12}{i}$ many ways of choosing the $i$ used pentominoes. For the current coordinates $(x, y)$ we can have all grid points reachable so far, which can be roughly estimated as $50i^2$, the approximate size of a Manhattan ball of radius $5i$ (since each pentomino roughly contributes $5$ cells in Manhattan distance). In Table~\ref{tab:dp-sizes}, we illustrate that this estimate (third column) matches the actual number of reachable positions during the algorithm (fourth column) well.

Then the upper bound for $|T_i|$ is obtained by multiplying the $4$ orientations with the $\sbinom{12}{i}$ used pentomino sets and finally the number of reachable positions. This means that the fifth column of Table~\ref{tab:dp-sizes} is the product of the second and fourth columns times $4$. It shows that no $|T_i|$ ever contains more than a few million entries,
which is also confirmed by the actual observed $|T_i|$ in the last column. The large discrepancy between actual size and upper bound for small $i$ is caused by the fact that some pentomino subsets can reach far fewer coordinates.

\begin{table}[h]
  \centering
  \begin{tabular}{|r||r|r|>{\ \ }r<{\hskip 2em}|>{\hskip 2em}r<{\hskip 2em}|>{\ \ }r<{\ \ }|}
    \hline
    $i$ & $\binom{12\mathstrut}{i\mathstrut}$
    & $50i^{2\mathstrut}$ & \multicolumn1{c|}{actual $|\{(x, y)\}|$}
    & \multicolumn1{c|}{upper bound $|T_i|$}
    & \multicolumn1{c|}{actual $|T_i|$} \\
    \hline
     1 &    12 &     50 &     57 &      2\,736 &        730 \\
     2 &    66 &    200 &    212 &     55\,968 &     29\,519 \\
     3 &   220 &    450 &    456 &    401\,280 &    245\,919 \\
     4 &   495 &    800 &    784 &  1\,552\,320 &  1\,028\,126 \\
     5 &   792 &  1\,250 &  1\,172 &  3\,712\,896 &  2\,621\,789 \\
     6 &   924 &  1\,800 &  1\,628 &  6\,017\,088 &  4\,451\,545 \\
     7 &   792 &  2\,450 &  2\,092 &  6\,627\,456 &  5\,230\,563 \\
     8 &   495 &  3\,200 &  2\,604 &  5\,155\,920 &  4\,292\,220 \\
     9 &   220 &  4\,050 &  3\,172 &  2\,791\,360 &  2\,424\,160 \\
    10 &    66 &  5\,000 &  3\,796 &  1\,002\,144 &    900\,744 \\
    11 &    12 &  6\,050 &  4\,476 &    214\,848 &    198\,688 \\
    \hline
  \end{tabular}
  \caption{The upper bound and actual size of the DP table $T_i$ at each step.}
  \label{tab:dp-sizes}
\end{table}

Note that the actual $|T_i|$ values reported in the last column are measured with the following symmetry breaking: After numbering the pentominoes, we require that the second placed pentomino has the lowest id of any pentomino used (to eliminate rotational symmetry) and that the first pentomino has a lower id than the third pentomino (to eliminate mirror symmetry).

This all means that except for the final check that the fence is closed at $i = 12$, we do not need any pruning. If instead we kept track of the order of pentominoes, the $\sbinom{12}{i}$ would turn into an $i!$ factor and make $|T_i|$ far too large to keep in memory.

\paragraph{Solution reconstruction.}
Solution reconstruction carries over unchanged in structure. The back-pointer table $\mathrm{Pred}$ now additionally records which pentomino was chosen at each step, so that the traceback recovers both the sequence of pentominoes and their placements. The self-intersection argument requires slightly more care: in the fixed-order case, an inflated Shoelace value from an overlapping placement can only displace a correct value for the same order. Here, it could in principle alter which pentomino appears to have been placed at which step, reordering the reconstructed sequence. However, the same verification argument holds: if no overlaps are detected in the reconstructed placement, then the placement is geometrically valid regardless of order, and if BFS shows that all enclosed cells belong to the same region, $\frac{1}{2} \hat A_{\max}$ is confirmed as the true optimum. If overlaps or multiple regions are found, $\frac{1}{2}\hat A_{\max}$ remains an upper bound.

Running the search single-threaded on the same Ryzen 7 machine mentioned in Section~\ref{sec:min-max}, it only takes 2.5 minutes to complete. It should be mentioned that when updating $\hat A_{\max}$ with unrestricted orders we can optionally drop the condition $i = 12$ (line~\ref{algline:amax-check}). Then each time we are finished building a given $T_{i}$ and arrive at line~\ref{algline:finish-iter}, $\frac{1}{2} \hat A_{\max}$ corresponds to the largest area of a fence achievable by using any $i$ of the $12$ pentominoes. We present those intermediate values in Table~\ref{tab:max-subset}. As already discussed at the end of Section~\ref{sec:def-area}, we observe that it is never advantageous to leave out a pentomino.

\begin{table}[htb]
  \centering
  \begin{tabular}{|c||c|c|c|c|c|c|c|c|c|c|c|c|}
    \hline
    size $i$
      & 1 & 2 & 3 & 4 & 5 & 6 & 7 & 8 & 9 & 10 & 11 & 12 \\
    \hline
    max.\ fence area $A_{\max}(i)$
      & 0 & 1 & 6 & 13 & 24 & 34 & 47 & 62 & 78 & 94 & 111 & 128 \\
    \hline
     $A_{\max}(i)-A_{\max}(i-1)$
      & {--} & 1 & 5 & 7 & 11 & 10 & 13 & 15 & 16 & 16 & 17 & 17 \\
    \hline
  \end{tabular}
  \caption{Largest enclosed area that can be obtained with $i$ pentominoes.}
  \label{tab:max-subset}
\end{table}

\subsection{Identification of all optimal pentomino fences}\label{sec:all-5}
As discussed before, solution reconstruction can also be run in the all-solutions mode, enumerating every optimal placement. This yields exactly $5760$ different solutions, which we can quickly verify to be valid and distinct. We call two solutions distinct if no combination of translation, rotation, and reflection maps one onto the other with every pentomino covering the same cells.
Figure~\ref{fig:optimal-fence} shows
one optimal solution.

\begin{figure}[h]
  \centering
  \begin{minipage}[t]{0.48\textwidth}
    \centering
    \begin{tikzpicture}[scale=0.4, y=-1cm]
      \begin{scope}[shift={(5,0)}]
        \foreach \x/\y in {1/0, 0/1, 1/1, 2/1, 1/2}
          \draw[fill=colorX] (\x,\y) rectangle ++(1,1);
      \end{scope}
      \begin{scope}[shift={(2,1)}]
        \foreach \x/\y in {0/0, 1/0, 2/0, 0/1, 0/2}
          \draw[fill=colorV] (\x,\y) rectangle ++(1,1);
      \end{scope}
      \begin{scope}[shift={(8,1)}]
        \foreach \x/\y in {0/0, 1/0, 2/0, 2/1, 3/1}
          \draw[fill=colorN] (\x,\y) rectangle ++(1,1);
      \end{scope}
      \begin{scope}[shift={(12,2)}]
        \foreach \x/\y in {0/0, 1/0, 1/1, 1/2, 2/2}
          \draw[fill=colorZ] (\x,\y) rectangle ++(1,1);
          \draw[dashed] (0.2,2.8) -- (2.8,0.2);
    \end{scope}
      \begin{scope}[shift={(0,4)}]
        \foreach \x/\y in {0/0, 1/0, 2/0, 1/1, 1/2}
          \draw[fill=colorT] (\x,\y) rectangle ++(1,1);
      \end{scope}
      \begin{scope}[shift={(14,5)}]
        \foreach \x/\y in {0/0, 0/1, 1/1, 0/2, 0/3}
        \draw[fill=colorY] (\x,\y) rectangle ++(1,1);
        \draw[dashed] (-1,2) -- (3,2);
      \end{scope}
      \begin{scope}[shift={(0,7)}]
        \foreach \x/\y in {0/0, 1/0, 0/1, 0/2, 1/2}
          \draw[fill=colorU] (\x,\y) rectangle ++(1,1);
      \end{scope}
      \begin{scope}[shift={(13,9)}]
        \foreach \x/\y in {0/2, 1/0, 1/1, 1/2, 2/1}
          \draw[fill=colorF] (\x,\y) rectangle ++(1,1);
      \end{scope}
      \begin{scope}[shift={(1,10)}]
        \foreach \x/\y in {0/0, 0/1, 1/1, 0/2, 1/2}
          \draw[fill=colorP] (\x,\y) rectangle ++(1,1);
      \end{scope}
      \begin{scope}[shift={(11,12)}]
        \foreach \x/\y in {2/0, 1/1, 2/1, 0/2, 1/2}
          \draw[fill=colorW] (\x,\y) rectangle ++(1,1);
      \end{scope}
      \begin{scope}[shift={(2,13)}]
        \foreach \x/\y in {0/0, 0/1, 1/1, 2/1, 3/1}
          \draw[fill=colorL] (\x,\y) rectangle ++(1,1);
      \end{scope}
      \begin{scope}[shift={(6,14)}]
        \foreach \x/\y in {0/0, 1/0, 2/0, 3/0, 4/0}
          \draw[fill=colorI] (\x,\y) rectangle ++(1,1);
      \end{scope}
    \end{tikzpicture}
    \subcaption{One of the $5760$ optimal pentomino fence solutions,
    enclosing an area of $128$ grid cells. If we reflect the
    gray Z-pentomino at its diagonal and/or
    the brown Y-pentomino at its horizontal axis, we get $3$ further optimal
    solutions with the same order of pentominoes.}
    \label{fig:optimal-fence}
  \end{minipage}%
  \hfill
  \begin{minipage}[t]{0.48\textwidth}
    \centering
    \begin{tikzpicture}[scale=0.4, y=-1cm]
      \colorlet{colorF}{colorF!35}
      \colorlet{colorI}{colorI!35}
      \colorlet{colorL}{colorL!35}
      \colorlet{colorN}{colorN!35}
      \colorlet{colorP}{colorP!35}
      \colorlet{colorT}{colorT!35}
      \colorlet{colorU}{colorU!35}
      \colorlet{colorV}{colorV!35}
      \colorlet{colorW}{colorW!35}
      \colorlet{colorX}{colorX!35}
      \colorlet{colorY}{colorY!35}
      \colorlet{colorZ}{colorZ!35}

      \begin{scope}[shift={(5,0)}]   
        \foreach \x/\y in {1/0, 0/1, 1/1, 2/1, 1/2}
          \draw[fill=colorX] (\x,\y) rectangle ++(1,1);
      \end{scope}
      \begin{scope}[shift={(2,1)}]   
        \foreach \x/\y in {0/0, 1/0, 2/0, 0/1, 0/2}
          \draw[fill=colorV] (\x,\y) rectangle ++(1,1);
      \end{scope}
      \begin{scope}[shift={(8,1)}]   
        \foreach \x/\y in {0/0, 1/0, 2/0, 2/1, 3/1}
          \draw[fill=colorN] (\x,\y) rectangle ++(1,1);
      \end{scope}
      \begin{scope}[shift={(12,2)}]  
        \foreach \x/\y in {0/0, 1/0, 1/1, 1/2, 2/2}
          \draw[fill=colorZ] (\x,\y) rectangle ++(1,1);
      \end{scope}
      \begin{scope}[shift={(0,4)}]   
        \foreach \x/\y in {0/0, 1/0, 2/0, 1/1, 1/2}
          \draw[fill=colorT] (\x,\y) rectangle ++(1,1);
      \end{scope}
      \begin{scope}[shift={(14,5)}]  
        \foreach \x/\y in {0/0, 0/1, 1/1, 0/2, 0/3}
          \draw[fill=colorY] (\x,\y) rectangle ++(1,1);
      \end{scope}
      \begin{scope}[shift={(0,7)}]   
        \foreach \x/\y in {0/0, 1/0, 0/1, 0/2, 1/2}
          \draw[fill=colorU] (\x,\y) rectangle ++(1,1);
      \end{scope}
      \begin{scope}[shift={(13,9)}]  
        \foreach \x/\y in {0/2, 1/0, 1/1, 1/2, 2/1}
          \draw[fill=colorF] (\x,\y) rectangle ++(1,1);
      \end{scope}
      \begin{scope}[shift={(1,10)}]  
        \foreach \x/\y in {0/0, 0/1, 1/1, 0/2, 1/2}
          \draw[fill=colorP] (\x,\y) rectangle ++(1,1);
      \end{scope}
      \begin{scope}[shift={(11,12)}] 
        \foreach \x/\y in {2/0, 1/1, 2/1, 0/2, 1/2}
          \draw[fill=colorW] (\x,\y) rectangle ++(1,1);
      \end{scope}
      \begin{scope}[shift={(2,13)}]  
        \foreach \x/\y in {0/0, 0/1, 1/1, 2/1, 3/1}
          \draw[fill=colorL] (\x,\y) rectangle ++(1,1);
      \end{scope}
      \begin{scope}[shift={(6,14)}]  
        \foreach \x/\y in {0/0, 1/0, 2/0, 3/0, 4/0}
          \draw[fill=colorI] (\x,\y) rectangle ++(1,1);
      \end{scope}

      \begin{scope}[red, line width=1.6pt, line cap=round]
        \draw (5,1)   -- (5,2);     
        \draw (8,1)   -- (8,2);     
        \draw (12,2)  -- (12,3);    
        \draw (14,5)  -- (15,5);    
        \draw (14,9)  -- (15,9);    
        \draw (13,12) -- (14,12);   
        \draw (11,14) -- (11,15);   
        \draw (6,14)  -- (6,15);    
        \draw (2,13)  -- (3,13);    
        \draw (1,10)  -- (2,10);    
        \draw (1,7)   -- (2,7);     
        \draw (2,4)   -- (3,4);     
      \end{scope}

      \coordinate (mXV) at (5,1.5);
      \coordinate (mXN) at (8,1.5);
      \coordinate (mNZ) at (12,2.5);
      \coordinate (mZY) at (14.5,5);
      \coordinate (mYF) at (14.5,9);
      \coordinate (mFW) at (13.5,12);
      \coordinate (mWI) at (11,14.5);
      \coordinate (mIL) at (6,14.5);
      \coordinate (mLP) at (2.5,13);
      \coordinate (mPU) at (1.5,10);
      \coordinate (mUT) at (1.5,7);
      \coordinate (mTV) at (2.5,4);

      \begin{scope}[blue, thick, -stealth]
        \draw (mXN) -- (mXV);
        \draw (mXV) -- (mTV);
        \draw (mTV) -- (mUT);
        \draw (mUT) -- (mPU);
        \draw (mPU) -- (mLP);
        \draw (mLP) -- (mIL);
        \draw (mIL) -- (mWI);
        \draw (mWI) -- (mFW);
        \draw (mFW) -- (mYF);
        \draw (mYF) -- (mZY);
        \draw (mZY) -- (mNZ);
        \draw (mNZ) -- (mXN);
      \end{scope}

      \foreach \m in {XV,XN,NZ,ZY,YF,FW,WI,IL,LP,PU,UT,TV}
        \fill[blue] (m\m) circle (0.1);
    \end{tikzpicture}
    \subcaption{The attachment edges (red) and the reference vectors (blue),
      which connect the red edges between neighboring pentominoes in
      counterclockwise direction.}
    \label{fig:optimal-fence-refpoly}
  \end{minipage}
\end{figure}

\begin{figure}[htb]
  \centering
  \includegraphics[page=2,scale=1]{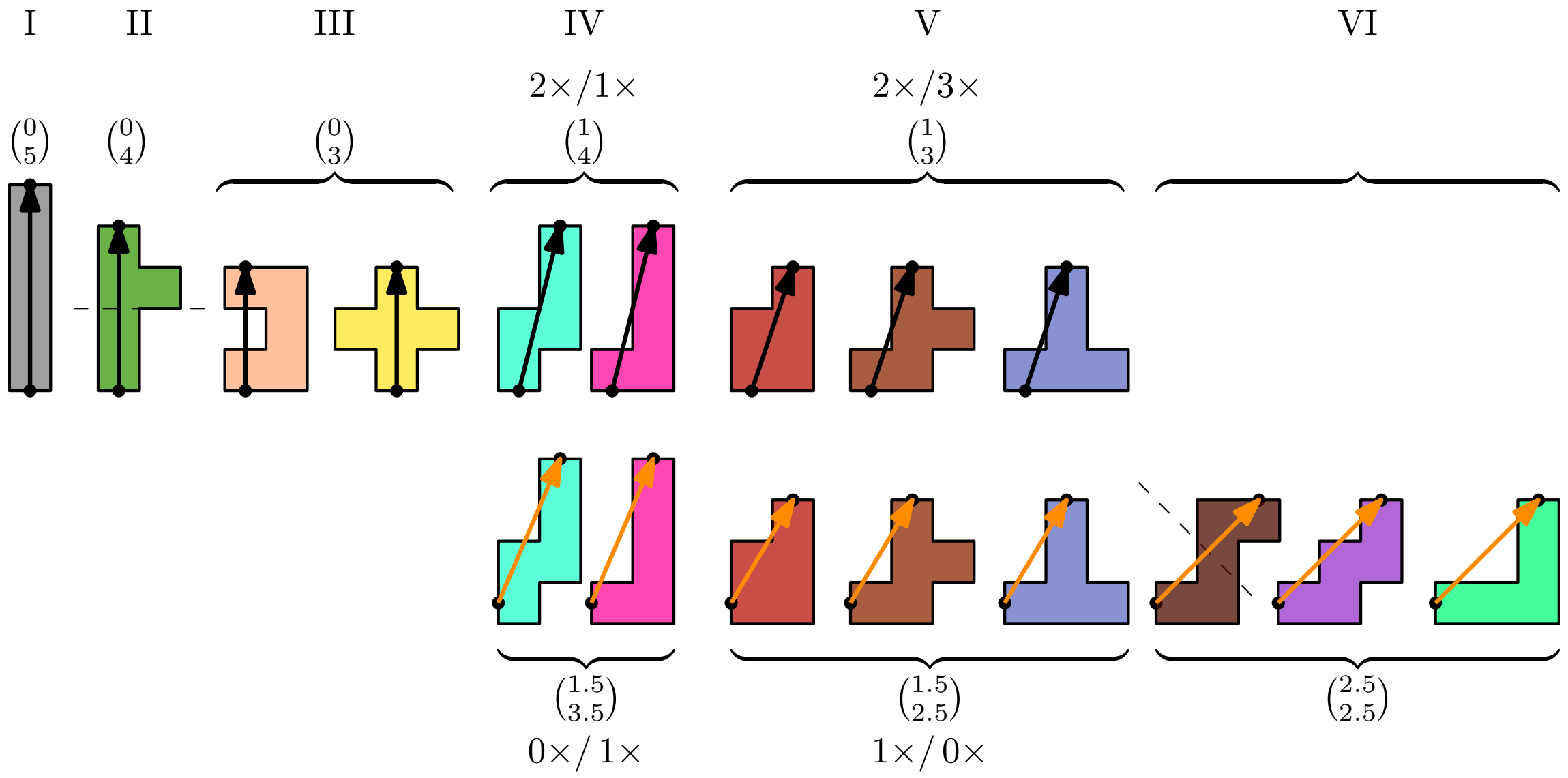}
  \caption{Usage of the parts in optimal pentomino fences.
The lower row shows the ways in which a pentomino can serve as one of
the four corner pieces, and the upper row shows the straight pieces.
  }
  \label{fig:use-opt}
\end{figure}

\begin{figure}[htb]
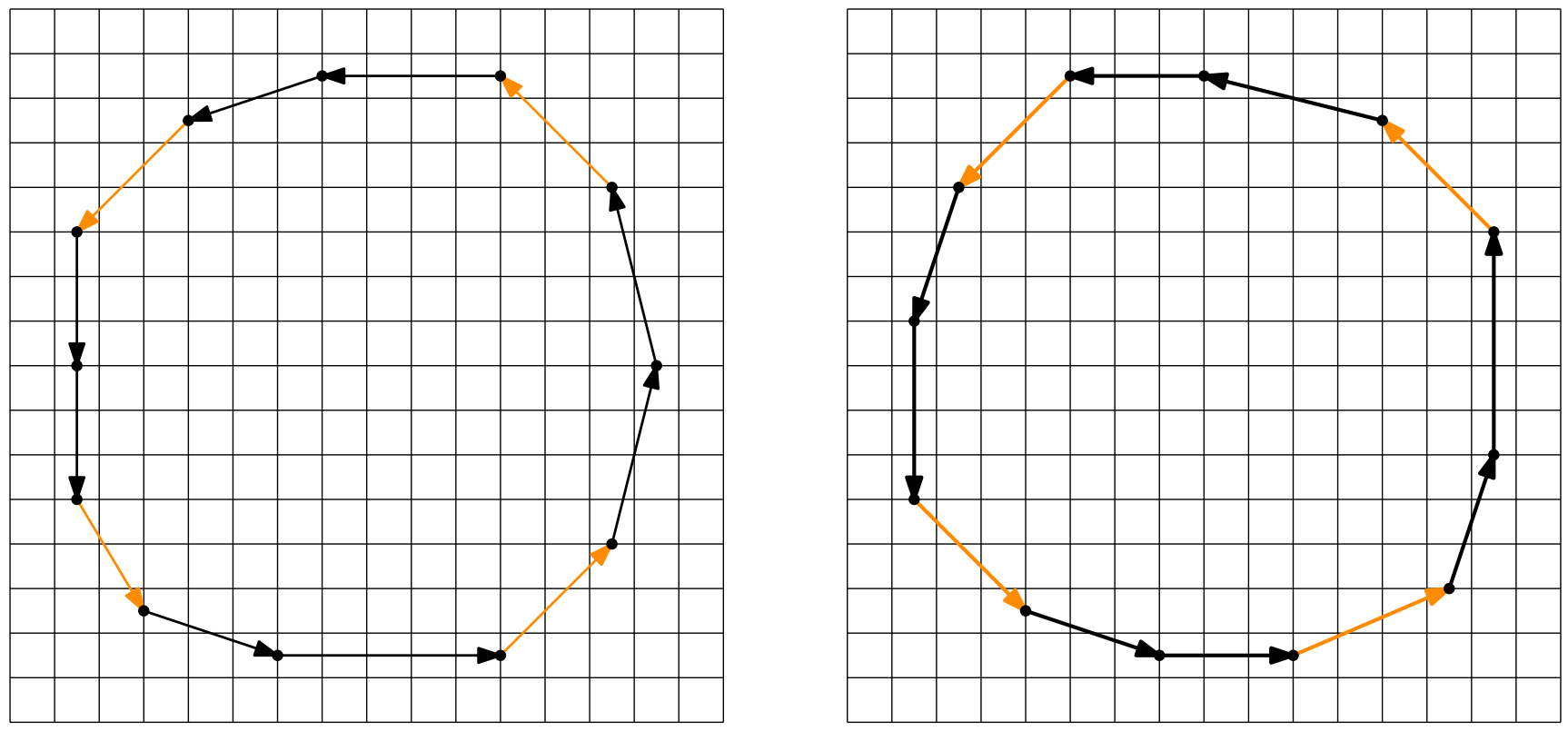

  \centering
  \noindent
  \includegraphics[page=7,scale=0.65]{best-fence-pentominos-with-colors}
  \qquad
  \includegraphics[page=3,scale=0.65]{best-fence-pentominos-with-colors}
  \caption{Two optimal reference polygons out of ten.
The four corner vectors are highlighted in orange.
  }
    \label{fig:refpol-opt}
\end{figure}

Despite this large number, it is possible to get a systematic overview
of the solutions.
For this purpose, it is useful to abstract from the actual shapes of
the pieces and consider only their \emph{reference vectors}, which run
between the midpoints of consecutive attachment edges, see \autoref{fig:optimal-fence-refpoly}.
(This
concept will play a prominent role when we
consider hexominoes in Section~\ref{sec:hexomino}, see Definition~\ref{def:reference-vector}.)
The reference vectors of a fence, taken together, form the
\emph{reference polygon.}

It turns out that every optimal fence has exactly four \emph{corner
  pieces}, where the attachment edges change direction between
horizontal and vertical. The remaining pieces are used in a
``straight'' way. This property is not surprising; it is also true for all optimal hexomino
fences, and most likely also for $n$-omino fences with larger~$n$.

\autoref{fig:use-opt} shows how the twelve pentominoes are used in an optimal
fence.
They 
can be classified by their reference vectors into groups I--VI of sizes
$1,1,2,2,3,3$.
Each pentomino of groups IV and V can be used in two possible ways,
either as a corner piece or as a straight piece.
All pieces of the same group are exchangeable.

The three pieces of group VI, with reference vector
$\binom{2.5}{2.5}$,
can only be used as corner pieces.
The fourth corner piece is either of the form
 $\binom{1.5}{2.5}$ from group~V,
or
$\binom{1.5}{3.5}$ from group~IV.
(When we say ``of the form
$\binom{1.5}{3.5}$'', we mean any of the eight rotated or
reflected vectors
$\binom{\pm1.5}{\pm3.5}$
or $\binom{\pm3.5}{\pm1.5}$.)

If we just look at the reference polygons,
it turns out that there are just ten different reference polygons that
 occur in optimal
solutions, two of which are shown in
\autoref{fig:refpol-opt}.
They are of two types:
\begin{enumerate}
\item
  There are six reference polygons that contain
two reference vectors of the form
$\binom{1}{4}$
and
two reference vectors of the form
$\binom{1}{3}$,
see \autoref{fig:refpol-opt} (left) for an example.
The fourth corner vector is of the form
 $\binom{1.5}{2.5}$.

\item
  There are four reference polygons that contain
one reference vector of the form
$\binom{1}{4}$
and
three reference vectors of the form
$\binom{1}{3}$,
see \autoref{fig:refpol-opt} (right) for an example.
The fourth corner vector is of the form
 $\binom{1.5}{3.5}$.
\end{enumerate}

For a fixed reference polygon, we have to
assign the pieces to the reference vectors.
The pieces of each group can be arbitrarily permuted among the
positions for which they fit.
This leads to
$2!\cdot 2!\cdot 3!\cdot 3!=144$ possibilities, for a total of
$10\times 144=1440$ different optimal pentomino orders,
which is in accordance with the data reported in
Figure~\ref{fig:pent-ord}.

In addition, the Z-pentomino can always be reflected in its
diagonal, and the Y-pentomino along its shorter axis, without changing
the enclosed area, as indicated by
the dotted axis lines in
\autoref{fig:use-opt}.
This leads to the total of $1440\times 2\times 2= 5760$ optimal fences.

\section{Hexomino fences}\label{sec:hexomino}

For the $35$ hexominoes, the dynamic program of Section~\ref{sec:dp-unrestricted} is no longer feasible: tracking the set of used pieces alone introduces a factor of up to $\binom{35}{17} \approx 4.5 \cdot 10^{9}$, and the number of reachable positions grows as well. In this section we therefore exploit the geometric structure of good fences to reduce the search space drastically. The main result is the following.

\begin{theorem}\label{thm:1597}
  The largest area that can be enclosed by a fence of the $35$ hexominoes is $1597$.
\end{theorem}

\begin{figure}[htb]
  \centering
  \includegraphics[scale=0.22]{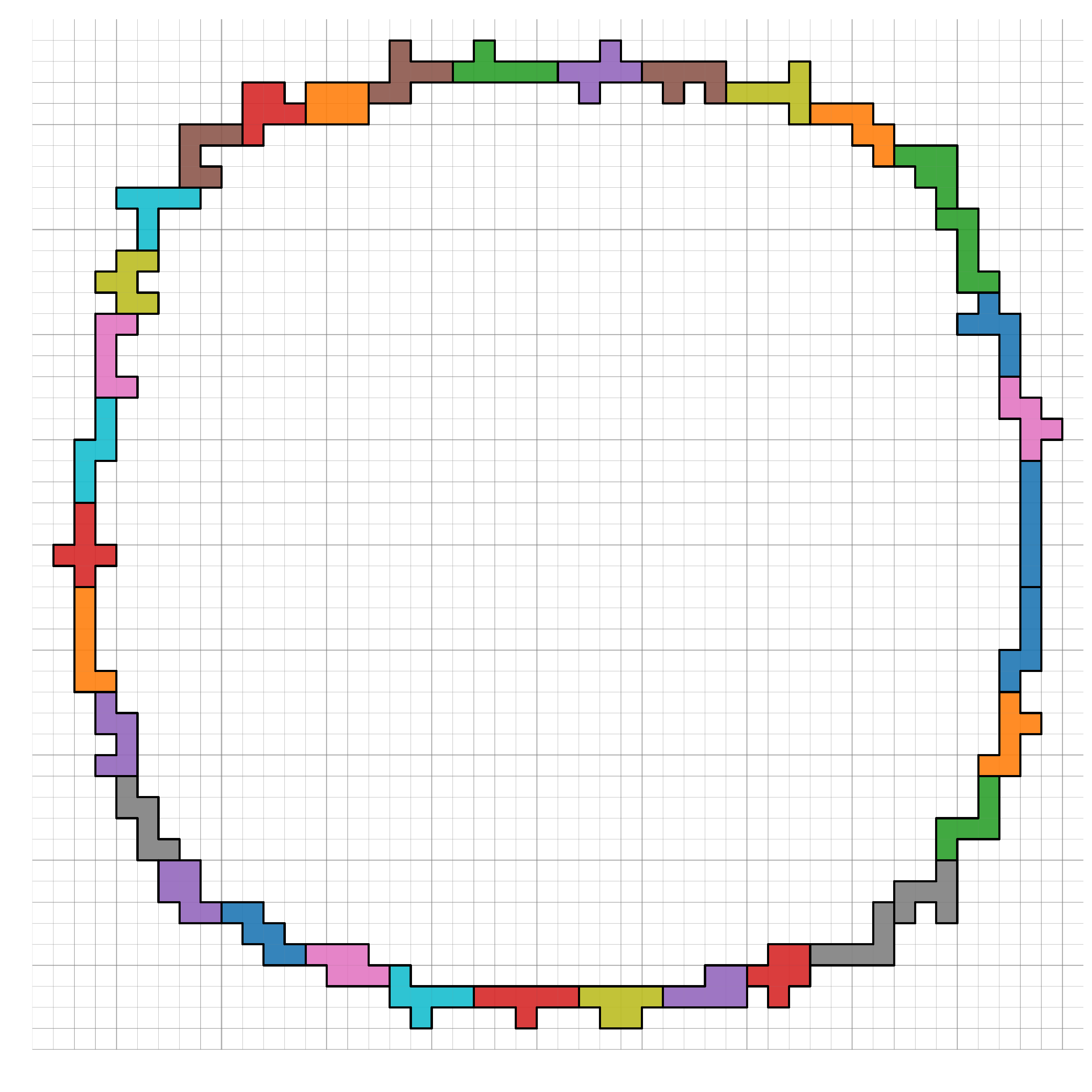}
  \caption{A hexomino fence enclosing the maximum possible area of $1597$ cells.}
  \label{fig:1597}
\end{figure}

A fence enclosing $1597$ cells is shown in Figure~\ref{fig:1597}. It
was found with a stochastically adapted linear programming approach~\cite{github}, improving the previously best known area of
$1586$~\cite{previous-german}.
The remainder of this
section proves that no fence encloses $1598$ or more cells and then
determines all optimal fences.
In Section~\ref{sec:opt-dp}, we begin with a simpler problem
for the dynamic program:
We assume that it has already been decided
\emph{how} each hexomino is used, i.e., which pair of attachment edges
is used,
and we only search for the optimal order of
the hexominoes.
We develop various optimizations to make this simplified
problem feasible.
In Section~\ref{sec:hexomino-proof}, we describe a Branch-and-Bound
preprocessing step that, for a given target area, narrows down the
possible ways of using the hexominoes to a short list of candidates
that might reach this area. For these candidates, we run the dynamic programming
search
of Section~\ref{sec:opt-dp}
to
confirm that for none of them can the hexominoes be ordered to achieve an area of at least~$1598$.

A similar dynamic programming approach
has been
  recently used by 
  Har-Peled, Harb, and He~\cite{harpeled2026catchkgridpoints} for a
  problem of a different flavor: enclosing a given number of grid
  points by the shortest rope, or more precisely, finding
a set of grid points whose convex hull has the smallest possible perimeter.

\paragraph{Interactive lab.} In the \emph{hexomino lab} (\url{https://erikaroldanroa.github.io/fence-challenge/hexomino-lab/}) one can try building a hexomino fence and interact with one optimal solution.

\subsection{Optimizations for dynamic programming
}\label{sec:opt-dp}

\subsubsection{Reference vectors}\label{sec:reference-vectors}

\begin{figure}[htb]
  \centering
  \includegraphics[page=1]{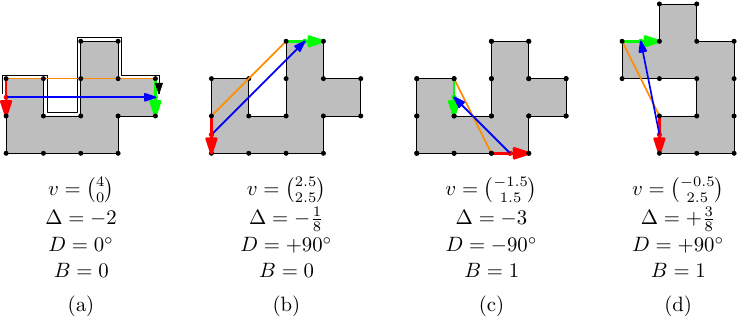}
  \caption{Some pairs of attachment edges, incoming in red and
    outgoing in green, together with the corresponding reference
    vectors $v$ (in blue), corrections $\Delta$, rotations $D$, and parity bits~$B$.}
  \label{fig:refvec-examples}
\end{figure}

By Lemma~\ref{lem:interval},
consecutive pieces around the largest region 
share an attachment edge, and each piece has an \emph{incoming} and an \emph{outgoing} attachment edge. We orient every attachment edge from the enclosed region towards the outside. The midpoints of the attachment edges are the key to a more compact description of a fence, see Figure~\ref{fig:optimal-fence-refpoly}.
The following definition builds on the geometric observations at the end of Section~\ref{sec:boundary}, and Figure~\ref{fig:refvec-examples}
illustrates it with a few examples.

\begin{definition}\label{def:reference-vector}
  Consider a piece that is used with an incoming attachment edge with
  midpoint $m_{\mathrm{in}}$ and an outgoing attachment edge with
  midpoint $m_{\mathrm{out}}$, and let $f_p, \ldots, f_q$ be the edges
  of the piece that lie on the boundary of the enclosed region between
  the two attachment edges, oriented
clockwise around the piece
  as in Section~\ref{sec:boundary}.
  \begin{enumerate}
    \item The \emph{reference vector} of the piece is $v := m_{\mathrm{out}} - m_{\mathrm{in}}$.
    \item The \emph{correction} $\Delta$ of the piece is the signed
      area enclosed by the closed path that starts at
      $m_{\mathrm{in}}$, follows the incoming attachment edge inwards
      to the starting point $s_p$ of $f_p$, then the boundary edges
      $f_p, \ldots, f_q$, 
      then the outgoing attachment edge outwards to $m_{\mathrm{out}}$, and finally returns to $m_{\mathrm{in}}$ along $-v$.
    \item The \emph{rotation} $D$ of the piece is the signed turning angle (counterclockwise turns counted positively) accumulated when walking from the incoming attachment edge along $f_p, \ldots, f_q$ to the outgoing attachment edge, plus $180^\circ$.
    \item The \emph{parity bit} $B$ is the Manhattan distance between
      the start points of the two attachment edges (that is, their inner endpoints), modulo~$2$.
    \end{enumerate}
  Concatenating the reference vectors of the pieces $\pi_1,\ldots,\pi_k$ around the largest region of a fence gives a closed polygon, the \emph{reference polygon}. It is a (not necessarily simple) $k$-gon, which we traverse counterclockwise.
\end{definition}

The quantities $v$, $\Delta$, $D$, and $B$ are invariant under translation, so they can be
 precomputed for every piece and every choice of attachment edges,
 just like the Shoelace contributions $\Delta A_{p,q}$ in
 Section~\ref{sec:boundary}. Since the vectors are attached to the midpoints of the attachment edges, which have half-integer coordinates, all corrections are multiples of~$\frac18$.

The area of the reference polygon is a good approximation of the area of the fence, and the corrections account exactly for the difference:

\begin{lemma}\label{lem:area-decomposition}
  Let $A_{\mathrm{ref}}$ be the signed area of the reference polygon
  of a fence with pieces $\pi_1, \ldots, \pi_k$ around the largest region
  and corresponding corrections $\Delta_1, \ldots, \Delta_k$. Then the area of the fence is
  \begin{equation}\label{eq:area-decomposition}
    A = A_{\mathrm{ref}} + \Delta_1 + \Delta_2 + \dots + \Delta_k.
  \end{equation}
\end{lemma}

\begin{proof}
The signed area of a closed path is half its Shoelace sum, and the
Shoelace sum is additive over the edges of the path. We add up the
closed paths of Definition~\ref{def:reference-vector} for the pieces
$\pi_1,\ldots,\pi_k$. The reversed reference vectors contribute
$-A_{\mathrm{ref}}$. Every attachment edge is shared by two
consecutive pieces and has the same inner endpoint for both: the
first piece traverses the half edge outwards to the midpoint, and the
second piece traverses it inwards, so these contributions cancel.
What remains are the sections $f_p,\ldots,f_q$ of all pieces, which
together form the boundary walk of the largest region~$\mathcal R$,
by Lemma~\ref{lem:interval} and the discussion after it. By the
Shoelace formula, their contribution is the area enclosed by this walk,
which is $|\mathcal R|=A$, because no piece lies inside~$\mathcal R$
(Section~\ref{sec:pieces-around}). Hence
$\Delta_1+\dots+\Delta_k=A-A_{\mathrm{ref}}$.
\end{proof}

\begin{figure}[htb]
  \centering
  \includegraphics[page=2]{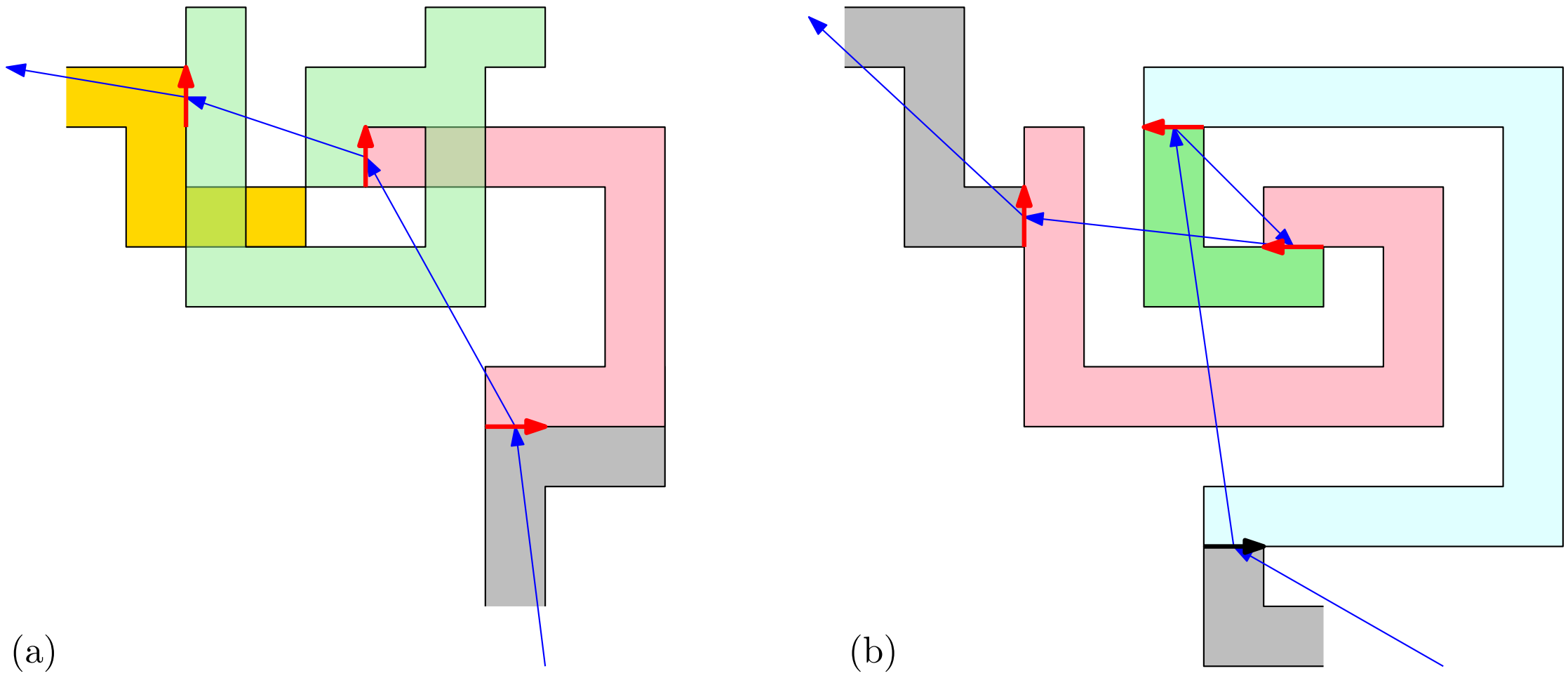}
  \caption{(a) Joining pieces at common attachment edges can lead to overlaps. (b) On the other hand, a reference polygon can intersect itself and still correspond to a valid fence.}
  \label{fig:refvec-overlaps}
\end{figure}

Note that the lemma is a statement about signed areas and does not require the reference polygon to be simple. As Figure~\ref{fig:refvec-overlaps}b shows, a valid fence can have a self-intersecting reference polygon. Conversely, as in Section~\ref{sec:pentomino}, joining pieces at common attachment edges may produce overlapping placements (Figure~\ref{fig:refvec-overlaps}a), for which formula~\eqref{eq:area-decomposition} still yields a well-defined, but meaningless value.

The correction is usually negative. For example, a straight piece whose reference vector runs parallel to the inner boundary at distance $\frac12$ loses a strip of width $\frac12$, which gives $\Delta = -2$ in Figure~\ref{fig:refvec-examples}a.

\paragraph{The rotation.}
The rotation $D$ is a multiple of $90^\circ$. Most pieces in a good
fence have parallel attachment edges and rotation $D = 0^\circ$, while
the \emph{corner pieces} have $D = +90^\circ$.
Unusual pieces can also have $D = 180^\circ$ or $D = -90^\circ$ (an inward-turning piece), although it is hard to imagine them in an optimal fence. Swapping the incoming and outgoing edge changes the sign of $D$. We do not exclude any possibility a priori but consider all pairs of edges as potential attachment edges, and for hexominoes this results in rotations between $-360^\circ$ and $+360^\circ$. Note that $D$ is determined as a signed angle, not only modulo $360^\circ$. At each attachment edge, the walk along the boundary makes a U-turn from the outgoing edge of one piece to the incoming edge of the next, which the added $180^\circ$ in the definition of the rotation counts as a left turn. With this convention the turning angles along the edge boundary add up correctly, and since the edge boundary is a closed curve traversed counterclockwise, the rotations of all pieces of a fence satisfy
\begin{equation}\label{eq:rotation-sum}
  D_1 + D_2 + \dots + D_k = 360^\circ.
\end{equation}

\paragraph{The parity bit.}
The vectors $v_i'$ between the start points (inner endpoints) of the attachment edges
form a closed
polygon. 
For a closed integer $k$-gon with edge vectors
$\binom{x_i}{y_i}$,
$\sum_{i=1}^k (x_i+y_i)$ is zero, and therefore even.
We
allow the pieces to be reflected and rotated, which may swap
$x_i$ and $y_i$ and negate them.
Negating an integer does not change its parity, and swapping them
does not change
the parity
of $x_i+y_i$.
We can therefore conclude:
The sum of the parity bits must be even, no
matter how the pieces are reflected or rotated:
\begin{equation}\label{eq:parity-sum}
  B_1 + B_2 + \dots + B_k \equiv 0 \pmod 2.
\end{equation}
This constraint will allow us to exclude combinations of reference
vectors, because they cannot form a closed reference polygon.

\paragraph{The repertoire of reference quadruples.}

Algorithm~\ref{alg:reference-repertoire} computes, for a single
polyomino, all possible quadruples $(v, \Delta, D,B)$. It enumerates the
pairs of attachment edges exactly as the precomputation in
Section~\ref{sec:boundary}; here, $w_1, \ldots, w_t$ are the vertices
of the polyomino boundary in clockwise order. For hexominoes, the
boundary is always a simple cycle, so this vertex sequence is well
defined. When this is no longer the case, which can happen from heptominoes
on, we use a ``right-hand rule'', which means that we always choose the
rightmost edge from the current vertex. Starting with octominoes, we
can have $n$-ominoes
that have one or more separate boundaries on the inside,
e.g., the $3 \times 3$ square with a missing center. In such cases we
only use the outer
boundary, since it makes no sense to put pieces into the holes of
other pieces.

\begin{algorithm}[H] 
  \caption{\textsc{ReferenceRepertoire}$(p)$}
  \label{alg:reference-repertoire}
  \begin{algorithmic}[1]
    \State \textbf{Input:} polyomino $p$
    \State $\mathcal{V} \gets \emptyset$
    \For{each distinct orientation $P$ of $p$}
      \State $(w_1, \ldots, w_t) \gets$ boundary vertices of $P$ in clockwise order, indices taken modulo $t$
      \For{$s = 1$ to $t$} \Comment{incoming attachment edge $(w_s, w_{s+1})$}
        \For{$e = s+1$ to $s+t-1$} \Comment{outgoing attachment edge $(w_e, w_{e+1})$}
          \State $m_{\mathrm{in}} \gets \frac12 (w_s + w_{s+1})$ and $m_{\mathrm{out}} \gets \frac12 (w_e + w_{e+1})$
          \State $v \gets m_{\mathrm{out}} - m_{\mathrm{in}}$
          \State $(u_0, u_1, \ldots, u_r) \gets (m_{\mathrm{in}}, w_{s+1}, w_{s+2}, \ldots, w_e, m_{\mathrm{out}})$ and $u_{r+1} \gets u_0$ \Comment{closed path}
          \State $\Delta \gets \frac12 \sum_{k=0}^{r} \det(u_k, u_{k+1})$ \Comment{Shoelace formula}
          \State $D \gets \sum_{k=s}^{e-1} \operatorname{turn}(w_{k+1}
          - w_k,\, w_{k+2} - w_{k+1})+ 180^\circ$
          \Comment {rotation}
          \State
          $\binom{x'}{y'}\gets  w_{s+1}-w_e$; \ 
        $B\gets ( x' +   y')
        \bmod 2$ \Comment{parity bit}
          \State add $(v, \Delta, D,B)$ to $\mathcal{V}$
        \EndFor
      \EndFor
    \EndFor
    \State \Return $\mathcal{V}$
  \end{algorithmic}
\end{algorithm}

Here $\operatorname{turn}(f,f')$ for two unit grid vectors $f$ and $f'$ is $+90^\circ$ if $\det(f,f') =1$,
$-90^\circ$ if $\det(f,f')=-1$ and $0^\circ$ if $f$ and $f'$ point in
the same direction. The case that they point in opposite directions is
excluded. For the straight piece in
Figure~\ref{fig:refvec-examples}a, the walk turns right twice, so
the sum of the turns is
$
-180^\circ$, and $D = 0^\circ$.

The sums in lines 10 and 11 need not be computed from scratch as
written;
for better efficiency,
they could, for example,
be evaluated incrementally, although the computational effort
of this preprocessing step
is negligible in the context of the whole process.

\subsubsection{Eight-way partitioning}\label{sec:8way}

Lemma~\ref{lem:area-decomposition} separates the area of a fence into
two parts. The sum of the corrections does not depend on the order in
which the pieces are joined, only on how each piece is used. The area
of the reference polygon, on the other hand, depends on the order, but
only through the reference vectors. If we have chosen a reference vector of every piece, it is then enough to arrange these vectors into a polygon of maximum area, and at the end we simply add the fixed sum of the corrections.

This is promising because optimal fences have a very regular
structure. There are four corner pieces, and all other pieces have
parallel attachment edges. Moreover, the attachment edges of almost
all pieces are chosen so that the reference vector is as long as
possible:
For a ``roundish'' polygon such as the reference polygon that we are
trying to build,
the area grows quadratically with its perimeter.
Thus, a slightly longer reference vector usually gains more than
enough area to compensate a deterioration through a different
correction. This is not a strict rule; the U-pentomino, for example,
is not used with its longest reference vector in optimal pentomino
fences,
see \autoref{fig:use-opt}.
 Hexominoes are a similar case, see \autoref{sec:overlap}.
When all reference vectors are fixed up to rotation and reflection, we arrive at the following auxiliary problem.

\goodbreak

\begin{quote}
   \textbf{Problem~\maxpolygon: the largest polygon from prescribed edge vectors, up to the 8 grid symmetries.}
  Given vectors $\binom{a_1}{b_1}, \ldots, \binom{a_N}{b_N}$, each
  vector may be used in any of its up to eight rotated and reflected
  variants $\binom{c_i}{d_i} \in \bigl\{\binom{\pm a_i}{\pm b_i},
  \binom{\pm b_i}{\pm a_i}\bigr\}$. Form a closed $N$-gon $\Omega$
  with these
  edge vectors that encloses the largest possible area.
\end{quote}

\goodbreak

Once a variant $\binom{c_i}{d_i}$ is selected for every vector, the problem is easy: the polygon closes if and only if $\sum_{i=1}^N \binom{c_i}{d_i} = \binom00$, and among all polygons with a given multiset of edge vectors, the convex one obtained by sorting the vectors by direction has the largest signed area. We write $F(v_1, \ldots, v_N)$ for the area of this convex polygon, with $F(v_1, \ldots, v_N) = -\infty$ if $v_1 + \dots + v_N \neq \binom00$. Since there are $8^N$ possible choices of variants, the difficulty lies in choosing them.

\begin{figure}[htb]
  \centering
  \tabskip=0pt
  \halign to \textwidth{\hfil#\hfil\tabskip=0pt plus 1fil & \hfil#\hfil\tabskip=0pt\cr
  \includegraphics{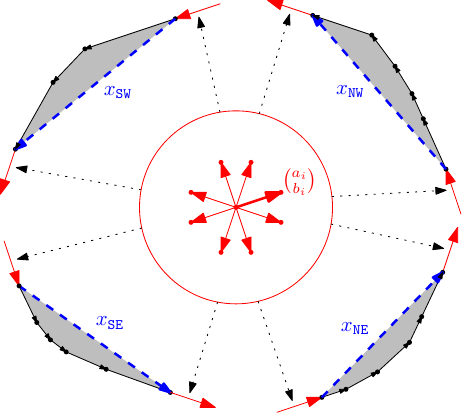}&
  \includegraphics[page=2]{data/img/add-one-vector}\cr\noalign{\smallskip}
  (a)&(b)\cr}
  \caption{(a) The $i$-th vector can extend the partial solution in one of eight possible ways. (b) At the end, the four chains are joined into a closed polygon.}
  \label{fig:adding}
\end{figure}

We solve \maxpolygon{} with dynamic programming. The optimal polygon is
convex, so we can build it by adding the vectors one at a time in the
order from the ``steepest'' (closest to a $\pm 45^\circ$ diagonal) to
the ``flattest'' (closest to a coordinate axis). The polygon then
grows in four separate convex chains, one around each of the four
diagonal directions, see Figure~\ref{fig:adding}. Following the
dominant direction of their edges, we call them the NE, NW, SW and SE
chains, and we denote the vectors from the first to the last vertex of
each chain by $x_{\mathrm{NE}}, x_{\mathrm{NW}}, x_{\mathrm{SW}},
x_{\mathrm{SE}}$, and call them the \emph{sum vectors} of the chains.

The eight rotated and reflected variants of a vector correspond to the eight places where it can be attached: at the beginning or at the end of one of the four chains. If we assume $0 \le b_i \le a_i$, then the variant $\binom{a_i}{b_i}$ is prepended to the NE chain and $\binom{b_i}{a_i}$ is appended to it, and the other three chains are obtained by rotating by multiples of $90^\circ$: the NW chain receives $\binom{-b_i}{a_i}$ at its beginning and $\binom{-a_i}{b_i}$ at its end, the SW chain $\binom{-a_i}{-b_i}$ and $\binom{-b_i}{-a_i}$, and the SE chain $\binom{b_i}{-a_i}$ and $\binom{a_i}{-b_i}$. Since the vectors are processed in order of decreasing steepness, each chain stays convex and sorted by direction.

Each chain, together with the segment between
its endpoints, encloses a certain area
(shaded gray in Figure~\ref{fig:adding}). The algorithm stores, for every combination of sum vectors, the largest total area
\begin{equation*}
  \mathrm{OPT}_i[x_{\mathrm{NE}}, x_{\mathrm{NW}}, x_{\mathrm{SW}}, x_{\mathrm{SE}}]
\end{equation*}
of the four chains that can be obtained from the first $i$ vectors
such that these four vectors are the resulting sum vectors. When
adding the variant $e = \binom{c_i}{d_i}$ to a chain with sum vector
$x$, the area of that chain grows by the triangle spanned by $x$ and
$e$, i.e.,\ by $\frac12 \det(x, e)$ if $e$ is appended at the end and
by $\frac12 \det(e, x)$ if it is prepended at the beginning, while the
sum vector becomes $x + e$. Thus, each entry of $\mathrm{OPT}_{i-1}$
produces eight candidate values for different entries of $\mathrm{OPT}_i$, of which we keep only the best value per entry. As in Section~\ref{sec:pentomino}, the crucial point is that we do not need to know the orientations and the order of the first $i-1$ vectors to extend the solution; the four sum vectors and the corresponding area are all the information we need.

After all $N$ vectors have been added, we consider the entries of $\mathrm{OPT}_N$ with $x_{\mathrm{NE}} + x_{\mathrm{NW}} + x_{\mathrm{SW}} + x_{\mathrm{SE}} = \binom00$. These are exactly the ones where the four chains close up to a polygon (Figure~\ref{fig:adding}b). To the stored value we add the area of the quadrilateral formed by the four sum vectors,
\begin{equation*}
  \tfrac12 \bigl[\det(x_{\mathrm{NE}}, x_{\mathrm{NW}}) + \det(x_{\mathrm{NE}} + x_{\mathrm{NW}},\, x_{\mathrm{SW}})\bigr],
\end{equation*}
to obtain the area of the whole polygon, and the largest value obtained in this way is the optimum.

\paragraph{Memory requirement.}
In our applications all vectors have integer or half-integer
coordinates. With $\ell := \sum_{i=1}^N (|a_i| + |b_i|)$, there are at
most $O(\ell^8)$ possible quadruples of sum vectors, so the running time
and the space requirement are polynomially bounded by $O(N\ell^8)$. If only the optimal value is needed and not the solution itself, the tables of previous steps can be discarded, and $O(\ell^8)$ memory suffices. For hexominoes this is still far too large, and we need to restrict the search space further.

\subsubsection{Using a prescribed target area}\label{sec:target-area}

The size of the tables is the main bottleneck.  In practice, most of
their entries cannot lead to an optimal solution. We can discard
entries if we specify in advance a target $A_{\mathrm{target}}$ that we
want to reach.

For a table entry $\mathrm{OPT}_i[x_{\mathrm{NE}}, x_{\mathrm{NW}},
x_{\mathrm{SW}}, x_{\mathrm{SE}}]$, the final polygon consists of the
four chains built so far and the $N-i$ remaining vectors. Replacing
each chain by its sum vector removes exactly the stored area, so the
area of any completion equals the stored value plus the area of a
polygon whose edges are the four sum vectors and the remaining
vectors. Ignoring the directions and only considering the lengths of
these $N - i + 4$ edges, we can bound the latter area by the largest
area of any polygon with the given side lengths.
It is well known (see, for example,
  \cite[Theorem 12.5a]{niven1981maxima}) that
this maximum is
attained if the polygon is inscribed in a circle.
Its radius $R$ can be determined as the solution of a nonlinear equation in
one variable, namely $\sum_k 2\arcsin\bigl(\ell_k / (2R)\bigr) = 2\pi$
for given side lengths $\ell_k$ (with an appropriate modification if the
center of the circle lies outside the polygon), and once $R$ is determined, the area
can be easily computed.

Adding the stored value gives an upper bound for the area that can be
reached from this partial solution. If this bound, plus the known sum
of the corrections, is smaller than the
target area $A_{\mathrm{target}}$, we do not need to store the entry. This reduces the size of the tables and, as a consequence, also the running time. To prove that no fence is better than the best known one with area $A^*$, we can choose $A_{\mathrm{target}} = A^* + 1$, since fence areas are integers; in our computations we use the slightly smaller threshold $A^* + 0.999$ to guard against rounding errors.

We determine the radius $R$ by a simple binary-search procedure, narrowing down the interval in which $R$ is contained.
The sum of the angles (which should equal~$2\pi$
in the above formula)
is a decreasing function of $R$, whereas the sum of the areas of the triangles spanned by the chords of length~$\ell_k$ is an increasing function of $R$. Thus, if $R$ is too big for the angle condition,
but the area is smaller than the threshold, we can safely discard the entry.
On the other hand, if $R$ is too small for the angle condition,
but the area is already bigger than the threshold, we already know that the entry will pass
the test. We hardly ever need to determine the radius $R$ very precisely.
 On average, a pruning test takes less than 8 binary search iterations.
The formulas involve
arcsine 
and square-root computations, plus the addition of $N$ numbers.
There is no propagation of floating-point rounding errors over several levels, and
the safety buffer of $0.001$ is generous in comparison to these errors. 

To give the reader an impression of how effective this pruning technique is: we can
confirm in $150$ minutes of computing time and with $130$\,GB of
memory that the reference vectors of the fence in
Figure~\ref{fig:1597} cannot be arranged into a fence of area larger
than $1597$. (By itself, this does not prove optimality, because
other choices of reference vectors might give better results.)
For heptominoes and octominoes, this pruning procedure is not sufficient, and we will additionally restrict the search space heuristically in Section~\ref{sec:hepto-okto}.

\subsection{Proof of optimality with DP after Branch-and-Bound (B\&B) preprocessing}\label{sec:hexomino-proof}

The dynamic program of Section~\ref{sec:8way} assumes that a fixed
reference vector has been chosen for every piece. To prove
Theorem~\ref{thm:1597}, we therefore first reduce the enormous number
of possible choices of reference vectors to a handful of
candidates. We apply the Branch-and-Bound technique, using an upper
bound based on Minkowski averages of polygons. The proof consists of
four steps.

\subsubsection{Step 0: Computing the reference vector repertoire of a hexomino}

For each hexomino $i = 1, \ldots, 35$, Algorithm~\ref{alg:reference-repertoire} determines the \emph{repertoire}
\begin{displaymath}
  \mathcal{V}_i = \bigl\{(v_i^1, \Delta_i^1, D_i^1,B_i^1), (v_i^2,
  \Delta_i^2, D_i^2, B_i^{2}), \ldots, (v_i^{l_i}, \Delta_i^{l_i}, D_i^{l_i}, B_i^{l_i}
  )\bigr\}
\end{displaymath}
of all possible reference vectors with their corrections,
rotations, and parity bits, where all rotated and reflected variants are listed
explicitly. A fence selects one quadruple $(v_i, \Delta_i, D_i, B_i) \in
\mathcal{V}_i$ for every piece and arranges the pieces in a suitable
order.
To allow a piece to be unused, which is harmless by Section~\ref{sec:def-area}, we add the null option
$(v, \Delta, D, B)=(0,0,0,0)$ to the repertoire of every piece.
Adding this option does not change the outcome of the computations below:
Step~1 returns exactly the same candidate sequences with and without it,
both for the threshold $1597.999$ and for the threshold $1596.999$ of Section~\ref{sec:all-6}.
By Lemma~\ref{lem:area-decomposition}, the maximality of the convex
arrangement, and conditions~\eqref{eq:rotation-sum} and~\eqref{eq:parity-sum}, every fence satisfies
\begin{multline}\label{eq:general-problem}
  \text{Area} \le S_0 := \max \bigl\{\, F(v_1, \ldots, v_{35}) +
  \Delta_1 + \dots + \Delta_{35} \;\bigm|\; (v_i, \Delta_i, D_i, B_i) \in
  \mathcal{V}_i,\\ D_1 + \dots + D_{35} = 360^\circ,\
    B_1 + B_2 + \dots + B_{35} \equiv 0 \pmod 2\,\bigr\}.
\end{multline}
This is already a rather precise formulation of the fence problem: The
only condition that is not checked is that the resulting fence is
free of overlaps.

\subsubsection{An upper bound via Minkowski averages}\label{sec:minkowski}

\begin{figure}[htb]
  \centering
  \includegraphics[page=2]{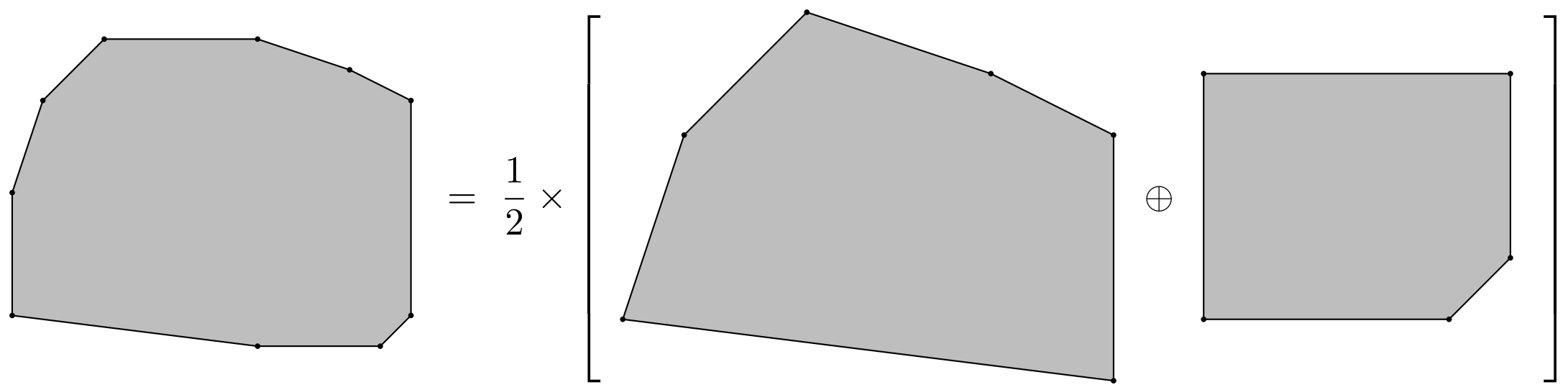}
  \caption{The Minkowski average of two polygons.}
  \label{fig:minkowski2}
\end{figure}

We can compute an easy upper bound for the problem \maxpolygon{} by the technique of
Minkowski averages.
The \emph{Minkowski average} of two polygons $\Omega_1$ and $\Omega_2$ is
\begin{equation*}
  \Omega' = \tfrac12 [\Omega_1 \oplus \Omega_2] := \bigl\{\, \tfrac12 (u_1 + u_2) \bigm| u_1 \in \Omega_1,\ u_2 \in \Omega_2 \,\bigr\},
\end{equation*}
where the term in brackets is the \emph{Minkowski sum}, the element-wise sum of the two sets. For two convex polygons, the Minkowski average is easy to construct, see Figure~\ref{fig:minkowski2}: we merge the counterclockwise directed edge vectors of $\Omega_1$ and $\Omega_2$ into one list sorted by direction, scale all vectors by~$\frac12$, and concatenate them to form $\Omega'$.

\begin{figure}[htb]
  \centering
  \includegraphics{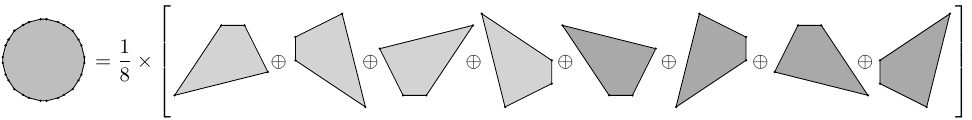}
  \caption{The Minkowski average of the eight rotated and reflected copies of a polygon.}
  \label{fig:minkowski8}
\end{figure}

The famous Brunn--Minkowski inequality\footnote{See, for example,
\url{https://en.wikipedia.org/wiki/Brunn-Minkowski_theorem}.}
implies
that, for convex sets $\Omega_1$ and $\Omega_2$ of equal area, the area of $\Omega'$
is at least as large as the area of $\Omega_1$.
Moreover, the inequality is strict unless
 $\Omega_1$ is a positive homothet (scaled and translated copy) of $\Omega_2$.
This also holds for
averages of more than two sets, and we apply it to the eight rotated
and reflected copies of a hypothetical optimal solution $\Omega$ of
\maxpolygon{}, see Figure~\ref{fig:minkowski8}. Let $\rho_1, \ldots, \rho_8$ be
the eight symmetries of the grid that fix the origin, i.e.,\ the
rotations by multiples of $90^\circ$ and the reflections in the axes
and the diagonals, and consider
\begin{equation*}
  \Omega' = \tfrac18 \bigl[ \rho_1 \Omega \oplus \rho_2 \Omega \oplus \dots \oplus \rho_8 \Omega \bigr].
\end{equation*}
The nice property of this expression is that we know exactly what $\Omega'$
looks like, even though we do not know $\Omega$: every edge vector $\binom{a_i}{b_i}$ of $\Omega$
appears in $\Omega'$ scaled by $\frac18$ in each of its eight
variants. Since all polygons $\rho_k \Omega$ have the same area as $\Omega$, the
area of $\Omega'$ cannot be smaller than that of $\Omega$. We thus obtain an
easily computable upper bound for \maxpolygon.

We expect the optimum solution to be rather round and not as spiky as
in Figure~\ref{fig:minkowski8}. Thus, the Minkowski average $\Omega'$ will
be just a slightly better rounded version of $\Omega$ and the true area of
$\Omega$ will be close to the bound.
(If every given vector occurs with a multiplicity divisible by $8$, then $\Omega'$ consists of whole vectors and is itself an optimal solution.)

The area of $\Omega'$ can be expressed by a simple formula:
\begin{equation}\label{eq:minkowski-bound}
  F_M := \frac{1}{16} \sum_{i=1}^N \sum_{j=1}^N
        C\bigl[\tbinom{a_i}{b_i}, \tbinom{a_j}{b_j}\bigr]
\end{equation}
{with}
\begin{equation}
  \label{eq:termC}
  C\bigl[\tbinom{a_i}{b_i}, \tbinom{a_j}{b_j}\bigr]
  := \max\bigl\{ \max(|a_i|,|b_i|) \cdot (|a_j|+|b_j|),\ (|a_i|+|b_i|) \cdot \max(|a_j|,|b_j|) \bigr\}.
\end{equation}
So $F_M$ is a sum over all pairwise products of vector lengths, where
one vector is measured in the maximum norm
$\max(|a|,|b|)$
and the other in the $1$-norm  $|a|+|b|$, choosing the combination that gives the larger product. If all vectors are transformed into the first octant and sorted by slope, i.e.,\ $0 \le b_i \le a_i$ and $\frac{b_1}{a_1} \le \frac{b_2}{a_2} \le \dots \le \frac{b_N}{a_N}$, the formula simplifies to
\begin{equation}\label{eq:minkowski-bound-sorted}
  F_M = \frac18 \Biggl(\, \sum_{1 \le i < j \le N} a_i (a_j + b_j) + \frac12 \sum_{i=1}^N a_i (a_i + b_i) \Biggr),
\end{equation}
since the flatter vector $\binom{a_i}{b_i}$ is always measured in the maximum norm $a_i$ and the steeper vector $\binom{a_j}{b_j}$ in the $1$-norm $a_j + b_j$.

\begin{figure}[htb]
  \centering
  \includegraphics{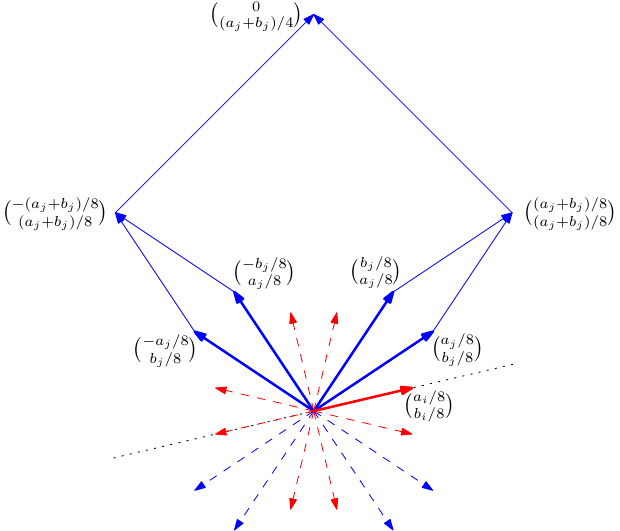}
  \caption{Combining $\frac18\binom{a_i}{b_i}$ with all rotated and reflected copies of $\frac18\binom{a_j}{b_j}$.}
  \label{fig:combine-i-j}
\end{figure}

\begin{theorem}
  The optimum area for the problem \textup{\maxpolygon} is less than or equal to
  $F_M$ from formula~\eqref{eq:minkowski-bound-sorted}.
\end{theorem}
\begin{proof}
  We have already argued why the
  area of the convex polygon $\Omega'$ is an upper bound.
  It remains to prove that this area equals~\eqref{eq:minkowski-bound-sorted}.
  
The polygon $\Omega'$ is centrally symmetric.
  A centrally symmetric
  convex polygon with edge vectors $\pm w_1, \pm w_2, \ldots,
  \allowbreak \pm w_h$
  is a translate of a \emph{zonogon}, a Minkowski sum
   $\bigoplus_k[0,w_k]$ of line segments.
  The area of such a zonogon is
  $\sum_{1 \le k < l \le h} |\det(w_k,  w_l)|$.
  (This follows easily from the fact that the zonogon can be tiled by
  the parallelograms spanned by the pairs $w_k,w_l$\footnote
  {\url{https://en.wikipedia.org/wiki/Zonogon}}, see~\cite[\thetag{56--57}]{Shephard-1974}.)
  Every vector $\binom{a_i}{b_i}$ contributes to $\Omega'$ the four pairs of opposite vectors $\pm w_i^{(1)}, \ldots, \pm w_i^{(4)}$ with $w_i^{(1)} = \frac18\binom{a_i}{b_i}$, $w_i^{(2)} = \frac18\binom{b_i}{a_i}$, $w_i^{(3)} = \frac18\binom{-b_i}{a_i}$ and $w_i^{(4)} = \frac18\binom{-a_i}{b_i}$ in the upper half-plane. Consider a pair $i \le j$. Since the whole configuration is symmetric under the eight grid symmetries, and $|\det|$ is invariant under them, it suffices to compute the contribution of $w_i^{(1)}$ and multiply it by $4$. As $\frac{b_i}{a_i} \le \frac{b_j}{a_j}$, all four determinants $\det(w_i^{(1)}, w_j^{(k)})$ are nonnegative (see Figure~\ref{fig:combine-i-j}), so we can combine them:
  \begin{equation*}
    \sum_{k=1}^4 \det\bigl(w_i^{(1)}, w_j^{(k)}\bigr)
    = \det\Bigl(w_i^{(1)},\, \sum_{k=1}^4 w_j^{(k)}\Bigr)
    = \begin{vmatrix} a_i/8 & 0 \\ b_i/8 & (a_j + b_j)/4 \end{vmatrix}
    = \frac{a_i (a_j + b_j)}{32}.
  \end{equation*}
  Multiplying by $4$ gives the term $\frac18 a_i (a_j + b_j)$ in
  \eqref{eq:minkowski-bound-sorted}. For $i = j$ this counts every
  determinant $|\det(w_i^{(k)}, w_i^{(k')})|$ with $k \neq k'$ twice
  (and for $k = k'$ the determinant is $0$), so in this case the
  result is weighted by the factor $\frac12$.
\end{proof}

This bound is often remarkably sharp. For example, the reference polygon of the optimal hexomino fence in Figure~\ref{fig:1597} has area $1663$, while the bound gives $1666.5625$.

\subsubsection{Step 1: Branch-and-Bound over vector sequences}\label{sec:bnb}

For now, we ignore the rotation condition $D_1 + \dots + D_{35} =
360^\circ$
and the parity condition
$    B_1 + \cdots\allowbreak + B_{35} \equiv 0 \pmod 2$, and we also no longer distinguish between rotated and reflected reference vectors. This leads to the reduced repertoire
\begin{equation*}
  \mathcal{W}_i = \bigl\{(\bar v_i^1, \Delta_i^1), (\bar v_i^2, \Delta_i^2), \ldots, (\bar v_i^{m_i}, \Delta_i^{m_i})\bigr\},
\end{equation*}
where $\bar v_i^j$ stands for an equivalence class of up to eight
rotated and reflected vectors, and $\Delta_i^j$ is the largest
correction of any copy of this vector in $\mathcal{V}_i$. For the
hexominoes, the sets $\mathcal{W}_i$ contain between $12$ and $23$
elements, including the null option. Our goal is to find a small list of vector sequences $(\bar
v_1, \ldots, \bar v_{35})$ with $\bar v_i \in \mathcal{W}_i$ that
contains all optimal ones. Writing $\bar F(\bar v_1, \ldots, \bar
v_{35})$ for the optimal value of
\maxpolygon{} with the given vector classes, we have
\begin{equation}\label{eq:S1}
  S_0 \le S_1 := \max \bigl\{\, \bar F(\bar v_1, \ldots, \bar v_{35}) + \Delta_1 + \dots + \Delta_{35} \bigm| (\bar v_i, \Delta_i) \in \mathcal{W}_i \,\bigr\},
\end{equation}
and by Section~\ref{sec:minkowski},
\begin{equation}\label{eq:maxFbar}
  \bar F(\bar v_1, \ldots, \bar v_{35}) \le F_M(\bar v_1, \ldots, \bar v_{35}) = \frac18 \sum_{1 \le i < j \le 35} C[\bar v_i, \bar v_j] + \frac1{16} \sum_{i=1}^{35} C[\bar v_i, \bar v_i].
\end{equation}
We distribute each correction $\Delta_i$ evenly, with a share of $\frac1{35}$, over the $35$ terms of \eqref{eq:maxFbar} in which $\bar v_i$ occurs. Maximizing every term separately over the available choices then gives the bound
\begin{multline}\label{eq:wide-bound}
  \mathcal{F}(\mathcal{W}_1, \ldots, \mathcal{W}_{35}) := \sum_{1 \le i < j \le 35} \max \bigl\{\, \tfrac18 C[\bar v_i, \bar v_j] + \tfrac{\Delta_i + \Delta_j}{35} \bigm| (\bar v_i, \Delta_i) \in \mathcal{W}_i,\ (\bar v_j, \Delta_j) \in \mathcal{W}_j \,\bigr\} \\
  + \sum_{i=1}^{35} \max \bigl\{\, \tfrac1{16} C[\bar v_i, \bar v_i] + \tfrac{\Delta_i}{35} \bigm| (\bar v_i, \Delta_i) \in \mathcal{W}_i \,\bigr\} \;\ge\; S_1,
\end{multline}
which can be computed easily by combining all pairs of candidates in each term. If every set $\mathcal{W}_i$ contains a single element, it coincides with $F_M + \sum_i \Delta_i$. For larger sets the bound is weaker, since in each of the $35$ terms containing $\bar v_i$ the vector can be chosen independently, whereas a solution uses the same $\bar v_i$ in all terms. For the original sets $\mathcal{W}_i$, the bound~\eqref{eq:wide-bound} has the value $1633.964$.

Conceptually, we now enumerate all vector sequences in $35$ nested
loops, selecting $(\bar v_1, \Delta_1) \in \mathcal{W}_1$, then $(\bar
v_2, \Delta_2) \in \mathcal{W}_2$ and so on, and output every sequence
with $F_M(\bar v_1, \ldots, \bar v_{35}) + \Delta_1 + \dots +
\Delta_{35} \ge 1597.999$.
Done naively, this would take astronomically long. Instead, on level $i$, when $(\bar v_1, \Delta_1), \ldots, (\bar v_i, \Delta_i)$ are already fixed, we evaluate the bound
\begin{equation*}
  \mathcal{F}\bigl(\{(\bar v_1, \Delta_1)\}, \ldots, \{(\bar v_i, \Delta_i)\}, \mathcal{W}_{i+1}, \ldots, \mathcal{W}_{35}\bigr),
\end{equation*}
and if it is already smaller than $1597.999$, the threshold cannot be
exceeded anywhere in this branch of the enumeration tree and we prune
it. As explained above, the bound becomes more accurate the more the
choices are restricted, which makes the pruning very effective, see
Algorithm~\ref{alg:bnb}.

\begin{algorithm}[H]  
  \caption{\textsc{BranchAndBound}$(i, \mathcal{W}_1, \ldots, \mathcal{W}_{35})$}
  \label{alg:bnb}
  \begin{algorithmic}[1]
    \State \textbf{Input:} level $i$, sets $\mathcal{W}_1, \ldots, \mathcal{W}_{i-1}$ already reduced to single elements
    \If{$\mathcal{F}(\mathcal{W}_1, \ldots, \mathcal{W}_{35}) < 1597.999$}
      \State \Return \Comment{prune this branch}
    \EndIf
    \If{$i = 36$}
      \State output the candidate sequence $(\bar v_1, \ldots, \bar v_{35})$ \Comment{all sets are singletons}
      \State \Return
    \EndIf
    \For{each $(\bar v_i, \Delta_i) \in \mathcal{W}_i$}
      \State \textsc{BranchAndBound}$(i+1, \mathcal{W}_1, \ldots, \mathcal{W}_{i-1}, \{(\bar v_i, \Delta_i)\}, \mathcal{W}_{i+1}, \ldots, \mathcal{W}_{35})$
    \EndFor
  \end{algorithmic}
\end{algorithm}

Starting with \textsc{BranchAndBound}$(1, \mathcal{W}_1, \ldots,
\mathcal{W}_{35})$, this procedure
visits about $30$ million
nodes of the enumeration tree
and finds $7128$
candidate sequences.
In our implementation, the runtime was roughly two hours.\footnote{The running time could be reduced
  further by eliminating some vectors beforehand: If fixing a single $(\bar
  v_i, \Delta_i)$ while keeping all other sets $\mathcal{W}_j$
  variable already gives a bound below $1597.999$, this vector can be
  removed once and for all. For simplicity, we did not do this.
  We remark that, with the threshold $1594.999$ instead of $1597.999$, the search
  produces $136\,543$ candidates in $5$ hours; among these there are
  also solutions containing pieces with rotation
  $-90^\circ$. Fortunately, no such solutions occur with the original
  bound $1597.999$,
  because the subsequent steps would not work for them.}
The candidate sequences are written to a file and passed on to the next step.

\subsubsection{Step 2: Taking into account the rotation and the parity}
\label{sec:rotations}

Next, we subject each of the $7128$ vector sequences from Step~1 to a
stricter test that also takes the
rotation condition~\eqref{eq:rotation-sum} and the parity
condition~\eqref{eq:parity-sum}
into account.
Consider a fixed vector sequence 
$((\bar v_1,\Delta_1),\ldots,
(\bar v_{35},\Delta_{35}))$.
For each~$i$, let $\mathcal{U}_i$ be the set of triples $(\Delta_i^j,
D_i^j, B_i^j)$ of those quadruples $(v_i^j, \Delta_i^j, D_i^j, B_i^j)
\in \mathcal{V}_i$ whose vector $v_i^j$ coincides with the prescribed
class $\bar v_i$ after a suitable rotation or reflection. Then we
define the optimal total correction under these constraints:
\begin{multline*}
  \Delta^{\mathrm{rot}} := \max \bigl\{\, \Delta_1 + \dots +
  \Delta_{35} \bigm| (\Delta_i, D_i, B_i) \in \mathcal{U}_i,\\ D_1 + \dots + D_{35} = 360^\circ,\
    B_1 + B_2 + \dots + B_{35} \equiv 0 \pmod 2\,\bigr\}
\end{multline*}
This can be computed very quickly by dynamic programming over the partial
sums of the rotations and the parities. If $F_M(\bar v_1, \ldots, \bar v_{35}) +
\Delta^{\mathrm{rot}} < 1597.999$, the sequence can be excluded.

\begin{figure}[htb]
  \centering
  \includegraphics[page=3]{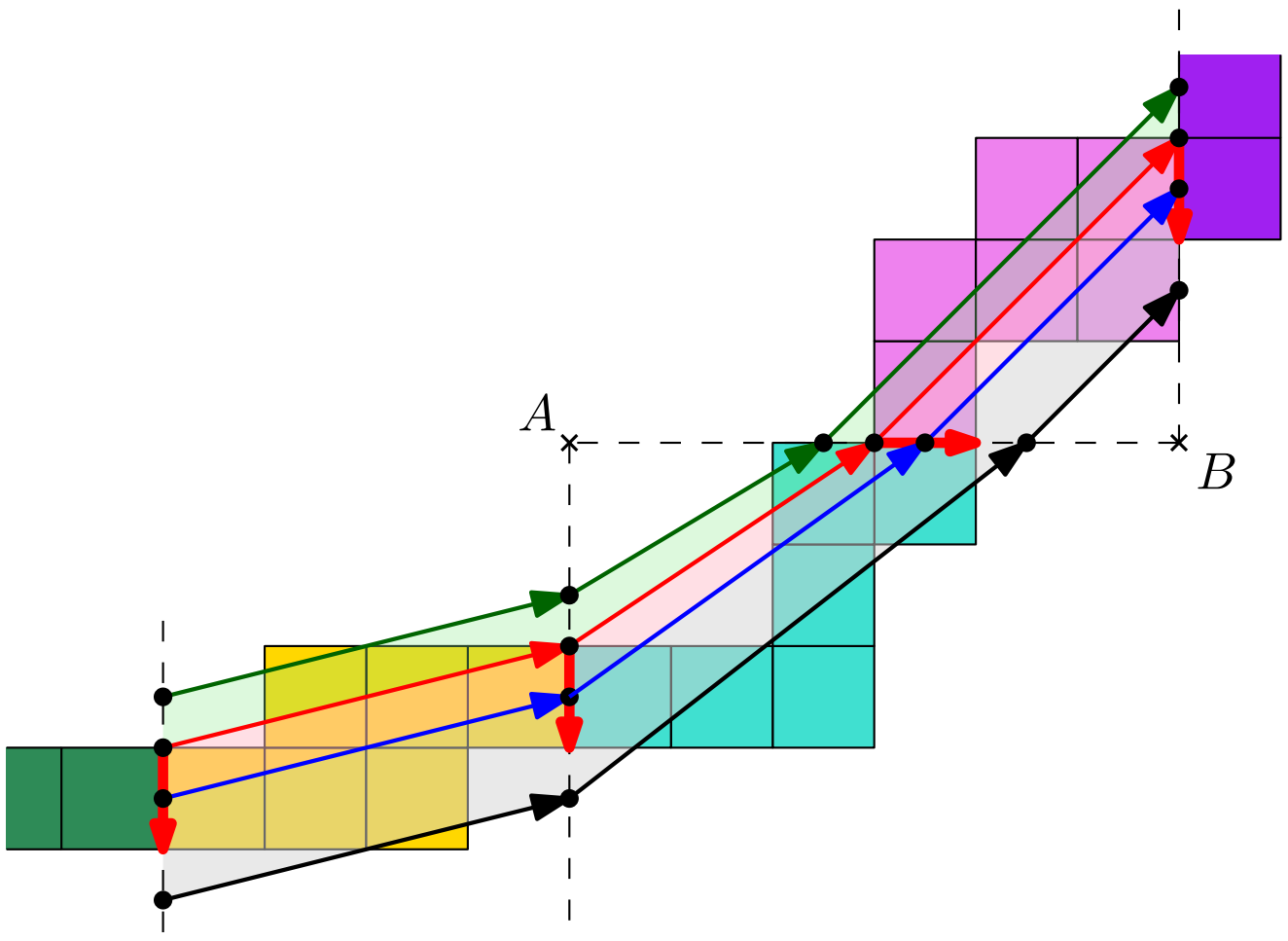}
  \caption{The reference vectors can be anchored anywhere
    on the line through the attachment edges, not just at the
    midpoints (blue) or at the inner endpoints (red).
    The endpoints can  even lie outside the polyominoes (green and black).
  The correction has to be adjusted by adding or subtracting the
  shaded area.
  We see a yellow straight
  piece (rotation $0^\circ$), followed by a turquois
``corner piece'' (rotation $+90^\circ$), followed by a violet
piece with rotation $-90^\circ$. Offsetting the endpoints
works differently for the three types of pieces.
}
  \label{fig:alternative-attachment}
\end{figure}

From now on, only the \emph{multiset} $\{\bar v_1, \ldots, \bar
v_{35}\}$ of vectors matters, and eight multisets remain. Two of them
can be excluded by the following trick:
We let the reference vectors start and end at some other points
instead of the midpoints of the attachment edges,
see \autoref{fig:alternative-attachment}.
%
For each offset between $-1.5$ and $+1.5$, in steps of $0.5$, we move
the endpoint on the line through the attachment edge
from its ``normal'' position at the midpoint of the attachment edge.
See
\autoref{fig:alternative-attachment} for examples of offsets $-1,
-0.5, 0, +1$.
This changes the
reference vectors of the corner pieces, and all corrections have to be
adjusted accordingly.
For the exact area, this makes no difference because
the two adjustments to the polygon area
and $\Delta^{\mathrm{rot}}$
cancel each other,
but the \emph{bounds} $F_M$ change (unless the four corner vectors
are equal $45^\circ$-vectors).
So we recompute the bound
$F_M(\bar v_1, \ldots, \bar v_{35}) +
\Delta^{\mathrm{rot}}$ with the adjusted vectors and penalties.
Indeed, for two of the eight vector sets, the bounds
become smaller than $1597.999$
when we choose $s=-0.5$, i.e., 
 the \emph{inner} endpoints of the adjacent edges (the red reference
 vectors in \autoref{fig:alternative-attachment}).

We made sure that
we don't ``overcross'' the reference vector and change the
sign of one of its components: We don't try an offset when
this would happen, i.e., when
an
endpoint would be moved to the point $A$ or $B$ marked in the
figure or beyond.
Otherwise the reference vectors would approach the attachment edges
from the wrong side, and there would automatically be an overlap
when the reference vectors form a convex polygon.
(Lemma~\ref{lem:area-decomposition} remains valid, no matter where we
put the endpoints on the attachment edge.)

Each of the remaining six multisets consists of four half-integral
corner vectors
otherwise only integral vectors.
The program checks
that the
corner pieces can be used only
with rotation $\pm 90^\circ$, and the integral pieces
 only with rotation
$0^\circ$.
This forces the corner pieces to have rotation $+ 90^\circ$,
which is important for the
next step: it tells us that the four corner pieces point in four
different directions, so the tables of the dynamic program can be
initialized accordingly. Step~2 takes only a few seconds.

\subsubsection{Step 3: Dynamic programming with fixed vectors}

Finally, we run the dynamic program from Section~\ref{sec:8way} with
the pruning from Section~\ref{sec:target-area} on each of the six
multisets, with each chain initialized with one of the corner
pieces. As in Step~2, the reference vectors are attached to the inner
endpoints of the attachment edges. This makes the reference vectors of
the corner pieces integral, while the reference vectors of all other
pieces remain unchanged, since their two attachment edges are
parallel. The main effect is that the work of the dynamic programming
algorithm
is reduced.

There is one exception, 
where this would make a corner piece
flatter than a ``straight'' piece with $D = 0^\circ$: the reference
vector $\binom{1.5}{2.5}$ would become $\binom12$, which has a
smaller slope than the vector $\binom23$. This would interfere with the
order of processing the vectors, so in this case we
went in the opposite direction:
we attached the reference vectors to the \emph{outer} endpoints of the
attachment edges instead, which gives the corner piece the reference
vector $\binom23$.
(By making the outward shift large enough, one can always make the
corner vectors steeper than the other vectors, provided
no other vector has slope $\pm1$. As argued in \autoref{sec:overlap},
vectors
with slope $\pm1$ other than the four corner vectors
are not likely to occur in optimum fences.)

This computation took about $16$ hours for all six multisets, with a
peak memory usage of about $130$\,GB on a machine with $300$\,GB of
main memory. No solution with an area larger than $1597.999$ was found, which completes the proof of Theorem~\ref{thm:1597}. Note that, unlike in Section~\ref{sec:pentomino}, no check for overlaps is needed here: every step only relaxes the problem, so the computations yield a valid upper bound, and the matching lower bound is given by the explicit fence in Figure~\ref{fig:1597}.
Section~\ref{sec:ilp-hex} confirms the upper bound $1597$ independently with integer programming, under an additional assumption on how the hexominoes are used.

\subsection{Step 4: Identification of all optimal solutions}\label{sec:all-6}

We ran Steps~1--3 with the weaker bound
$1596.999$,
and we added another step that finds all solutions with
 at least this area, i.e.,\ all optimal solutions, by backtracking
 through
 the dynamic-programming table.
Step~1 now produces
$21\,283$ candidate sequences (instead of $7128$),
Step~2 finds  $14$ multisets of reference vector classes (instead of 8).
All but two of these multisets 
consist of four corner pieces
with rotation $+90^\circ$ and otherwise only pieces with rotation
$0^\circ$.
The two exceptional multisets
have six half-integer vectors, which can be used with
rotation $\pm 90^\circ$:
$\left\{4\times\binom{2.5}{2.5}, 2\times\binom{3.5}{2.5}\right\}$
and
$\left\{5\times\binom{2.5}{2.5}, 1\times\binom{3.5}{2.5}\right\}$
To get a full rotation of
$+360^\circ$, five of these vectors must be used
with $+90^\circ$ and one with
$-90^\circ$. There are two possibilities for
that one vector, and we try each one individually.
 Three of these cases can be excluded with an appropriate choice of
 offset $s$, but
in one case
the
bound best meets the threshold 1597 exactly.
This is the multiset where the
vectors $\left\{3\times\binom{4}{0},2\times\binom{5}{0},1\times\binom{6}{0},4\times\binom{3}{1},
  11\times\binom{4}{1},3\times\binom{5}{1},5\times\binom{4}{2}\right\}$
are used with rotation $0^\circ$ five
vectors $\left\{4\times\binom{2.5}{2.5},1\times\binom{3.5}{2.5}\right\}$
are used with rotation $+90^\circ$,
and
one vectors $\binom{2.5}{2.5}$ is
 used with rotation $-90^\circ$.
For $s=-1.5$, the bound is 1597, and the other offsets yield larger bounds.\footnote
{Varying the offset $s$ in the vicinity of $-1.5$ did not improve
  the bound, and
  1597 is a local minimum. In fact,
the term $  C\bigl[\tbinom{a_i}{b_i}, \tbinom{a_j}{b_j}\bigr]$
from
\eqref{eq:termC}, which appears in the formula for $F_M$, is a maximum
of two expressions, and it changes from one expression to the other
when the order of the slopes of the reference vectors
$\tbinom{a_i}{b_i}$ and $\binom{a_j}{b_j}$ changes. This explains why
$F_M$ is not smooth at half-integral offsets, and
local minima naturally appear there.}

To deal with this case, we recall from Section~\ref{sec:minkowski} that the
Brunn--Minkowski inequality is strict unless the summands are
positively homothetic. Since the eight symmetric summands of which the Minkowski
average is taken (\autoref{fig:minkowski8}) have the same area, this
can
only happen if they are translates of each other; in other words, if
the reference polygon has all eight symmetries of the square.
However, this would require that the 3 vectors of
class $\binom51$
 should be distributed
evenly over the eight octants, which is of course impossible. (Also, none of the
half-integral reference vectors
is changed to $\binom51$ or a multiple of this vector when $s=1.5$ is added or
subtracted from its components).

We need not be worried about rounding errors in the computation of the
bound 1597: The computation of $F_M$ according to
formula~\eqref{eq:minkowski-bound}
involves multiplication of half-integer or integer coordinates in
pairs,
additions and subtractions, and finally a division by 16.
Since the size of the operands is very moderate,
all these calculations come out exactly, even with ordinary
floating-point arithmetic.
The corrections are multiples of $\frac18$ and are also exact.

We were lucky that this case could be just barely excluded, since 
Step~3 would not have been prepared to deal with the case when there
are not precisely four pieces with rotation $+90^\circ$, to be used as
corner pieces.

Eventually, Step~2 reduces the 14 options
 to $10$ multisets of reference vector classes
(instead of 6 with the bound 1598).
All of them have only four vectors with rotation $+90^\circ$,
which can be used as corner vectors, and all remaining pieces with
rotations $0^\circ$;
so Step~3 can proceed.

Step~3 reduces the 10 multisets
to a single
multiset 
\begin{displaymath}
  1 \times \tbinom60,\ 2 \times \tbinom50,\ 3 \times \tbinom40,\ 3 \times \tbinom51,\ 11 \times \tbinom41,\ 4 \times \tbinom31,\ 6 \times \tbinom42,\ 1 \times \tbinom32,\ 4 \times \tbinom{2.5}{2.5},
\end{displaymath}
running for about two days.

Step~4 then runs a modified dynamic-programming algorithm that takes care
not to count the same solution several times.
For example, if we place a vector $\binom42$ into the second quadrant
and
then  another vector $\binom42$ into the fifth quadrant, or vice
versa, our original dynamic-programming algorithm would consider this as two
different solutions. To avoid such overcounting, we insist that equal vectors
are inserted into their quadrants in order:
After inserting $\binom42$ into the fifth quadrant,
future copies of the same vector $\binom42$ can only be inserted into the
fifth or higher-numbered quadrants.
The indices of the dynamic-programming table must be extended to
record the last quadrant into which a piece was inserted.
We additionally maintain a counter of the number of solutions for each
entry,
which gives us then directly the number of optimal solutions (as an
independent check for the subsequent phase).

For the starting configuration, we are in the lucky situation that
the four corner vectors are equal and lie on the diagonal direction, but
the ``steepest'' vector $\binom32$, which is the first one to be
inserted, is unique. Thus we insert
it
into some fixed quadrant. In this way, we have eliminated all
symmetries from the start.

After the dynamic-programming table is filled, we process it in the
reverse direction and retrace all paths that lead to the optimum.
The tables become very big, so we store them on disk to free 
main memory, and read them back in the reverse phase.
Even so, the peak memory requirement is about 300 GByte.

In this way, after running for about one day,
Step~4 produces $244$ different optimal reference polygons up to symmetry.


There is one family of vectors, where
 different \emph{parallel} vectors exist, namely the axis-parallel ones
$  1 \times \tbinom60,\ 2 \times \tbinom50,\ 3 \times \tbinom40$.
Our program inserts these three families in a fixed order.
However, consecutive parallel vectors
can be permuted without changing the area,
giving rise to a different reference polygon
if the vectors are different.
In our 244 polygons, the six straight vectors are
always distributed to the four sides as $6=1+1+2+2$ (not necessarily
in this order).
The most frequent pattern (228 cases) has
pairs two edges of length 4 on one side, and one edge of length 4 with one edge of
length 5 on another side. The last two vectors can be swapped, giving rise
to $ 228 \times2=456$ different orders.
There are 4 similar cases, where
an edge of length 4 lies together with the edge of
length 6, giving rise
to $ 4 \times2=8$ different orders.

There are
8 polygons that have two sides with
edges of length $4+5$ respectively.
(\autoref{fig:1597} shows such a polygon.)
They give rise to $8\times 2\times 2=32$ different orders.
Only 4 polygons
pair two edges of length 4 on one side, and
two edges of length 5 on another side, and each of these needs to be counted
only once.
In total, there are $456+8+32+4 = 500$ different reference polygons.


%
%
Finally, as in the case of pentominoes (Section~\ref{sec:all-5}),
hexominoes that realize the same reference vector are
interchangeable.
 For the last two classes $1 \times \tbinom32$ and $ 4
\times \tbinom{2.5}{2.5}$ (the four corner pieces), the same set of
five hexominoes can be used, which gives a factor of $5!$.
Apart from these five pieces, no piece can be used for two different
  reference vectors.
Thus, we have
\begin{displaymath}
  2! \cdot 3! \cdot 3! \cdot 11! \cdot 4! \cdot 6! \cdot 5!
  \times 500 = 2\,979\,773\,153\,280\,000\,000
\end{displaymath}
different optimal hexomino orders. Finally, we determine all optimal solutions, which, as defined in Section~\ref{sec:all-5}, are distinct if they differ in their cells. Similarly to the Y- and Z-pentominoes, now $10$ of the hexominoes can be rotated or reflected without changing their reference vector and correction. Thus, there are
\begin{displaymath}
  2^{10} \times
  2! \cdot 3! \cdot 3! \cdot 11! \cdot 4! \cdot 6! \cdot 5!
  \times 500 =
  3\,051\,287\,708\,958\,720\,000\,000 
  \approx 3.05 \cdot 10^{21} 
\end{displaymath}
optimal solutions to the hexomino fence problem.

The above calculation rests on the assumption that different reference
polygons give rise to different fences. It is conceivable that two
pieces share more than one edge, as in \autoref{fig:shoelace-multiple},
and then the fence would not have a
unique
reference polygon.
The following section excludes this possibility.

\subsection{Overlapping pieces and multiple shared edges}
\label{sec:overlap}

\begin{figure}[htb]
  \centering
  \includegraphics{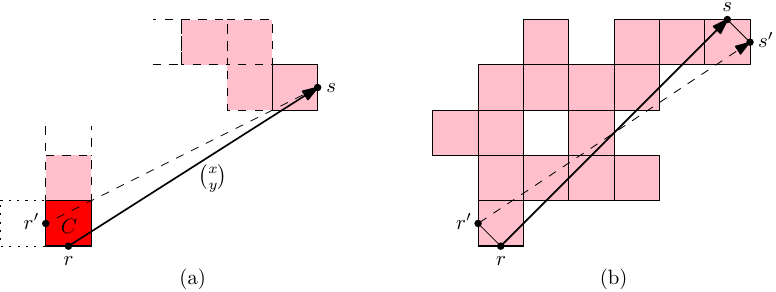}
  \caption{(a) The largest distance between two edge midpoints $r$ and
    $s$. (b) When a reference vector of slope $\pm1$ is available, there
  is always a longer vector.}
  \label{fig:longest-edge}
\end{figure}

We need an elementary lemma about longest reference vectors.

\begin{lemma}\label{lem:integer-vector}
  Let $P$ be a finite, nonempty set of unit cells of the square grid.
Then a longest vector between two 
edge midpoints of the boundary of $P$ lies in 
$
\mathbb{Z}^2$,
and the two edges on which these midpoints lie are parallel.
\end{lemma}

\begin{proof}
  A boundary edge is vertical or horizontal; we call its midpoint a
  \emph{V-point} in $\mathbb{Z}\times(\mathbb{Z}+\tfrac12)$ or an
  \emph{H-point} in $(\mathbb{Z}+\tfrac12)\times\mathbb{Z}$,
  respectively. The difference of two points of the same type lies in
  $\mathbb{Z}^2$, whereas the difference of a V-point and an H-point
  lies in $(\mathbb{Z}+\tfrac12)^2$.
  Suppose that $r,s$ are of different types, and consider the vector
  $s-r=\binom xy$ with $x,y\in\mathbb{Z}+\tfrac12$.
We will show that $r,s$ cannot have the largest distance.
  
Without loss of generality, we assume that $0<y\le x$, $r$ is an
H-point, and $s$ is a V-point, as in \autoref{fig:longest-edge}a.
(The other cases are symmetric.)
Consider the red cell $c_0$ on whose boundary $r$ lies.
Assume first that the cell to the left of $c_0$
(the dotted cell in \autoref{fig:longest-edge}a) does not belong to $P$.
Then we take the midpoint $r'$ of the left boundary edge of $c_0$,
and
\begin{displaymath}
  \left\lVert s-r'\right\rVert^2 = (x+\tfrac12)^2 + (y-\tfrac12)^2
  = x^2+y^2 + (x-y) + \tfrac12
  \ge x^2+y^2 + \tfrac12
  > x^2+y^2 =
  \left\lVert s-r\right\rVert^2,
\end{displaymath}
that is, $r,s$ cannot have the maximum distance.
If the dotted cell to the left of $c_0$ belongs to $P$, we walk left until we
hit a vertical boundary edge of $P$. The midpoint of that edge will have an
even larger distance from $s$ than $r'$.
Hence the two endpoints of a longest vector are of the same type: they lie on
parallel edges, and their difference lies in~$\mathbb{Z}^2$.
\end{proof}

\begin{figure}[htb]
  \centering
  \includegraphics[page=2]{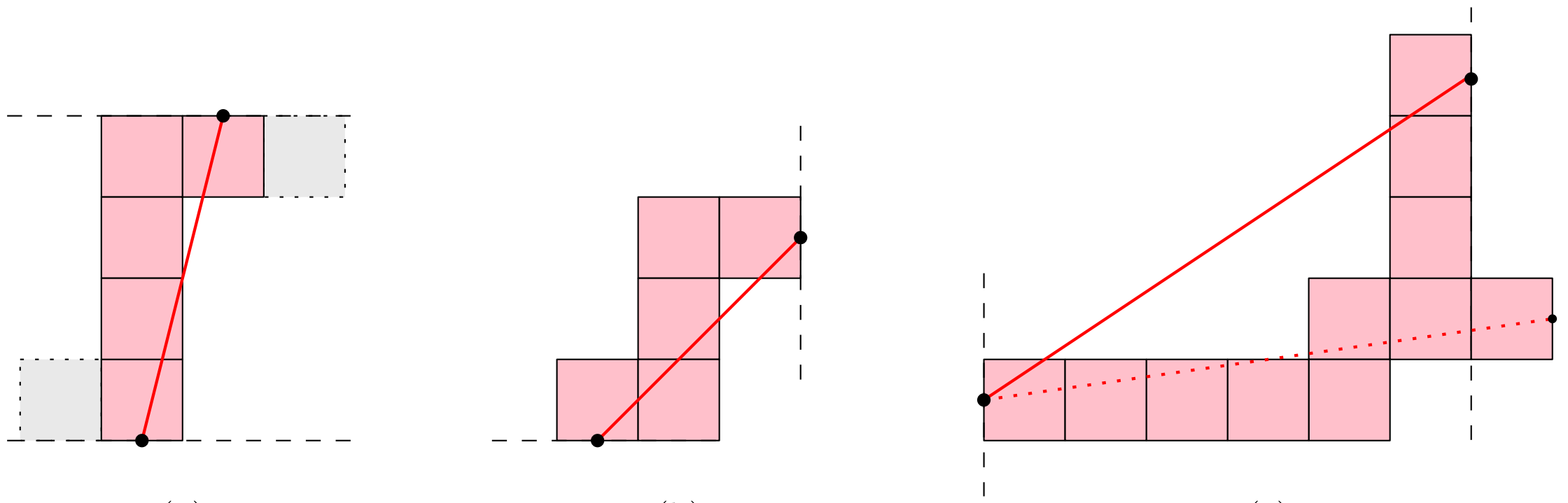}
  \caption{(a) A separated piece with two parallel attachment edges.
    If the red reference vector is the longest possible reference vector,
    the two shaded cells cannot be part of the piece.
    (b) A separated corner piece.
    (c) The longest reference vector, of length $\sqrt{6^2+4^2}=\sqrt{52}$, does not make
    the piece separated:
 The
alternative
choice of reference vector, which would lead to a separated piece, has length
$\sqrt{7^2+1^2}=\sqrt{50}$.}
  \label{fig:sandwiched}
\end{figure}

We call a piece $P$ with reference vector $v$ 
\emph{separated} if $P$ does not cross the lines through
the edges at the endpoints of $v$, see
\autoref{fig:sandwiched}.

\begin{lemma}\label{non-overlap}
  Suppose that:
  \begin{enumerate}
  \item 
    The reference polygon is convex.
  \item It
    contains
  exactly four half-integral ``corner vectors'' in
 $(\mathbb{Z}+\tfrac12)^2$, one pointing into each of the four quadrants.
\item It contains at least one axis-parallel vector in
  each of the four directions.
\item   Every piece is separated.
\item For each corner piece,
  the reference vector is the longest possible vector between
midpoints of a horizontal and a vertical
 boundary edge.
\item For each of the remaining pieces,
 the reference vector is the longest possible vector between two
 boundary-edge midpoints.
  \end{enumerate}
  Then
    consecutive pieces do not overlap, and their common attachment edge
    is their only common edge.  
\end{lemma}

\begin{proof}
 By assumption~2, the four corner vectors divide
  the cyclic sequence of reference vectors into four {runs} of
  consecutive non-corner vectors. Since the polygon is convex, the
  vectors are ordered by direction. Hence the run between the corner
  vectors pointing into the first and into the second quadrant consists of the
  vectors whose directions lie strictly between the directions of
  these two corner vectors. We call it the \emph{right run}. Its
  vectors $\binom xy$ have directions in the interval $(0^\circ,180^\circ)$, and
  so $y>0$.
We concentrate on the right run; the other three runs are treated similarly.

  By assumption~3, the right run contains a vector
  $\binom 0y$; 
  in
  particular, no run is empty, and no two corner pieces are
  consecutive.

By Lemma~\ref{lem:integer-vector}, every piece in the run has parallel attachment
edges, and hence all attachment edges in the right run are parallel.
In fact, since the run contains a vertical vector, they are horizontal.
The horizontal lines through the attachment edges
cut the plane into strips, and
since the pieces are separated, they are confined to their strips and
cannot overlap.

We now show that two pieces cannot share more than one edge.
Consider the initial part of the right run, up to the first vertical
vector.
These vectors are of the form $\binom xy$ with $x\ge0$ (and $y>0$),
see \autoref{fig:sandwiched}a.
If $P$ is such a piece, no cell of $P$ can touch the bottom supporting
line to the left of the lower endpoint of the vector, and no cell can
touch the top supporting line to the right of the upper endpoint,
because otherwise there would be a longer reference vector.
The next piece $P'$ above $P$ has the same
properties. Thus, on the line separating them, $P$ touches only at or
to the left of the chosen attachment edge, and $P'$ only at or to the
right of it, and so that edge is the unique common edge of $P$ and~$P'$.
The only case where this argument fails is when
a right-leaning vector
$\binom xy$ with $x>0$ 
is followed by a left-leaning vector
$\binom {x'}{y'}$ with $x'<0$, 
but since there is at least
one vertical vector with $x=0$, this case does not arise.

The argument extends to the two corner pieces at the lower and
upper end of the right run: The maximality condition 5 applies
analogously,
and we conclude that the whole run including the corner pieces
is overlap-free, and consecutive pieces intersect only in one edge.
\end{proof}

As mentioned above,
the six axis-parallel vectors 
for hexominoes are distributed over the four directions in such a way
that each direction gets at least one vector.
Thus, Assumption~3 is fulfilled.

The remaining properties are clear, with one exception:
Assumptions 4 and 5 are unfortunately not fulfilled for one piece:
The cyan T-shaped piece in the upper left side of
\autoref{fig:1597} is not treated as required:
Its longest reference vector is $\binom40$, as opposed to the
reference vector $\binom23$ of length $\sqrt{13}$ with which it is
used.
It can also be used with
reference vector $\binom{5/2}{5/2}$ as a corner piece; in fact, it is
one of the five hexominoes mentioned above that can be corner pieces.

Nevertheless, it can be checked that in the position in which it is
used in
\autoref{fig:1597},
it has the properties needed in
the proof: It is separated, and besides the attachment edge, it
touches each separating line only on the side permitted in the proof.
However, the T-piece can also be used in a reflected position where
the notch points upwards and not to the left. In that position it
is not separated. Here we have to make an ad-hoc argument: The piece $P'$
above cannot intersect the notch, because otherwise it would have to
touch the separating line on the wrong side of the attachment edge.
It cannot even touch the notch: It cannot touch the notch from the
right (or from the left), because then it would touch the separating line.
$P'$ must touch the notch from the top. On the other
hand, 
its reference vector $\binom{x'}{y'}$ must satisfy $\frac{x'}{y'}\ge
\frac 23$, because otherwise the reference polygon would not be
convex, and then one can argue that the top edge of the notch would
give rise to a longer reference vector.

When the T-piece is used as a corner piece, the assumptions of the
lemma are fulfilled.
All other pieces are used with their longest possible reference
vectors, and they are separated.
 All the remaining assumptions of the lemma are fulfilled as well, for all
500 reference polygons.

We conclude that the reference polygon is unique for all solutions.
We also conclude that none of the pieces in the right run can overlap,
or in any other run.
What we have proved does not exclude the possibility that, say, pieces
from the right run and the bottom run intersect.
However, these runs move away from each other so clearly, and
we refrain from proving formally that they cannot intersect.

We mention that the polyomino pieces
for all choices of reference vectors shown in
\autoref{fig:use-opt} are separated, as can be checked by inspection. This confirms that the
5760 pentomino fences in
\autoref{sec:all-5}
are all valid,
independent of the computer check that was run by for this small number of fences.
\paragraph{Larger polyominoes.}

\autoref{fig:sandwiched}c is an example where the longest possible
reference vector does not lead to separation.
It can be checked that
this is a smallest example in terms of the size
of the polyomino: one needs at least 11 cells for this phenomenon to
occur with a strict maximum. One can construct an example with a tie
in
length between $\binom 43$
and $\binom 50$ with eight cells, but not with seven.
For heptomino fences, Lemma~\ref{non-overlap} can be applied
if longest reference vectors are chosen.

\subsection{Independent confirmation of hexomino optimality with integer programming}
\label{sec:ilp-hex}

We now want to demonstrate that the ILP encoding, which we used in our previous paper~\cite{previous-german}, is also able to determine that there is no hexomino fence with an area of more than $1597$. To this end, we reduce each of the $35$ hexominoes to a straight or L-shaped piece with the same effect on the fence. A replacement piece does not capture the re-entrant corners of the original hexomino, so the reduction requires a correction (penalty) for the cells in these corners. Figure~\ref{fig:reduction-examples} shows three examples. In total, there are eight shapes to which the 35 pieces can be reduced. Some hexominoes can be feasibly used in two different ways, and in this case both reductions have to be considered. Since different hexominoes and different choices of reduction can lead to the same combination of L-shaped pieces, only $24$ distinct combinations remain.

\providecolor{hexcol}{gray}{0.55}
\providecolor{lcol}{RGB}{50,110,180}

\newcommand{\sq}[3][]{\filldraw[#1, draw=white, line width=.7pt] (#2,#3) rectangle ++(1,1);}
\newcommand{\hexo}[1]{\foreach \x/\y in {#1}{\sq[fill=hexcol]{\x}{\y}}}
\newcommand{\stick}[2][]{\foreach \i in {1,...,#2}{\sq[#1]{0}{\i-1}}}
\newcommand{\Lshape}[3][]{%
  \foreach \i in {1,...,#2}{\sq[#1]{0}{\i-1}}%
  \foreach \i in {1,...,#3}{\sq[#1]{\i}{0}}}

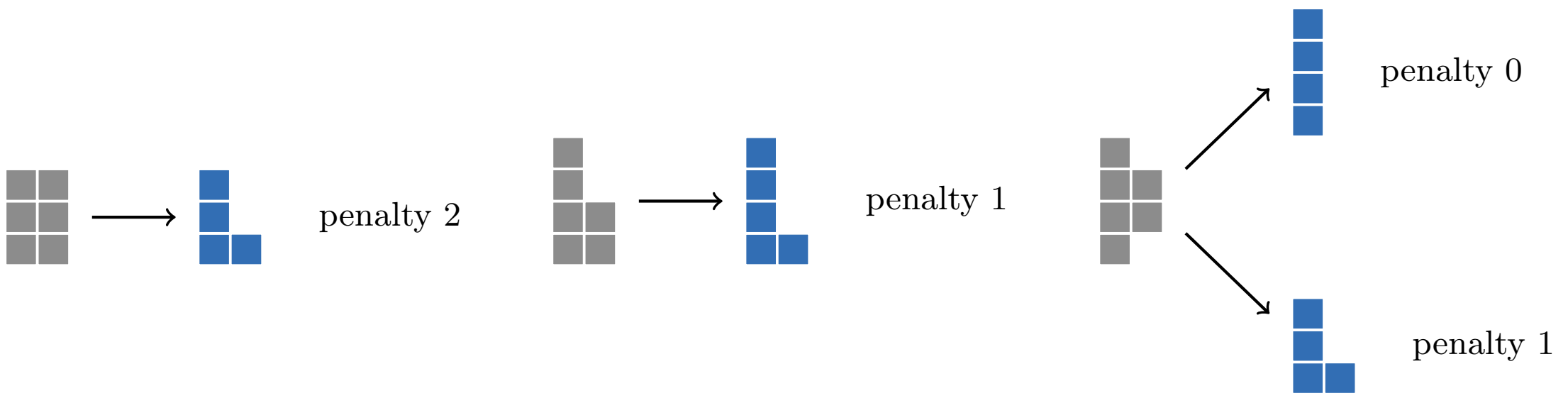
\begin{figure}[htp]
    \centering
    \begin{tikzpicture}[x=3mm,y=3mm,every node/.style={font=\footnotesize}]
        \begin{scope}
            \hexo{0/0,1/0,0/1,1/1,0/2,1/2}
            \draw[->,thick] (2.7,1.5) -- (5.3,1.5);
            \begin{scope}[shift={(6,0)}]\Lshape[fill=lcol]{3}{1}\end{scope}
            \node[anchor=west] at (9.3,1.5) {penalty $2$};
        \end{scope}
        \begin{scope}[shift={(17,0)}]
            \hexo{0/0,0/1,0/2,0/3,1/0,1/1}
            \draw[->,thick] (2.7,2) -- (5.3,2);
            \begin{scope}[shift={(6,0)}]\Lshape[fill=lcol]{4}{1}\end{scope}
            \node[anchor=west] at (9.3,2) {penalty $1$};
        \end{scope}
        \begin{scope}[shift={(34,0)}]
            \hexo{0/0,0/1,0/2,0/3,1/1,1/2}
            \draw[->,thick] (2.7,3) -- (5.3,5.5);
            \draw[->,thick] (2.7,1) -- (5.3,-1.5);
            \begin{scope}[shift={(6,4)}]\stick[fill=lcol]{4}\end{scope}
            \node[anchor=west] at (8.3,6) {penalty $0$};
            \begin{scope}[shift={(6,-4)}]\Lshape[fill=lcol]{3}{1}\end{scope}
            \node[anchor=west] at (9.3,-2.5) {penalty $1$};
        \end{scope}
    \end{tikzpicture}
    \caption{Three examples of the reduction. Left and middle: a hexomino with a single replacement piece. Right: a hexomino that can be replaced in two different ways. The penalty is the correction for re-entrant cells. It is $0$ for the upper possibility of the right piece, because the hexomino can simply be used reflected in the $y$-axis.}
    \label{fig:reduction-examples}
\end{figure}

\def\size{0.25}
\def\lv{\begin{tikzpicture}
    \foreach \y in {0,1,2,3,4} { \draw (0, \y*\size) rectangle (\size, \y*\size + \size); } \draw (\size, 0) rectangle (2*\size, \size);
\end{tikzpicture}
}
\def\liv{\begin{tikzpicture}
    \foreach \y in {0,1,2,3} { \draw (0, \y*\size) rectangle (\size, \y*\size + \size); } \draw (\size, 0) rectangle (2*\size, \size);
\end{tikzpicture}
}

\def\liii{\begin{tikzpicture}
    \foreach \y in {0,1,2} { \draw (0, \y*\size) rectangle (\size, \y*\size + \size); } \draw (\size, 0) rectangle (2*\size, \size);
\end{tikzpicture}
}
\def\liiiii{\begin{tikzpicture}
    \foreach \y in {0,1,2} { \draw (0, \y*\size) rectangle (\size, \y*\size + \size); }
	\draw (\size, 0) rectangle (2*\size, \size);
	\draw (2*\size, 0) rectangle (3*\size, \size);
\end{tikzpicture}
}
\def\livii{\begin{tikzpicture}
    \foreach \y in {0,1,2,3} { \draw (0, \y*\size) rectangle (\size, \y*\size + \size); }
	\draw (\size, 0) rectangle (2*\size, \size);
	\draw (2*\size, 0) rectangle (3*\size, \size);
\end{tikzpicture}
}

\def\stickiv{\begin{tikzpicture}
    \foreach \y in {0,1,2,3} { \draw (0, \y*\size) rectangle (\size, \y*\size + \size); } 
\end{tikzpicture}
}
\def\stickv{\begin{tikzpicture}
    \foreach \y in {0,1,2,3,4} { \draw (0, \y*\size) rectangle (\size, \y*\size + \size); } 
\end{tikzpicture}
}
\def\stickvi{\begin{tikzpicture}
    \foreach \y in {0,1,2,3,4,5} { \draw (0, \y*\size) rectangle (\size, \y*\size + \size); } 
\end{tikzpicture}
}

\begin{table}[ht]
\centering
\begin{tabular}{ccccccccc}
\toprule
\stickiv & \stickv & \stickvi & \lv & \liv & \livii & \liii & \liiiii & $A_0 - \textrm{correction} = A$ \\ \midrule
3& 2& 1& 3& 11& 6& 4& 5& $1634 - 37 = 1597$\\
3& 2& 1& 3& 10& 6& 4& 6& $1627 - 35 = 1592$\\
4& 2& 1& 3& 11& 6& 4& 4& $1628 - 37 = 1591$\\
4& 2& 1& 3& 10& 6& 4& 5& $1622 - 35 = 1587$\\
3& 2& 1& 3& 9& 6& 4& 7& $1619 - 35 = 1584$\\
4& 2& 1& 3& 9& 6& 4& 6& $1616 - 35 = 1581$\\
2& 2& 1& 3& 11& 6& 5& 5& $1617 - 36 = 1581$\\
3& 2& 1& 3& 11& 6& 5& 4& $1614 - 36 = 1578$\\
2& 2& 1& 3& 10& 6& 5& 6& $1611 - 34 = 1577$\\
3& 2& 1& 3& 10& 6& 5& 5& $1608 - 34 = 1574$\\
2& 2& 1& 3& 9& 6& 5& 7& $1603 - 34 = 1569$\\
3& 2& 1& 3& 9& 6& 5& 6& $1602 - 34 = 1568$\\
1& 2& 1& 3& 11& 6& 6& 5& $1600 - 35 = 1565$\\
2& 2& 1& 3& 11& 6& 6& 4& $1598 - 35 = 1563$\\
1& 2& 1& 3& 10& 6& 6& 6& $1594 - 33 = 1561$\\
2& 2& 1& 3& 10& 6& 6& 5& $1593 - 33 = 1560$\\
2& 2& 1& 3& 9& 6& 6& 6& $1586 - 33 = 1553$\\
1& 2& 1& 3& 9& 6& 6& 7& $1586 - 33 = 1553$\\
0& 2& 1& 3& 11& 6& 7& 5& $1583 - 36 = 1547$\\
1& 2& 1& 3& 11& 6& 7& 4& $1582 - 36 = 1546$\\
0& 2& 1& 3& 10& 6& 7& 6& $1577 - 34 = 1543$\\
1& 2& 1& 3& 10& 6& 7& 5& $1576 - 34 = 1542$\\
1& 2& 1& 3& 9& 6& 7& 6& $1570 - 34 = 1536$\\
0& 2& 1& 3& 9& 6& 7& 7& $1569 - 34 = 1535$\\
\bottomrule
\end{tabular}
\caption{\label{tab:lreduction} Sensible combinations of reducing the 35 hexominoes to different straight and L-shaped pieces, with corrections for re-entrant cells. The first eight columns give how many pieces of each type occur in the replacement problem. $A_0$ is the optimum of the replacement problem, found with Gurobi, and $A = A_0 - \textrm{correction}$ is the resulting upper bound.}
\end{table}

For each combination we keep the smallest correction, which gives the largest bound and is therefore safe. These $24$ simplified problems can all be solved exactly with the commercial ILP solver Gurobi, each within a few minutes. The results are listed in Table~\ref{tab:lreduction}. The maximum, $1597$, is attained in the first row. Hence no hexomino fence can enclose more than $1597$ cells, and the configuration shown in Figure~\ref{fig:1597} is optimal.

It is important to mention that this approach assumes that all $35$
hexominoes are only used in a ``sensible'' configuration, so that
their contribution is actually captured by the one or two assigned
L-shaped pieces. Otherwise, the number of combinations of possible
L-shaped pieces would be infeasibly large. In contrast, the complete
search described in the previous subsections allows each hexomino to
be used with any pair of incoming and outgoing attachment edges and
thus provides a completely rigorous bound.

\section{Heptomino and octomino fences}\label{sec:hepto-okto}

\subsection{Lower bound using pruned dynamic programming}\label{sec:hepto-okto-lower}

The target-area test of Section~\ref{sec:target-area}
is relatively expensive: It
requires iteratively solving a nonlinear equation to find
the radius $R$ of the circumcircle.
Therefore, we apply this test once per generation, and not for each
insertion.
During the $i$-th generation, an entry in $\mathrm{OPT}_i$ may receive
several insertions, and its value may be improved several times.
After the whole table $\mathrm{OPT}_i$ has been filled,
we apply the test to every entry and keep only the entries that pass.

However, filling the complete table
without any pruning
becomes infeasible for heptominoes and octominoes.
Therefore, we apply some extra heuristic
tests.
      We expect the four sum vectors
 $x_{\mathrm{NE}}, x_{\mathrm{NW}},\allowbreak
 x_{\mathrm{SW}}, \allowbreak x_{\mathrm{SE}}$ to be balanced:
 They should be directed roughly in their respective $\pm45^\circ$ direction, and they
 should approximately have the same length.
Therefore, at each table insertion, we test the following two conditions.
\begin{enumerate}
\item \emph{Balance.}
    After 
    folding all four sum vectors
    into the
    first quadrant,
    we discard an entry if any two vectors have 
    $L^1$-distance more than a fixed threshold $W$, betting that an
    entry with highly
    unbalanced vectors 
    is unlikely to be completed to a good polygon.
\item \emph{Reach.}
  In addition, we look at the sum $\binom{s_x}{s_y}
  = x_{\mathrm{NE}} + x_{\mathrm{NW}}
    + x_{\mathrm{SW}} + x_{\mathrm{SE}}$, the vector still missing
    to close the polygon.
    We discard an entry if $\max(|s_x|, |s_y|)$ exceeds a second threshold $W'$.
\end{enumerate}

\begin{figure}[tbp]
  \centering
  \includegraphics[width=\textwidth]{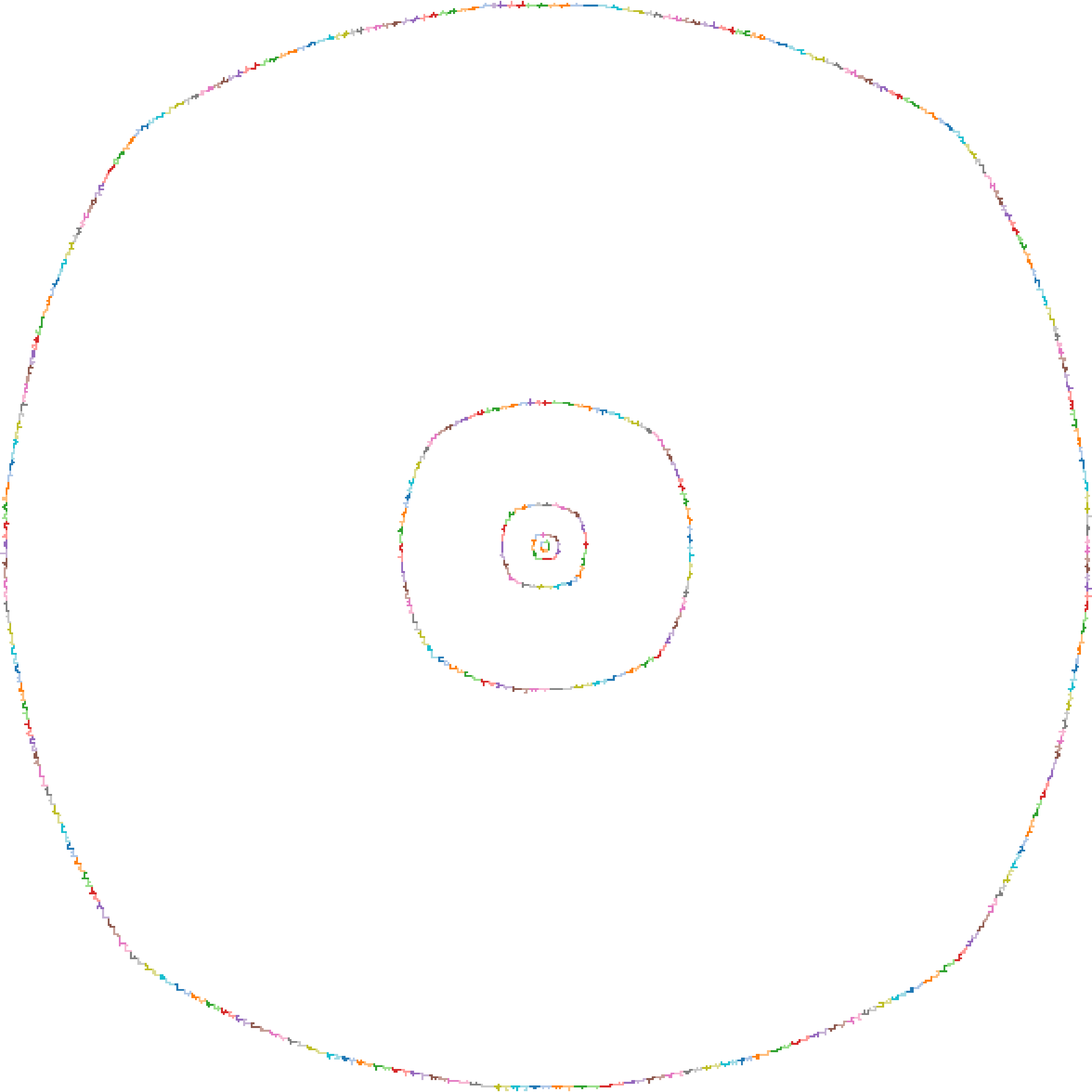}
  \caption{Our best solutions for the heptomino and octomino fences, alongside an optimal tetromino, pentomino, and hexomino fence.}
  \label{fig:fence-art}
\end{figure}

By trying various values for $W$ and $W'$ for a given $n$-omino fence, the dynamic programming search is able to find very good heuristic solutions in a reasonable amount of time. Figure~\ref{fig:fence-art} shows a heptomino fence with area $19\,733$ which we found in one hour and an octomino fence with area $279\,134$ that took seven hours of runtime.

Note that there is one octomino (a $3 \times 3$ square with a missing center) whose center hole is not 8-connected to the outside. By our definition in Section~\ref{sec:def-area}, this enclosed cell forms a separate region, which does not count towards the area. 

\subsection{Upper bounds with Branch-and-Bound}
\label{sec:ub-hex}

For an upper bound we can simply reuse the Branch-and-Bound method of
Section~\ref{sec:bnb}. If Step~1 is run with a bound
of $c$ instead of $1597.999$ and
returns no feasible candidate sequence,
we know that no fence can achieve an area $\geq c$.

This method is already quite powerful for the hexominoes: $c = 1601.999$
is the smallest value for which all possibilities can be ruled out by
Step~1. Thus, an upper bound of $1601$ for the hexomino fence, which
exceeds the true optimum by only 4, can be proved in around ten minutes (measured on an AMD Ryzen 7 7735HS CPU and 20\,GB of DDR5 memory).

On the same hardware, we can prove an upper bound of $19\,830$ for
heptominoes and $280\,389$ for octominoes in around 2 days of runtime
each. We compare them to the lower bounds from
Section~\ref{sec:hepto-okto-lower} in Table~\ref{tab:hepto-octo}, and
can observe that the relative difference between the bounds is less
than $0.5\%$. Note that the upper bounds are just the best possible
ones within that time limit, and are likely not the best bounds
achievable by this method. We stopped since the runtime of the B\&B
appears to grow significantly even for small decreases in~$c$;
for example, proving the heptomino bound $19\,850$ takes just one hour.

\begin{table}[htb]
  \centering
  \begin{tabular}{|r@{ }l|c|c|c|}
    \hline
    \multicolumn{2}{|l|}{$n$} & best DP solution & B\&B upper bound & relative difference\\
    \hline
     7 & & \phantom019\,733 & \phantom019\,830 & 0.49\% \\
     8 & & 279\,134 & 280\,389 & 0.45\% \\
     \hline
  \end{tabular}
  \caption{The best lower and upper bounds that we have found for the
    heptomino and octomino fences. The 
     relative
    difference
is taken with respect to
the lower bound.}
  \label{tab:hepto-octo}
\end{table}

\section{Larger fences}
\label{sec:larger-fences}

We now investigate $n$-omino fences for $n \geq 9$, where our dynamic programming approach is infeasible. We denote the number of free $n$-ominoes by $A^{\prime}(n)$; its exact values are known up to $n = 59$~\cite{OEIS_A000105}. For larger $n$, the following
estimate exists:
\begin{equation} \label{eq:an-asymptot}
    A^{\prime}(n) \approx \frac{1}{8} \frac{c \cdot \lambda^n}{n}
\end{equation}
where $\lambda = 4.0626$ and $c = 0.3169$~\cite{Iwan_Jensen_2000}. The original formulation
in~\cite{Iwan_Jensen_2000}
estimates the number of fixed $n$-ominoes (i.e.,\ rotated and reflected copies are counted as
distinct), but since almost no large $n$-ominoes have any symmetries,
we can just divide that estimate by $8$. Note that this estimate has a significant uncertainty: it is only proven that $4.0025 \leq \lambda \leq 4.5238$~\cite{poly-lambda,bui2026}.

\subsection[Fmax + Dmax upper bound]{$F_{\max} + \Delta_{\max}$ upper bound}
\label{sec:fmax}

We relax the bound $S_1$ from \eqref{eq:S1} by optimizing over the
choice of $\bar 
v_i$ and $\Delta_i$ separately:
\begin{align}\nonumber
\text{Area} &\le S_1
= \max \Bigl\{
\bar F\bigl(\bar v_1^{a_1},\bar v_2^{a_2},\ldots,\bar v_{A^{\prime}(n)}^{a_{A^{\prime}(n)}}\bigr)
+ \Delta_1^{a_1}+\Delta_2^{a_2}+\dots+\Delta_{A^{\prime}(n)}^{a_{A^{\prime}(n)}}
\Bigm|
1 \le a_i \le m_i\,
\Bigr\} \\
            &\le 
              \max \Bigl\{
\bar F\bigl(\bar v_1^{a_1},\bar v_2^{a_2},\ldots,\bar v_{A^{\prime}(n)}^{a_{A^{\prime}(n)}}\bigr)
\Bigm|
1 \le a_i \le m_i
              \Bigr\}
              + \max \Bigl\{
         \Delta_1^{a_1}+
         \dots+\Delta_{A^{\prime}(n)}^{a_{A^{\prime}(n)}}
\Bigm|
1 \le a_i \le m_i
\Bigr\} \label{eq:Dmax}
\end{align}
Recall that
$\bar F(\bar v_1, \bar v_2, \ldots, \bar v_{A^{\prime}(n)})$
is the maximum area of a convex $A^{\prime}(n)$-gon
obtained by rotating and reflecting the vectors $\bar v_1, \bar v_2, \ldots, \bar
v_{A^{\prime}(n)}$ and
 concatenating them
in a suitable
order. For this
we use as an easy upper bound the area of the largest
$A^{\prime}(n)$-gon with a given perimeter.
It is well known that
the $k$-gon of maximum area with given perimeter $u$ 
is the regular $k$-gon
\cite[Theorem 12.4a]{niven1981maxima}, and its area is
\[
G_k(u) = \frac{u^2}{4  k}\cot\frac{\pi}{k}.
\]
Thus,
\begin{displaymath}
\bar F(\bar v_1, \bar v_2, \ldots,\bar  v_{A^{\prime}(n)})
\le
G_{A^{\prime}(n)}( |\bar v_1| + |\bar v_2| + \dots + |\bar v_{A^{\prime}(n)}|),
\end{displaymath}
and we can bound the first term in \eqref{eq:Dmax}
by 
\begin{displaymath}
 F_{\max} :=  
G_{A^{\prime}(n)}\bigl(|\bar v_1^{\max}| + |\bar v_2^{\max}| + \dots + |\bar v_{A^{\prime}(n)}^{\max}|\bigr)
\end{displaymath}
where $\bar v_i^{\max}$ denotes a longest vector among
$\bar v_i^{1},\ldots,\bar v_i^{m_i}$.
We denote the maximum
in the second expression in \eqref{eq:Dmax}
by $\Delta_{\max}$. (This is the maximum without the rotation and parity conditions, in contrast to $\Delta^{\mathrm{rot}}$ of Step~2 in Section~\ref{sec:rotations}.) Note that the sets $\mathcal{W}_i$ contain the null option $(\bar v,\Delta)=(0,0)$ to allow a piece not to be used. This means that the maximum is $\Delta_{\max}=\sum_{i=1}^{A^{\prime}(n)}\max\{0,\Delta_i^{\mathrm{best}}\}$, where $\Delta_i^{\mathrm{best}}$ is the largest correction of piece~$i$ over all its proper (non-null) reference vectors. The null option does not affect $F_{\max}$, because leaving out a piece decreases the perimeter and thus can only lower the bound.
 Since the area is an integer, we have proven the upper bound
\[
\text{Area} \leq \lfloor F_{\max} + \Delta_{\max} \rfloor.
\]
This is easy to compute. For each $n$-omino we only need to determine
the longest reference vector and the maximal correction
separately. The bounds up to $n = 15$
are given in Table~\ref{tab:heur}. For $n=7$ and~$8$, this bound is
considerably weaker than the more expensive bound
of \autoref{tab:hepto-octo}.

\subsection{Improved upper bound by linearization}
\label{sec:linearization}

The bound $U(n)$ takes the longest reference vector and the largest
correction of each piece independently of each other, although they
usually come from different pairs of attachment edges. The best
correction that a piece can achieve together with a longest reference
vector is typically much smaller than its overall largest correction.
For the hexominoes, for example, $\Delta_{\max}=27$, but if every piece
uses a longest reference vector, the sum of the corrections is at
most $-73$.
We now show how the two terms can be coupled. The idea is to replace
the quadratic function $G_k(u)$ by a linear upper bound in the relevant
range of~$u$. Then the objective becomes a sum of separate terms for
the individual pieces, which can again be maximized for each piece
separately.

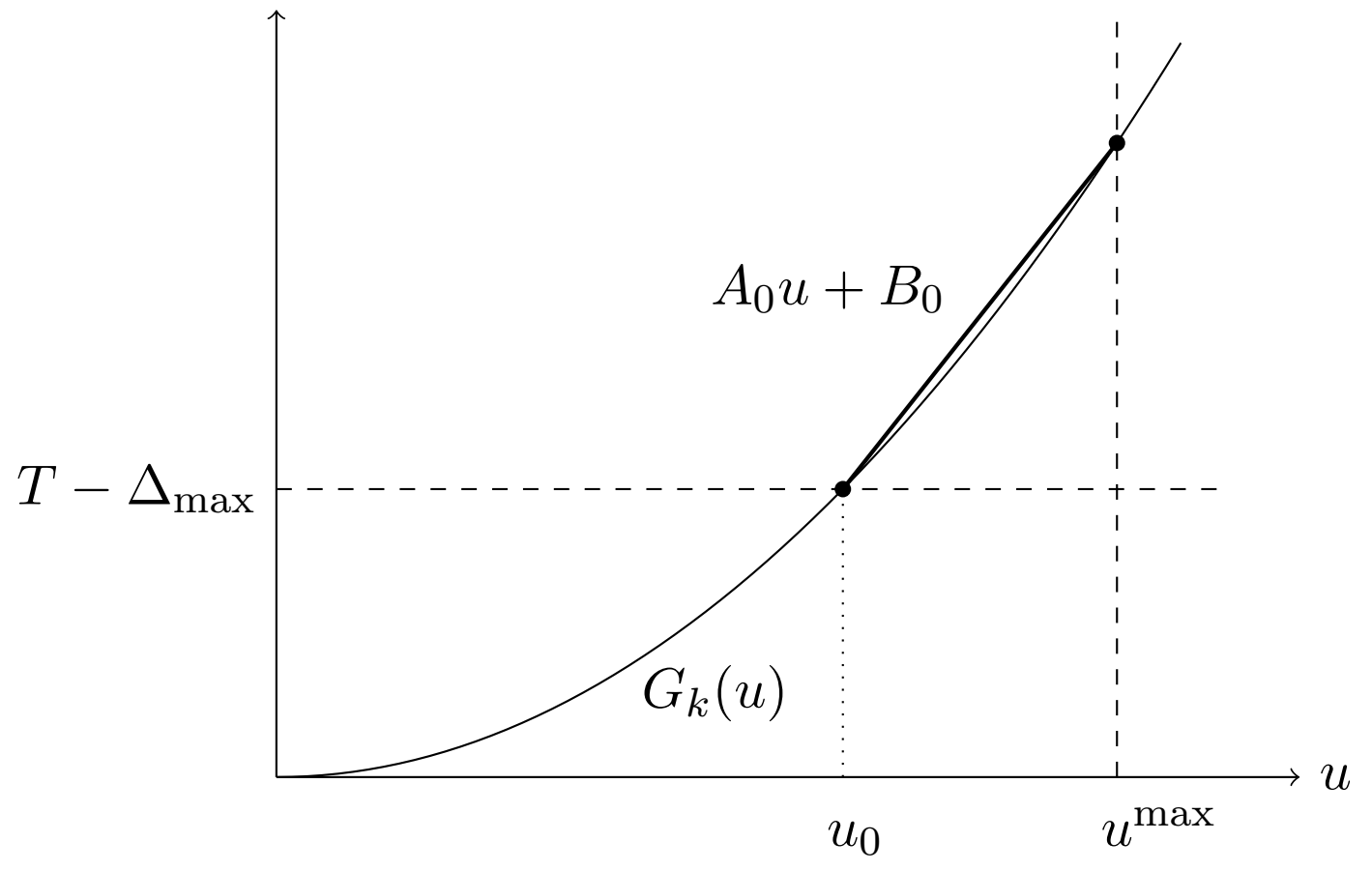
\begin{figure}[htb]
  \centering
  \begin{tikzpicture}[x=1.25cm,y=0.82cm]
    \def\uz{3.1}\def\um{4.6}
    \draw[->] (0,0) -- (5.6,0) node[right] {$u$};
    \draw[->] (0,0) -- (0,6.4);
    \draw[domain=0:4.95,samples=60] plot (\x,{0.25*\x*\x});
    \node at (2.4,0.7) {$G_k(u)$};
    \draw[dashed] (0,{0.25*\uz*\uz}) node[left] {$T-\Delta_{\max}$} -- (5.2,{0.25*\uz*\uz});
    \draw[dashed] (\um,0) node[below] {$\strut u\rlap{$^{\max}$}$} -- (\um,6.3);
    \draw[dotted] (\uz,0) node[below] {$\strut u\rlap{$_0$}$} -- (\uz,{0.25*\uz*\uz});
    \draw[thick] (\uz,{0.25*\uz*\uz}) -- (\um,{0.25*\um*\um})
      node[pos=0.45,above left] {$A_0u+B_0$};
    \fill (\uz,{0.25*\uz*\uz}) circle (1.6pt) (\um,{0.25*\um*\um}) circle (1.6pt);
  \end{tikzpicture}
  \caption{The linear upper bound $A_0u+B_0$ for the convex function
    $G_k(u)$ in the range $u_0\le u\le u^{\max}$ (Lemma~\ref{lem:linear}).}
  \label{fig:linear}
\end{figure}

Let $k=A'(n)$ and
$u^{\max} := |\bar v_1^{\max}| + |\bar v_2^{\max}| + \dots + |\bar
v_{k}^{\max}|$, so that $F_{\max}=G_k(u^{\max})$.
We choose a threshold~$T$ with
$\Delta_{\max} \le T\le F_{\max}+\Delta_{\max}$,
define $u_0\le u^{\max}$ by the equation $G_k(u_0)+\Delta_{\max}=T$,
and let $A_0u+B_0$ be the linear function that agrees with $G_k(u)$ at
$u=u_0$ and at $u=u^{\max}$, see \autoref{fig:linear}.

\begin{lemma}\label{lem:linear}
  For $0\le u\le u^{\max}$,
  \begin{equation*}
    G_k(u) \le \max\{\,T-\Delta_{\max},\ A_0u+B_0\,\}.
  \end{equation*}
\end{lemma}
\begin{proof}

  The function $G_k(u)$ is increasing and convex for $u\ge0$.
  For $u\le u_0$, we have
  $G_k(u)\le G_k(u_0)=T-\Delta_{\max}$.
  For $u_0\le u\le u^{\max}$, the convex function $G_k$ lies below its chord
  between $u_0$ and~$u^{\max}$, which is the line $A_0u+B_0$.
\end{proof}

Now consider an arbitrary choice $(\bar v_i^{a_i},\Delta_i^{a_i})$,
$1\le a_i\le m_i$, of a reference vector class and a correction for every
piece, and let
$u=\sum_{i=1}^k|\bar v_i^{a_i}|\le u^{\max}$ and
$\Delta=\sum_{i=1}^k\Delta_i^{a_i}\le\Delta_{\max}$.
By the bound $\bar F\le G_k$ of Section~\ref{sec:fmax}
and Lemma~\ref{lem:linear},
\begin{align*}
  \bar F\bigl(\bar v_1^{a_1},\ldots,\bar v_{k}^{a_{k}}\bigr)+\Delta
  &\le G_k(u)+\Delta
  \le \max\{\,T-\Delta_{\max}+\Delta,\ A_0u+B_0+\Delta\,\}
  \\&\le \max\Bigl\{\,T,\ \sum_{i=1}^{k}\bigl(A_0|\bar v_i^{a_i}|+\Delta_i^{a_i}\bigr)+B_0\Bigr\}.
\end{align*}
Taking the maximum over all choices, we obtain
\begin{equation}\label{eq:S3}
  \text{Area}\le S_1 \le \max\{T,S_3\}
  \quad\text{with}\quad
  S_3 := \sum_{i=1}^{k}\max\bigl\{\,A_0|\bar v_i^{a}|+\Delta_i^{a}
  \bigm| 1\le a\le m_i\,\bigr\}+B_0 .
\end{equation}
Hence, if $S_3\ge T$, we have proved the upper bound
$\text{Area}\le\lfloor S_3\rfloor$.
Like $U(n)$, the bound $S_3$ can be computed for each piece separately:
We try all reference vectors $\bar v_i^a$ of piece~$i$ and determine the
maximum of the weighted sum $A_0|\bar v_i^{a}|+\Delta_i^{a}$ of length and
correction.

The coefficient $A_0$ is the \emph{price of length}: shortening a reference
vector by one unit reduces the bound by about $A_0$ units of
area. For the hexominoes, with $T=1597.999$, we have $A_0\approx 22.80$
and $B_0\approx-1636.44$.  The line $A_0u+B_0$ exceeds $G_{35}(u)$ by
less than $0.7$ in the middle of the interval $[u_0,u^{\max}]$.
For larger $n$, $A_0$ grows rapidly, up to $A_0\approx 4.29\cdot10^6$ for
$n=15$, whereas the corrections of a single piece vary only by small
amounts.
In other words, a tiny improvement in length will more than compensate
any variation in penalty.
At least
in terms of the upper bound, 
only \emph{length matters}.
The correction is only a secondary criterion, to decide between
alternative reference vectors of the same largest length.

Consequently, for every $n$ from $6$ to~$15$ (with $T=1597.999$ for
$n=6$ and $T=L(n)+0.999$ for $n\ge7$) and for every piece, the
maximum of $A_0|\bar v_i^{a}|+\Delta_i^{a}$ is attained by a longest
reference vector $\bar v_i^{\max}$. For these vectors, the lemma is
applied at $u=u^{\max}$, where the upper bound is tight, and therefore
\begin{equation*}
  S_3 = A_0u^{\max}+B_0+\Delta^{\mathrm{top}}
      = F_{\max}+\Delta^{\mathrm{top}},
\end{equation*}
where $\Delta^{\mathrm{top}}$ is the sum, over all pieces, of the
largest correction that the piece achieves with a longest reference
vector.
In particular, $S_3$ does not change if the chosen threshold~$T$ is increased
(this increases~$A_0$), and the
improvement over $U(n)$ is
$\Delta_{\max}-\Delta^{\mathrm{top}}$, which is a multiple
of~$\frac18$. For the hexominoes, it is $27-(-73)=100$, and we get
$S_3\approx1631.61$, so the area is at most $1631$, compared to $U(6)=1731$.
(This is even a little better than the initial value $1633.964$ of
the bound~\eqref{eq:wide-bound} of Section~\ref{sec:minkowski},
although it is weaker than the bound $1601$ of
Section~\ref{sec:ub-hex}.)

We denote the improved upper bound by $U'(n):=\lfloor S_3\rfloor$.
In other words, for $6\le n\le 15$,
\begin{equation*}
  U'(n) = \lfloor F_{\max}+\Delta^{\mathrm{top}}\rfloor.
\end{equation*}
It is the bound that one would get by using every piece with its
longest reference vector and the best correction for this vector,
as in the greedy heuristic of the next section.
The point of the linearization is to justify this: Without it, a
piece might trade a little length for a larger correction.
For the pentominoes, this actually happens: With $T=128.999$, four of
the pentominoes achieve the maximum in~\eqref{eq:S3} with a shorter
vector, and $S_3\approx131.27$ gives the bound $131$. The value
$F_{\max}+\Delta^{\mathrm{top}}\approx128.06$ is not justified by
this argument. (The optimum is $128$.)
It would not be hard to prove that the ``longest-first'' rule,
which we observed to hold for $6\le n\le15$, holds
generally for larger~$n$, given that the gap between vectors of
different lengths (square roots of integers) cannot be very small, the
variation
of correction terms is bounded, and the fact that the number of pieces
and hence the slope $A_0$ grows exponentially with $n$.

The values of $U'(n)$ for $n\le15$ are shown in Table~\ref{tab:heur}.
For $n=7$ and~$8$, the gap to the lower bound shrinks considerably,
from $3.38\%$ to $1.52\%$ and from $1.81\%$ to $1.28\%$. For larger~$n$,
the improvement, although large in absolute terms, becomes
negligible in relative terms, because the corrections are negligible
compared to $F_{\max}$: For $n=12$ and $13$, the relative gap changes only in
the third decimal, and for larger $n$ is it smaller still. The remaining gap comes almost entirely from
the difference between the reference polygon and a regular
$A'(n)$-gon.
(The bound might be tightened by using the largest polygon area with
given individual edge lengths
$|\bar v_i^{\max}|$, as in
\autoref{sec:target-area}, instead of just using the perimeter
$u^{\max}=
\sum_i|\bar v_i^{\max}|$.)
Our program \texttt{upper-bounds.py}~\cite{github} computes $S_3$ in
addition to $U(n)$, using $T=L(n)+0.999$.

\subsection{Heuristic greedy solver}

For the lower bound, we heuristically attempt to construct an
$A^{\prime}(n)$-gon that is as regular as possible.
As a start, we simply choose the longest reference vector for each
$n$-omino.
When there is a tie, we use the correction as a secondary criterion.
 We then try to solve the
\maxpolygon\ problem of Section~\ref{sec:8way} for these vectors.
A solution can only exist if the parity condition
\eqref{eq:parity-sum} is satisfied. Therefore, if necessary,
we flip the parity by changing the reference vector of the longest L-piece so that it is used
like a stick of $n-1$ cells.

The four vectors closest to the $45^\circ$ direction are
selected beforehand as the polygon's four corners.
More precisely, these are $n$-ominoes in which the vector
between two diagonally opposite corners is $\binom{\lceil
  n/2\rceil}{\lceil n/2\rceil}$ (and thus their longest reference vector would be $\binom{\lceil n/2\rceil - 1/2}{\lceil n/2\rceil - 1/2}$). For $n\ge 6$, there are always
sufficiently many pieces of this kind.
(This changes four reference vectors, but leaves the parity condition
intact.)
In contrast to the dynamic programming approach
 of Section~\ref{sec:8way}, we process
 the remaining
vectors
in a greedy way, in decreasing order of length:
Each vector
is assigned 
to whichever of the eight arcs keeps the partial vector sums of the arcs
as balanced as possible, breaking ties in favor of keeping the running
total closest to zero.
 Most of the time, the greedy balancing
 will not result in an exact total of $(0, 0)$, i.e.,\ a closed loop.
We therefore make a small number of 
corrective swaps at the end: a residual of at least $2$ in one coordinate is removed by exchanging
two vectors that have the same other coordinate and whose coordinates
in the residual direction add up to $\pm1$: each of them is replaced by
the mirror image of the other. This does not change the lengths of the
vectors used, but reduces the residual by~$2$ in that coordinate.
A final residual of $(\pm1,\pm1)$ is eliminated by reflecting one
suitable vector at the diagonal (that is, by swapping the absolute
values of its coordinates).
All chosen vectors, corners included, are then
sorted once by slope into the cyclic order of the final $A'(n)$-gon.
The area is then computed 
via the Shoelace formula.

We remark that a parity discrepancy arises also for $n=5$ and $n=6$.
For pentominoes, it is the U-pentomino
that is treated exceptionally: it is not used with its longest reference
vector in the optimal fences (\autoref{fig:use-opt}).
For hexominoes, it is the T-shaped piece discussed in
Section~\ref{sec:overlap}.

A Python implementation can be found in \texttt{best-fence-greedy.py}~\cite{github}. Finding the greedy solutions itself is very quick, but generating all $n$-ominoes and determining their longest reference vectors can take a while for larger $n$. After a solution is constructed, our implementation validates that all $n$-ominoes are placed exactly once (can be skipped) and without overlap. Additionally, we verify, either with the Shoelace formula applied to the edge boundary of Section~\ref{sec:boundary} or with an interval sweep, that every region other than the largest one is a hole inside a single $n$-omino.

Note that we could also extend Lemma~\ref{non-overlap} to reference vectors with small angles (such as the angles that arise in the greedy construction) and then skip the overlap check entirely. However, since its runtime is negligible even for $n = 15$, we do not include this argument here.

We computed heuristic solutions for $n \in \{7, \ldots, 15\}$, which, including the polyomino generation and all verifications, took around 90 minutes on an AMD Ryzen 7 7735HS CPU. The greedy construction itself takes just 11 seconds even for $n = 15$, once the longest reference vectors are calculated. Besides the 20\,GB of DDR5 memory, we needed additional swap space for $n = 14$ and~$15$, which probably slowed down the verification further.

The areas of these solutions, which directly give a lower bound $L(n)$ on the optimum, are shown in Table~\ref{tab:heur}. One can observe that the relative difference to the upper bound 
decreases clearly
for $n \geq 9$.

\begin{table}[htb]
  \centering
  \small
  \setlength{\tabcolsep}{4pt}
  \begin{tabular}{|r@{ }l|r|r|r|r|c|c|}
    \hline
    \multicolumn{2}{|l|}{$n$} & $A^{\prime}(n)$ & $U(n)$ & $U'(n)$ & $L(n)$ & $\dfrac{U(n) - L(n)}{L(n)}$ & $\dfrac{U'(n) - L(n)}{L(n)}$ \\
    \hline
     7 & &         108 &                 20\,385 &                 20\,018 &                 19\,718 &    3.3827\% &    1.5215\% \\
     8 & &         369 &                284\,170 &                282\,687 &                279\,105 &    1.8147\% &    1.2834\% \\
     9 & &        1\,285 &            4\,052\,974 &            4\,047\,043 &            4\,003\,811 &    1.2279\% &    1.0798\% \\
    10 & &        4\,655 &           61\,484\,424 &           61\,459\,838 &           60\,832\,848 &    1.0711\% &    1.0307\% \\
    11 & &       17\,073 &          942\,469\,010 &          942\,367\,722 &          934\,201\,230 &    0.8850\% &    0.8742\% \\
    12 & &       63\,600 &       14\,736\,083\,080 &       14\,735\,660\,662 &       14\,625\,233\,629 &    0.7579\% &    0.7550\% \\
    13 & &      238\,591 &      231\,340\,860\,898 &      231\,339\,111\,595 &      229\,830\,309\,620 &    0.6572\% &    0.6565\% \\
    14 & &      901\,971 &    3\,657\,193\,899\,930 &    3\,657\,186\,626\,293 &    3\,635\,789\,660\,265 &    0.5887\% &    0.5885\% \\
    15 & &    3\,426\,576 &   57\,965\,740\,117\,847 &   57\,965\,710\,006\,639 &   57\,656\,581\,776\,211 &    0.5362\% &    0.5362\% \\
    \hline
  \end{tabular}
  \caption{For the maximum area of each $n$-omino fence for $n \leq
    15$, we present the upper bound $U(n) := \lfloor F_{\max} +
    \Delta_{\max} \rfloor$ of Section~\ref{sec:fmax}
    and the improved upper bound
    $U'(n) := \lfloor S_3 \rfloor = \lfloor F_{\max} + \Delta^{\mathrm{top}} \rfloor$
    of Section~\ref{sec:linearization}.
    The areas of the heuristic solutions are themselves a lower
    bound $L(n)$, which we present alongside the relative gaps $(U(n)
    - L(n)) / L(n)$ and $(U'(n) - L(n)) / L(n)$. For comparison, we also include the bounds for $n
    \leq 8$, for which we know better bounds or even the exact
    values. We also have $U(5) = 154$, $U(6) = 1\,731$,
    $U'(5) = 131$, and $U'(6) = 1\,631$.}
  \label{tab:heur}
\end{table}

\begin{figure}[ht!]
    \centering
    \includegraphics[width=0.8\textwidth]{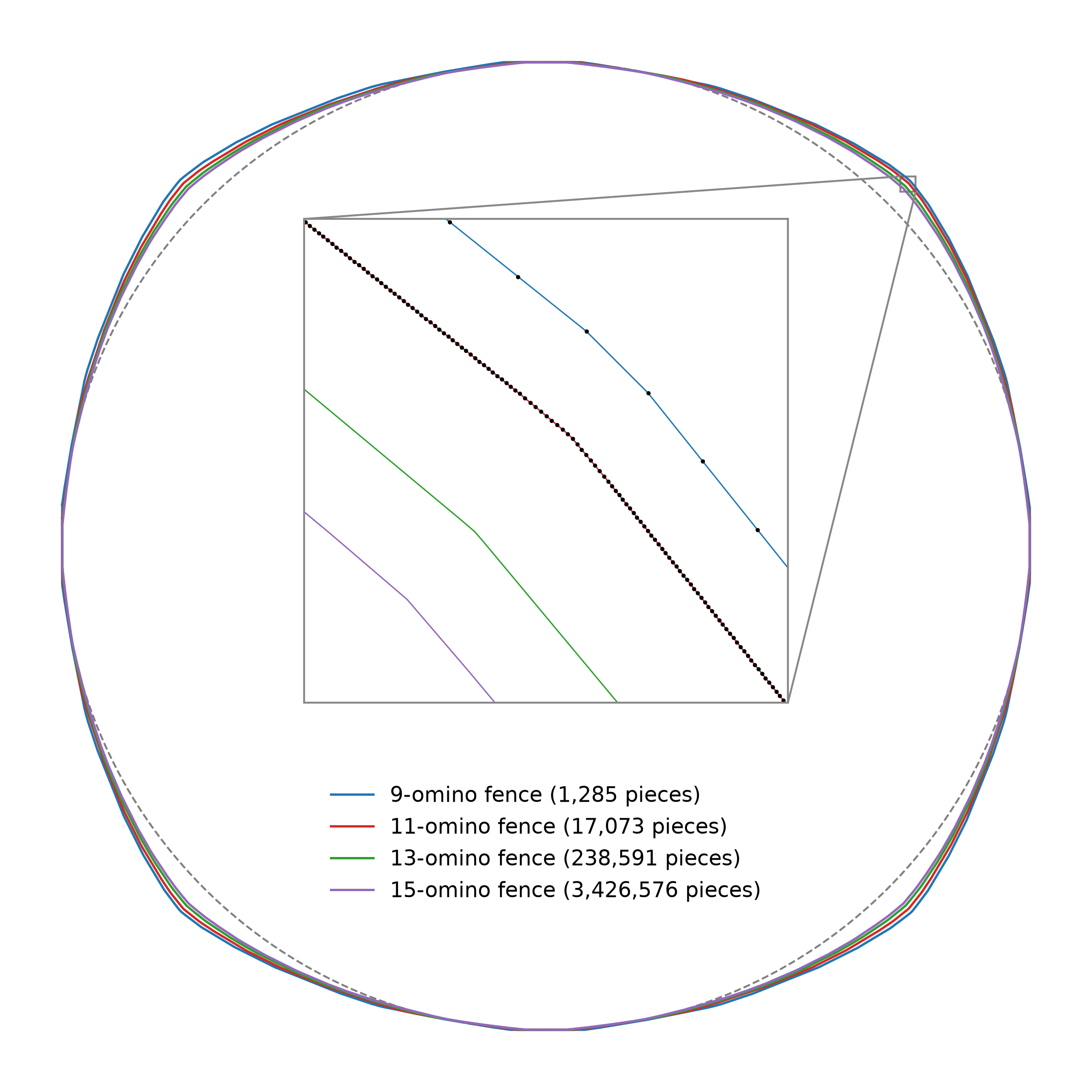}
    \caption{The reference polygons for $n = 9, 11, 13, 15$ found by
      the greedy approach. The polygons are scaled to fit in the same
      square box, alongside a gray reference circle. The enlarged inset shows the situation near the
      Northeast corner. Vertices
      for $n = 9, 11$ appear as small black dots.
      For larger~$n$, the vertices lie so densely on
 the edges that they could not be distinguished.}
    \label{fig:large-ref-polygons}
  \end{figure}

 Figure~\ref{fig:large-ref-polygons}.
shows the reference polygons for the greedy solutions for
 $n = 9, 11, 13, 15$.
A $n$ gets larger, the shape becomes
 rounder, but even with over $3$ million vectors at $n = 15$, the
 ``corners'' at the corner pieces 
 are still distinctly visible. Of course, our heuristic solver explicitly uses the eight arcs, but the same rounded octagon shape can be observed in the optimal hexomino fence and in the dynamic programming solutions for $n = 7$ and~$8$. Additionally, the upper bounds $U(n)$ and $U'(n)$ show that the true optimum cannot be much larger.

\subsection{The limit shape}
\label{sec:limit}

Looking at the results for $n \leq 15$ in Figure~\ref{fig:large-ref-polygons}, it remains unclear whether the four diagonal dents in the reference polygon will persist for larger $n$ or if the reference polygon will just become circular. To investigate this problem, we choose $n = 500$, which is a significant leap compared to the heuristic results but still leads to an acceptable runtime for the following sampling procedure. By formula~\eqref{eq:an-asymptot}, we estimate $A^\prime(500)$, the number of free $500$-ominoes, as $2.0 \cdot 10^{300}$. However, due to the uncertainty in the growth constant~$\lambda$ mentioned above, the value of $A^\prime(500)$ could be only around $10^{297}$ or roughly as high as $10^{323}$.

\begin{figure}[htb]
  \centering
  \raisebox{-.5\height}{\includegraphics[scale=0.062,angle=0]
    {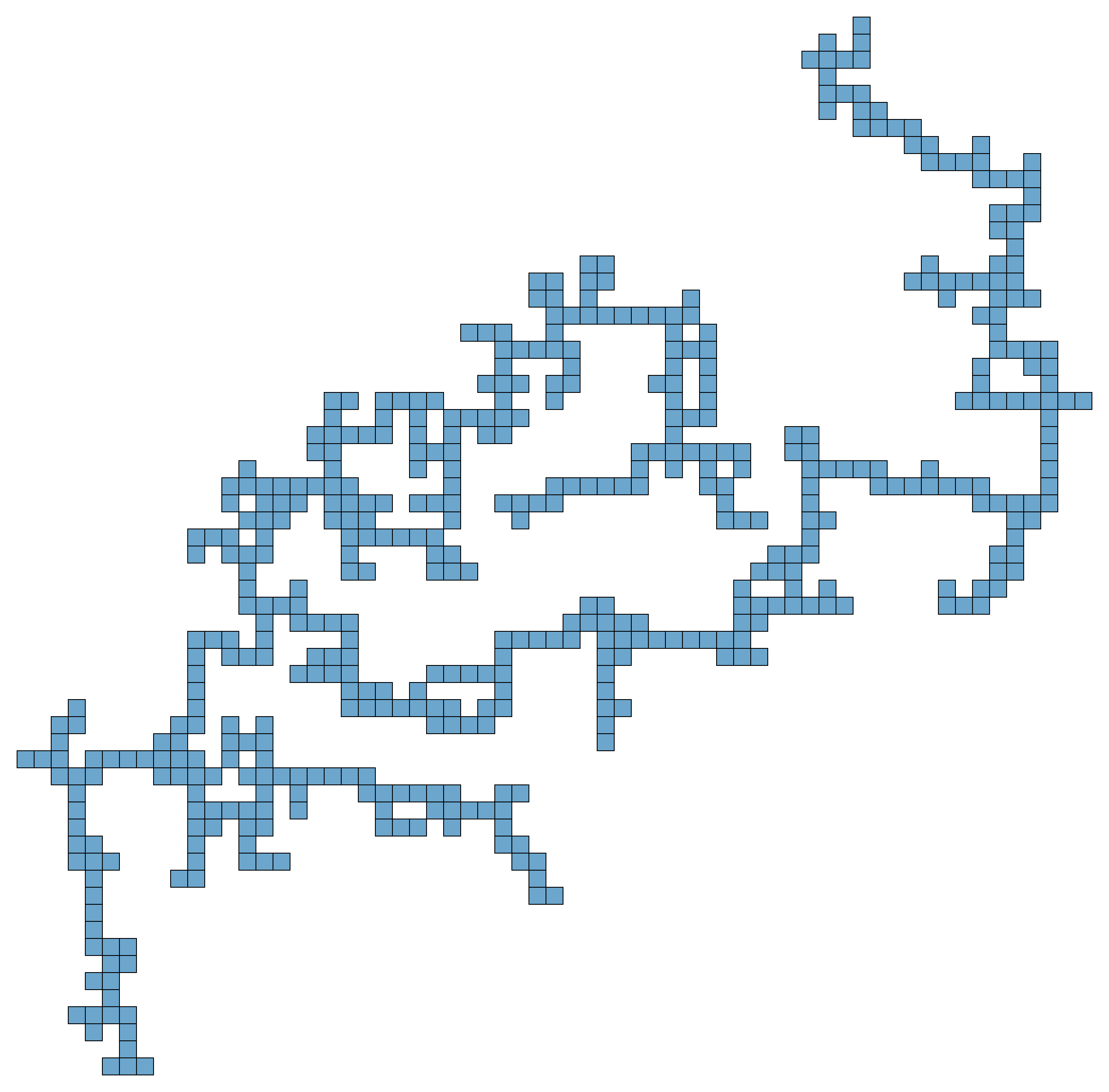}}
  \raisebox{-.5\height}{\includegraphics[scale=0.062]
    {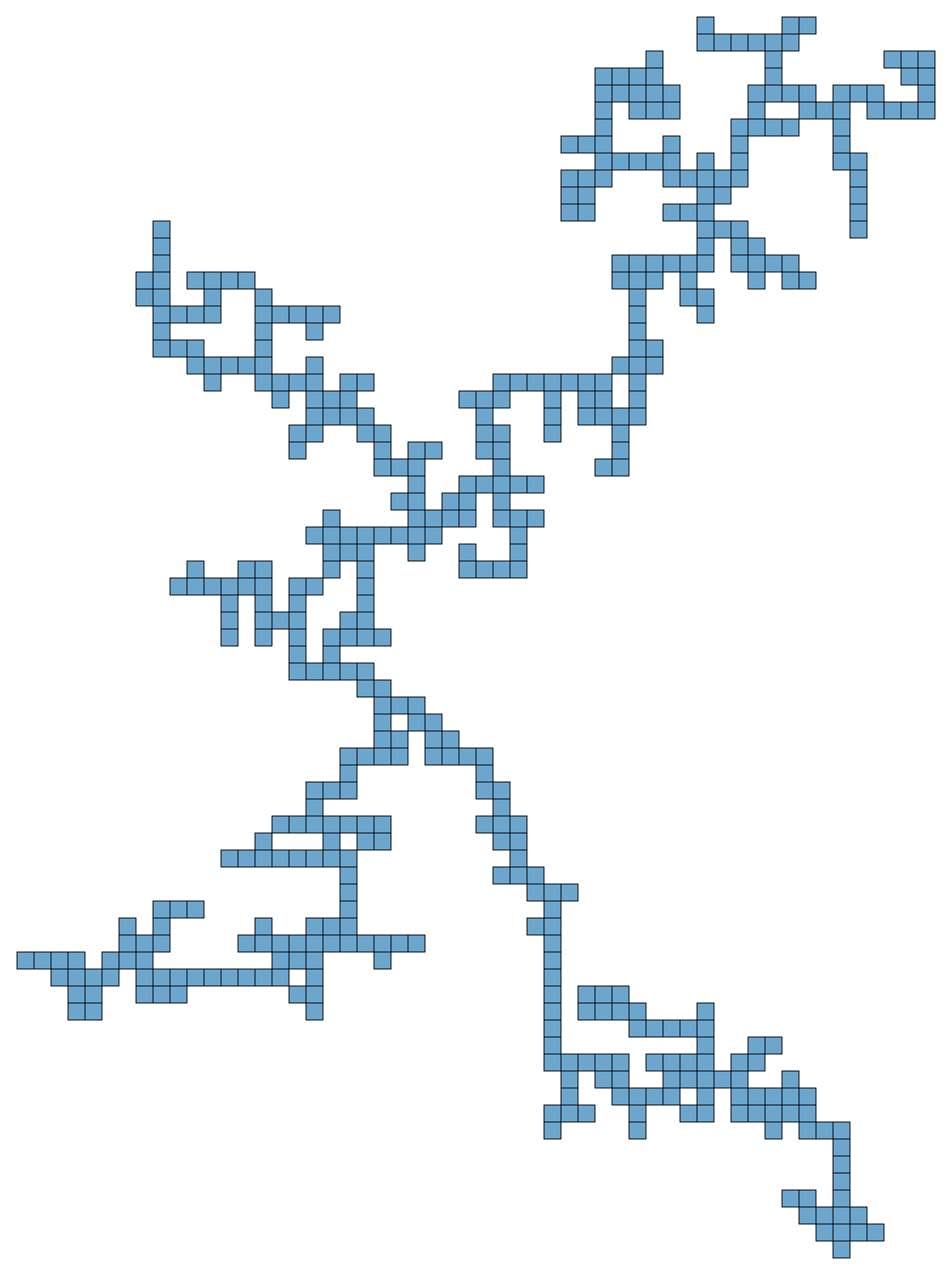}}
  \raisebox{-.5\height}{\includegraphics[scale=0.062,angle=90]
    {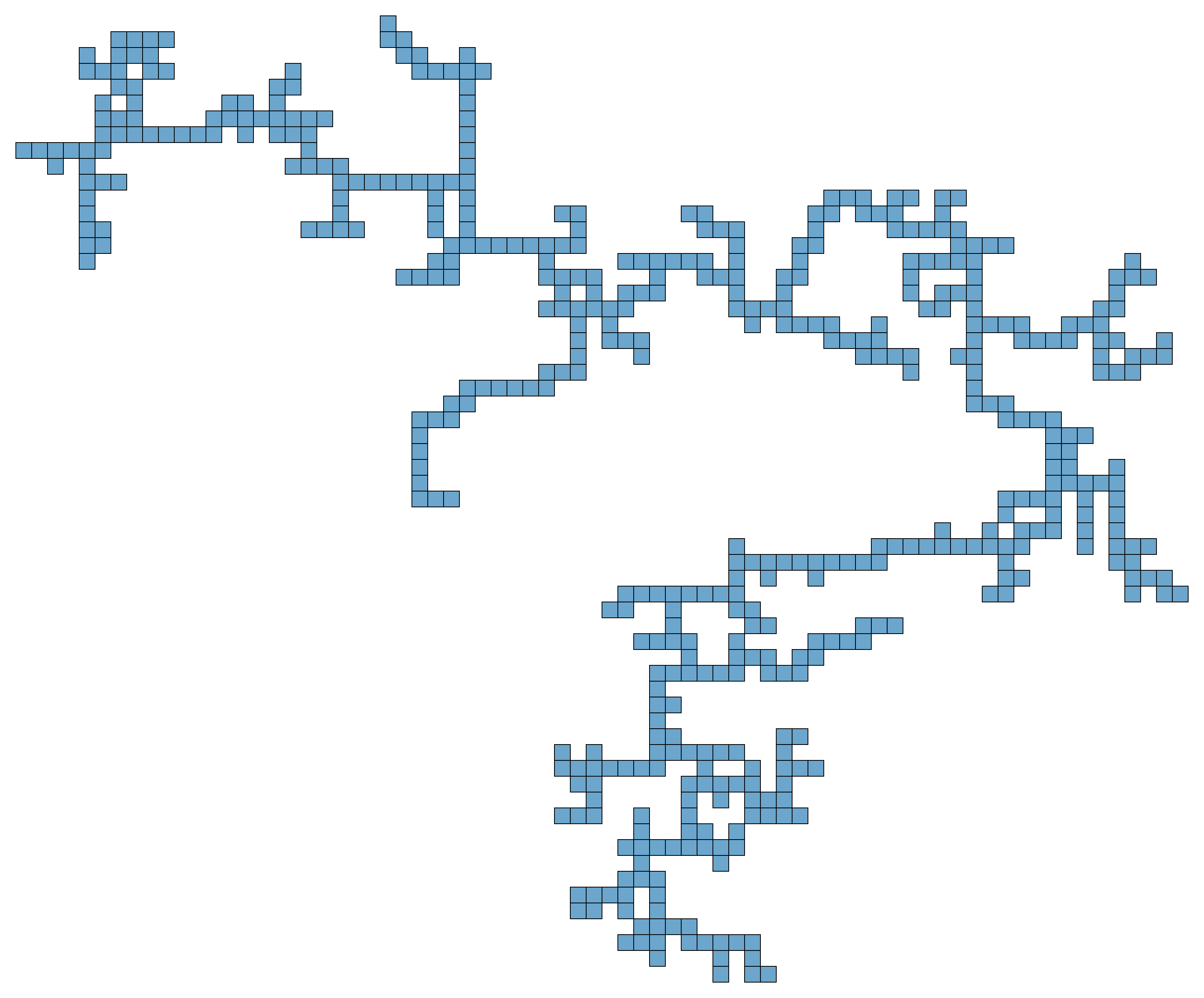}}

  \caption{Three random 500-ominoes}
  \label{fig:500}
\end{figure}

We use the Metropolis algorithm~\cite{hastings1970monte} to sample
random $500$-ominoes. A step of the chain removes a uniformly random cell~$x$ of the current polyomino~$A$ and adds a uniformly random cell~$y$ that is 4-adjacent to $A\setminus\{x\}$ (possibly $y=x$); if the result is connected, it becomes the new state, and otherwise the move is rejected and $A$ is kept.
For our experiments, we made $4n^2 =
10^6$ steps of the Markov chain for every sample; this number is a heuristic choice
for the mixing time. Our starting position for the first sample is the $1 \times n$ straight line, while for all later samples we simply choose the previous sample as the starting position.
 Figure~\ref{fig:500} shows three random $500$-ominoes.
Our algorithm samples
from the uniform distribution of
\emph{fixed} polyominoes, where rotated and reflected copies are distinct. Since the overwhelming majority of polyominoes have no symmetries, this is a good
approximation for the free polyominoes, which are used for our fences.

In this way, we generated $13.5$ million random $500$-ominoes,
using about 200 CPU days.
For each of them, we determine the longest reference vector in the straightforward way: among all midpoints of boundary edges of the $n$-omino, we take the pair with the largest Euclidean distance and use the vector connecting the two midpoints. Up to reflections and rotations by $90^\circ$, we normalize this vector to $(x,y)$ with $x \ge y \ge 0$. If multiple vectors share the same length, we record all of them.
By Lemma~\ref{lem:integer-vector},
 these vectors have integer
coordinates.

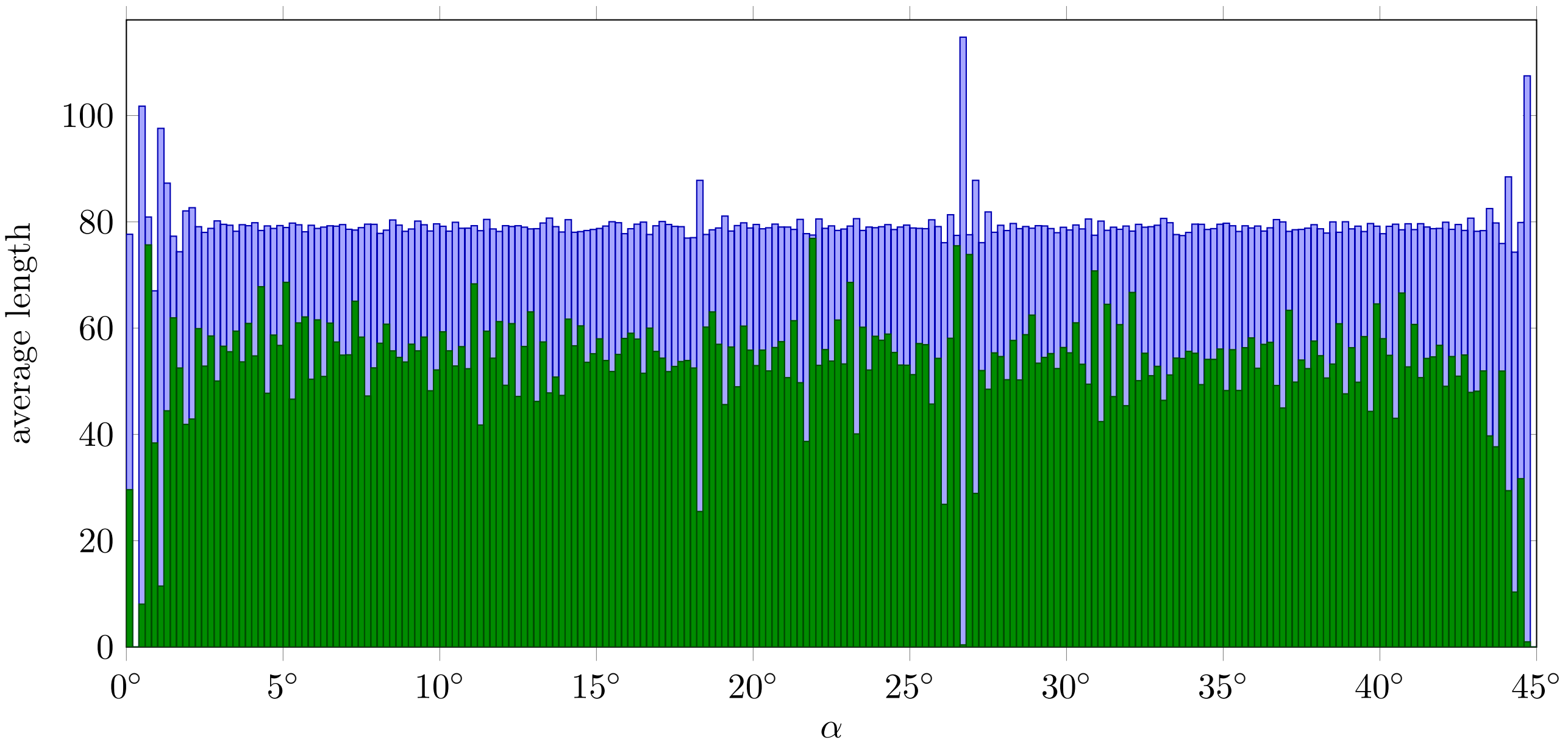
\begin{figure}[htb]
  \centering
\begin{tikzpicture}
\begin{axis}[scale only axis, width=0.9\linewidth, height=0.4\linewidth, ybar, bar shift=0pt, bar width=0.2, xmin=0, xmax=45, ymin=0, ymax=118, xtick={0,5,...,45}, xticklabel={$\pgfmathprintnumber{\tick}^\circ$}, scaled y ticks=false, xlabel={$\alpha$}, ylabel={average length}]
\addplot[fill=blue!35, draw=blue!70!black, thin] coordinates {
(0.1000,77.665)
(0.3000,0)
(0.5000,101.768)
(0.7000,80.91)
(0.9000,67.0201)
(1.1000,97.5901)
(1.3000,87.2834)
(1.5000,77.3)
(1.7000,74.3779)
(1.9000,82.0593)
(2.1000,82.6601)
(2.3000,79.086)
(2.5000,78.0146)
(2.7000,78.7711)
(2.9000,80.1806)
(3.1000,79.5335)
(3.3000,79.3639)
(3.5000,78.2333)
(3.7000,79.4646)
(3.9000,79.2749)
(4.1000,79.8341)
(4.3000,78.3746)
(4.5000,79.2374)
(4.7000,78.7625)
(4.9000,79.2997)
(5.1000,78.9371)
(5.3000,79.7628)
(5.5000,79.4471)
(5.7000,78.1337)
(5.9000,79.3562)
(6.1000,78.7612)
(6.3000,79.0133)
(6.5000,79.2562)
(6.7000,79.186)
(6.9000,79.4838)
(7.1000,78.6021)
(7.3000,78.4573)
(7.5000,78.9118)
(7.7000,79.5664)
(7.9000,79.5364)
(8.1000,77.8267)
(8.3000,78.4409)
(8.5000,80.3542)
(8.7000,79.3877)
(8.9000,78.2206)
(9.1000,78.656)
(9.3000,80.1466)
(9.5000,79.4499)
(9.7000,78.2803)
(9.9000,79.6233)
(10.1000,79.1656)
(10.3000,78.2485)
(10.5000,79.9211)
(10.7000,78.7978)
(10.9000,78.8153)
(11.1000,79.2741)
(11.3000,78.3383)
(11.5000,80.4616)
(11.7000,78.6679)
(11.9000,78.1986)
(12.1000,79.2792)
(12.3000,79.1232)
(12.5000,79.2682)
(12.7000,79.0029)
(12.9000,78.6754)
(13.1000,78.7135)
(13.3000,79.7826)
(13.5000,80.7111)
(13.7000,79.0975)
(13.9000,78.0897)
(14.1000,80.4057)
(14.3000,78.0383)
(14.5000,78.1693)
(14.7000,78.3684)
(14.9000,78.5668)
(15.1000,78.7596)
(15.3000,79.211)
(15.5000,80.0254)
(15.7000,79.8389)
(15.9000,77.7771)
(16.1000,78.7064)
(16.3000,79.5658)
(16.5000,79.9544)
(16.7000,77.6218)
(16.9000,79.2597)
(17.1000,80.0641)
(17.3000,79.5085)
(17.5000,79.1472)
(17.7000,79.0926)
(17.9000,76.9469)
(18.1000,77.0297)
(18.3000,87.811)
(18.5000,77.6257)
(18.7000,78.5015)
(18.9000,78.8461)
(19.1000,81.0894)
(19.3000,78.2705)
(19.5000,79.3005)
(19.7000,79.8162)
(19.9000,78.8462)
(20.1000,79.5092)
(20.3000,78.6946)
(20.5000,78.8839)
(20.7000,79.579)
(20.9000,79.0294)
(21.1000,79.0248)
(21.3000,78.5764)
(21.5000,80.4646)
(21.7000,77.7802)
(21.9000,77.5159)
(22.1000,80.5411)
(22.3000,78.777)
(22.5000,79.2436)
(22.7000,78.3864)
(22.9000,78.6463)
(23.1000,79.2088)
(23.3000,80.6021)
(23.5000,78.3807)
(23.7000,79.025)
(23.9000,78.8984)
(24.1000,79.0834)
(24.3000,79.4805)
(24.5000,78.5529)
(24.7000,79.0245)
(24.9000,79.3908)
(25.1000,78.8486)
(25.3000,78.7722)
(25.5000,78.7164)
(25.7000,80.3884)
(25.9000,79.1137)
(26.1000,76.0782)
(26.3000,81.3389)
(26.5000,77.452)
(26.7000,114.747)
(26.9000,77.5789)
(27.1000,87.8153)
(27.3000,76.0848)
(27.5000,81.8696)
(27.7000,78.0335)
(27.9000,79.3573)
(28.1000,78.3658)
(28.3000,79.6897)
(28.5000,78.7255)
(28.7000,79.1217)
(28.9000,78.8063)
(29.1000,79.2793)
(29.3000,79.2288)
(29.5000,78.7384)
(29.7000,77.958)
(29.9000,78.98)
(30.1000,78.4736)
(30.3000,79.4101)
(30.5000,78.6694)
(30.7000,80.5356)
(30.9000,77.4799)
(31.1000,80.1432)
(31.3000,78.4073)
(31.5000,78.9972)
(31.7000,78.5994)
(31.9000,79.2108)
(32.1000,78.2599)
(32.3000,79.5334)
(32.5000,78.9958)
(32.7000,79.091)
(32.9000,79.364)
(33.1000,80.6463)
(33.3000,79.8365)
(33.5000,77.6109)
(33.7000,77.4204)
(33.9000,78.0076)
(34.1000,79.572)
(34.3000,79.5406)
(34.5000,78.5697)
(34.7000,78.721)
(34.9000,79.5517)
(35.1000,79.7387)
(35.3000,79.2636)
(35.5000,78.1928)
(35.7000,79.2597)
(35.9000,78.849)
(36.1000,79.2142)
(36.3000,78.2681)
(36.5000,78.8775)
(36.7000,80.4278)
(36.9000,79.9898)
(37.1000,78.2093)
(37.3000,78.518)
(37.5000,78.5772)
(37.7000,78.8145)
(37.9000,79.4671)
(38.1000,78.691)
(38.3000,77.8937)
(38.5000,79.9916)
(38.7000,78.0382)
(38.9000,80.0178)
(39.1000,78.6786)
(39.3000,79.19)
(39.5000,78.1794)
(39.7000,79.684)
(39.9000,79.1834)
(40.1000,77.772)
(40.3000,79.1623)
(40.5000,79.5589)
(40.7000,78.5083)
(40.9000,79.6391)
(41.1000,78.5413)
(41.3000,79.6549)
(41.5000,79.0322)
(41.7000,78.7174)
(41.9000,78.7675)
(42.1000,79.9311)
(42.3000,78.5904)
(42.5000,79.5029)
(42.7000,78.3842)
(42.9000,80.6932)
(43.1000,78.2253)
(43.3000,78.3516)
(43.5000,82.4995)
(43.7000,79.7906)
(43.9000,75.9267)
(44.1000,88.4655)
(44.3000,74.2933)
(44.5000,79.8754)
(44.7000,107.47)
(44.9000,0)
};
\end{axis}
\begin{axis}[scale only axis, width=0.9\linewidth, height=0.4\linewidth, ybar, bar shift=0pt, bar width=0.2, xmin=0, xmax=45, ymin=0, ymax=118, hide axis]
\addplot[fill=green!55!black, draw=green!30!black, thin] coordinates {
(0.1000,29.6066)
(0.3000,0)
(0.5000,8.06894)
(0.7000,75.66)
(0.9000,38.3965)
(1.1000,11.4706)
(1.3000,44.4465)
(1.5000,61.9389)
(1.7000,52.4989)
(1.9000,41.8913)
(2.1000,42.9027)
(2.3000,59.9082)
(2.5000,52.8544)
(2.7000,58.522)
(2.9000,50.051)
(3.1000,56.5959)
(3.3000,55.5632)
(3.5000,59.4354)
(3.7000,53.6546)
(3.9000,60.8854)
(4.1000,54.7736)
(4.3000,67.7966)
(4.5000,47.7474)
(4.7000,58.6994)
(4.9000,56.7501)
(5.1000,68.6162)
(5.3000,46.6367)
(5.5000,60.9878)
(5.7000,62.1078)
(5.9000,50.3885)
(6.1000,61.5559)
(6.3000,50.913)
(6.5000,60.9499)
(6.7000,57.3884)
(6.9000,54.9276)
(7.1000,54.9788)
(7.3000,65.0646)
(7.5000,58.3145)
(7.7000,47.2317)
(7.9000,52.5231)
(8.1000,57.1596)
(8.3000,60.7436)
(8.5000,55.7254)
(8.7000,54.4934)
(8.9000,53.6276)
(9.1000,56.9793)
(9.3000,55.7293)
(9.5000,58.3186)
(9.7000,48.2181)
(9.9000,52.1364)
(10.1000,59.3122)
(10.3000,55.7346)
(10.5000,52.881)
(10.7000,56.5195)
(10.9000,52.364)
(11.1000,68.3286)
(11.3000,41.7722)
(11.5000,59.4236)
(11.7000,54.3491)
(11.9000,61.2312)
(12.1000,49.2486)
(12.3000,60.8499)
(12.5000,47.1515)
(12.7000,56.5459)
(12.9000,63.0633)
(13.1000,46.2124)
(13.3000,57.4017)
(13.5000,47.8118)
(13.7000,50.7924)
(13.9000,47.3441)
(14.1000,61.6974)
(14.3000,56.6778)
(14.5000,60.4284)
(14.7000,53.5734)
(14.9000,55.1875)
(15.1000,57.9894)
(15.3000,53.9096)
(15.5000,51.8392)
(15.7000,55.0546)
(15.9000,58.0571)
(16.1000,59.0403)
(16.3000,57.9457)
(16.5000,51.5132)
(16.7000,60.0031)
(16.9000,55.6229)
(17.1000,54.3531)
(17.3000,51.823)
(17.5000,52.7889)
(17.7000,53.7184)
(17.9000,53.9105)
(18.1000,52.4982)
(18.3000,25.5035)
(18.5000,60.1977)
(18.7000,63.0586)
(18.9000,56.9628)
(19.1000,45.6193)
(19.3000,56.4431)
(19.5000,48.9683)
(19.7000,60.3834)
(19.9000,55.8595)
(20.1000,52.954)
(20.3000,55.8686)
(20.5000,51.9517)
(20.7000,56.3545)
(20.9000,57.4579)
(21.1000,50.6859)
(21.3000,61.3988)
(21.5000,49.73)
(21.7000,38.6827)
(21.9000,76.8806)
(22.1000,53.0026)
(22.3000,55.9653)
(22.5000,53.8036)
(22.7000,61.5257)
(22.9000,53.2617)
(23.1000,68.6003)
(23.3000,40.0959)
(23.5000,60.1732)
(23.7000,52.1177)
(23.9000,58.4683)
(24.1000,57.7229)
(24.3000,58.88)
(24.5000,55.4192)
(24.7000,53.0547)
(24.9000,53.0244)
(25.1000,51.2511)
(25.3000,57.1)
(25.5000,56.8901)
(25.7000,45.715)
(25.9000,54.3152)
(26.1000,26.8513)
(26.3000,58.0949)
(26.5000,75.5013)
(26.7000,0.40736)
(26.9000,73.8522)
(27.1000,28.9049)
(27.3000,52.0199)
(27.5000,48.5047)
(27.7000,55.346)
(27.9000,54.6553)
(28.1000,50.2936)
(28.3000,57.6848)
(28.5000,50.2511)
(28.7000,58.7733)
(28.9000,62.4516)
(29.1000,53.4082)
(29.3000,54.4928)
(29.5000,55.2125)
(29.7000,52.3976)
(29.9000,56.355)
(30.1000,55.3528)
(30.3000,61.0024)
(30.5000,53.2041)
(30.7000,49.4531)
(30.9000,70.7691)
(31.1000,42.4186)
(31.3000,64.4841)
(31.5000,47.1204)
(31.7000,60.6597)
(31.9000,45.4075)
(32.1000,66.7146)
(32.3000,50.1188)
(32.5000,55.2768)
(32.7000,51.0542)
(32.9000,52.8086)
(33.1000,46.4145)
(33.3000,51.1901)
(33.5000,54.3679)
(33.7000,54.2907)
(33.9000,55.6347)
(34.1000,55.2761)
(34.3000,49.3742)
(34.5000,54.1265)
(34.7000,54.1259)
(34.9000,56.0702)
(35.1000,48.2575)
(35.3000,55.9419)
(35.5000,48.2742)
(35.7000,56.3204)
(35.9000,58.1641)
(36.1000,52.4539)
(36.3000,56.9495)
(36.5000,57.3261)
(36.7000,49.21)
(36.9000,44.9825)
(37.1000,63.3537)
(37.3000,49.8484)
(37.5000,53.9984)
(37.7000,52.3834)
(37.9000,57.5648)
(38.1000,54.8137)
(38.3000,50.6194)
(38.5000,53.2405)
(38.7000,60.8267)
(38.9000,47.6411)
(39.1000,56.3183)
(39.3000,49.8248)
(39.5000,58.4087)
(39.7000,44.3374)
(39.9000,64.5554)
(40.1000,58.0158)
(40.3000,54.8832)
(40.5000,43.0375)
(40.7000,66.6063)
(40.9000,52.7106)
(41.1000,60.6949)
(41.3000,50.6897)
(41.5000,54.302)
(41.7000,54.5937)
(41.9000,56.7647)
(42.1000,49.0577)
(42.3000,54.6551)
(42.5000,50.9598)
(42.7000,54.9482)
(42.9000,47.9079)
(43.1000,48.1368)
(43.3000,51.9464)
(43.5000,39.7242)
(43.7000,37.6748)
(43.9000,51.9266)
(44.1000,29.4139)
(44.3000,10.3417)
(44.5000,31.673)
(44.7000,0.97416)
(44.9000,0)
};
\end{axis}
\end{tikzpicture}

  \caption{A histogram of the frequency (green) and the average length
    (blue) of the 13.5 million random reference vectors of 500-ominoes,
    distributed into buckets of width $0.2^\circ$.
    The buckets
    $0.2^\circ$--$0.4^\circ$ and
    $44.8^\circ$--$45^\circ$ are empty.
    The bucket 
    $26.6^\circ$--$26.8^\circ$ contains 
    only 465 vectors, and
    these are atypically long.
  }
  \label{fig:angle-225}
\end{figure}

 Figure~\ref{fig:angle-225}
shows a histogram of 
the frequency and the average vector length
of the generated reference vectors,
after distributing them into buckets of
width $0.2^\circ$ according to their angle $\alpha=\arctan \frac yx$.
The average number of vectors per bucket
is $60\,000$.
In case of a tie between $k$ reference vectors of equal length, we have counted each of
these vectors with weight $\frac{1}{k}$.

Two observations stand out: Vectors around the main axes ($0^{\circ}$) are rare, and exact diagonal lines (slope $1$,
i.e.,\ $45^\circ$) are completely missing. Additionally, we observe
fewer but longer reference vectors slightly above slope $\frac{1}{2}$
($26.57^\circ$), and a weaker version of this effect near slope
$\frac{1}{3}$ ($18.43^\circ$).

The missing $45^\circ$ directions are explained by the fact that
we select longest distances between \emph{edge midpoints},
see \autoref{fig:longest-edge}.
A reference vector in diagonal direction can 
always be improved 
by ``nudging'' it off the diagonal.
In fact, the closest vector to a diagonal direction in our sample is
$\binom{95}{94}$, with angle $44.70^\circ$.

The irregularities near the ``simple'' slopes $0,\frac12,\frac13,$
etc.\ are due to the properties of the rational numbers and have nothing
to do with polyominoes per se.
The fractions $\frac{p}{q}$ with $p, q$ bounded by some constant are
not distributed evenly, but they leave gaps
near simple fractions.
For example,
the bucket $0^\circ$--$0.2^\circ$ contains all the axis-parallel
vectors
$\binom10,\binom20,\binom30,\ldots$,
but the shortest vector that is not of this form is the vector
$\binom{287}{1}$.
The shortest vector in the bucket
    $0.2^\circ$--$0.4^\circ$ is 
$\binom{144}{1}$. The average length of the reference vectors is around
$80$. There exist 500-ominoes with diameters as large as that, but they are
extremely rare, and they do not appear in our sample.
The vector with the smallest positive angle that does appear in our sample is
$\binom{140}1$ with angle $0.41^\circ$. Accordingly,
the vectors
 in the bucket
    $0.4^\circ$--$0.6^\circ$ are on average longer than average.

Let us look near the slope $\frac12$:
The shortest vector in the range
$26.4^\circ$--$26.6^\circ$ that is not a multiple of $\binom21$ is
$\binom{141}{70}$.
The shortest vectors in the adjacent range
$26.6^\circ$--$26.8^\circ$
are $\binom{99}{50},\binom{101}{51},\ldots$,
and there are few samples of this sort.
The irregularity near slope $\frac{1}{3}$ is less pronounced but still
visible.
Gaps near less simple fractions like $\frac{2}{3}$, $\frac{1}{4}$,
etc., are smaller; they would show up in a histogram with narrower buckets.

\begin{figure}[htb]
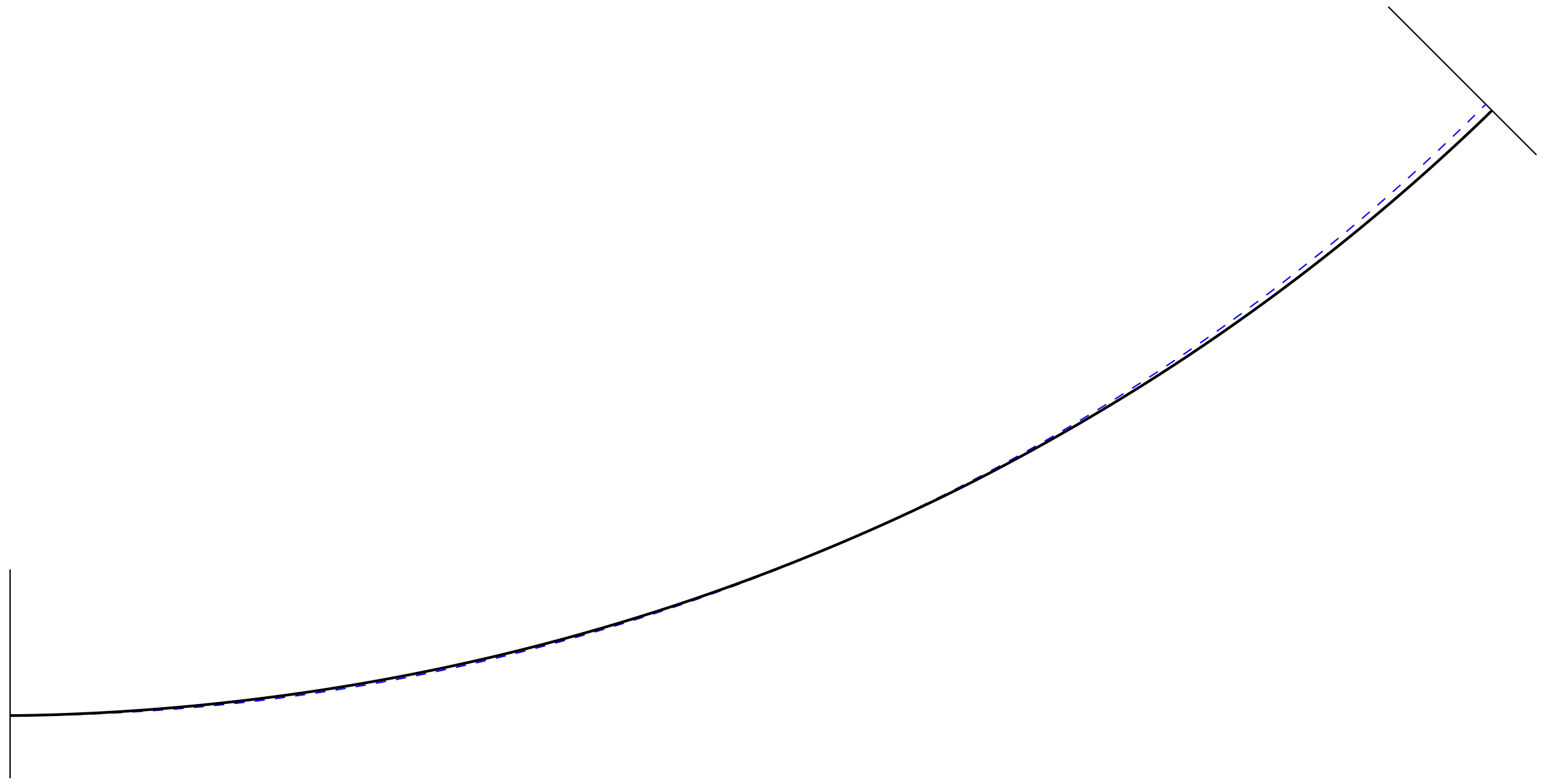

  \centering
\input chain 
%
\caption{The reference vectors of
13.5 million
  random $500$-ominoes, sorted by
slope and concatenated into a convex polygonal chain from $(0,0)$ to $(X,Y)$
(solid). The dashed blue curve is the arc of the circle of radius $R=X+Y$ around
$K=(0,X+Y)$, which starts at $(0,0)$ with a horizontal tangent. The
two little straight
segments lie on the rays from $K$ through the start and the end point of the
chain; they enclose an angle of $45^\circ$.
}
  \label{fig:chain}
\end{figure}

\begin{figure}[htb]
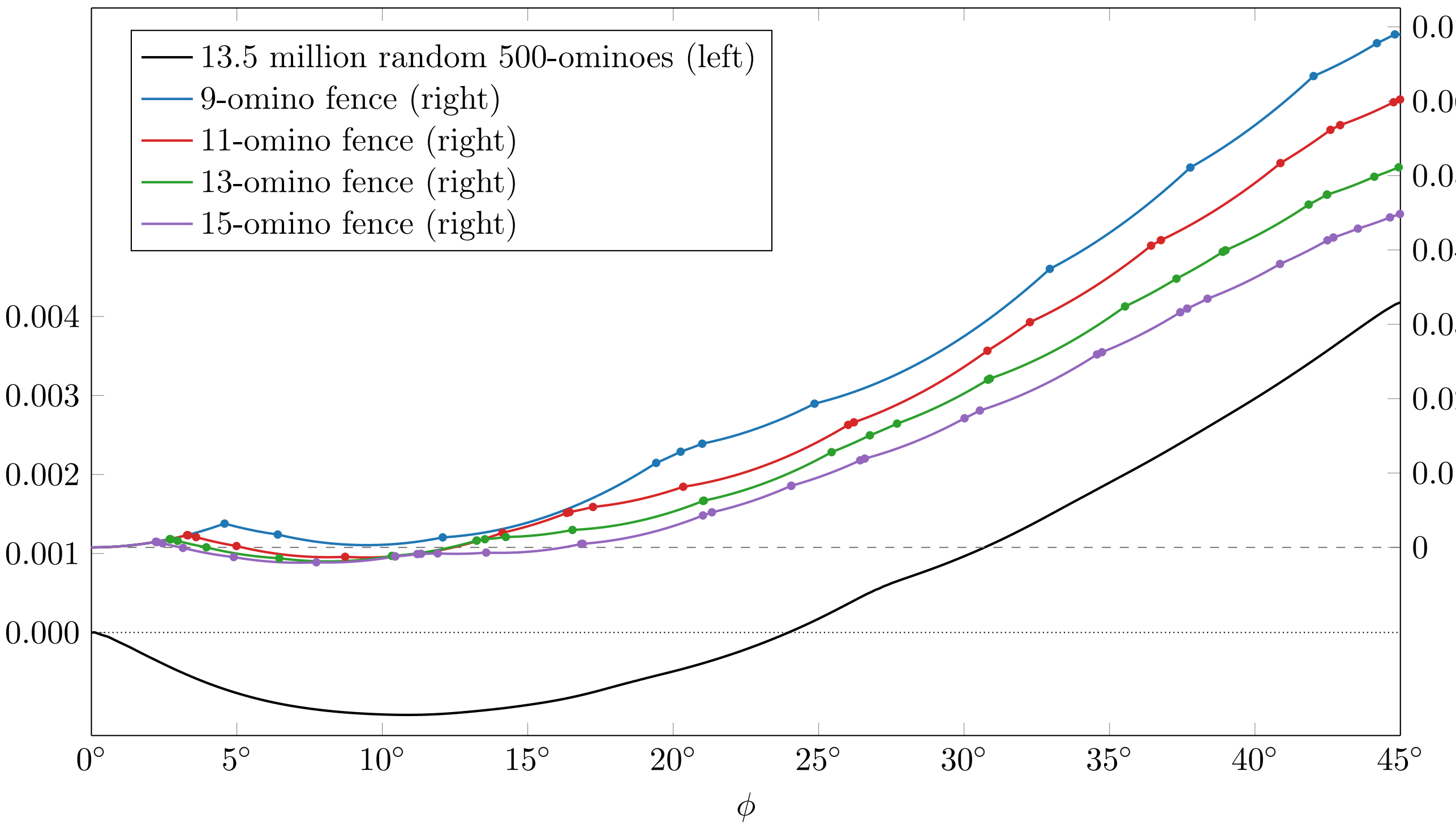

  \centering
\input deviation




\caption{
  The deviation 
  of the chain from the circle
  in Figure~\ref{fig:chain},
  relative to the radius $R$ of the circle.
  The parameter $\phi$ is the angle as seen from the center of the
  circle, starting from the downward direction.
  For comparison, the colored curves show the result from
  the same S-to-SE 1/8-sector of the polygons
of \autoref{fig:large-ref-polygons}, drawn with a larger vertical
scale that is shown on the right edge of the diagram.
The dots mark the true vertices, where a change of direction occurs.
  The straight polygon edges appear curved because the plot compares
  them with a circle.}
  \label{fig:deviation}
\end{figure}

Let us investigate the shape of the $500$-omino fence. For this, we
sort our 13.5 million samples by slope, concatenate them, and compare
them with an arc of $45^\circ$ of a circle, as shown in
Figure~\ref{fig:chain}. Since differences are hard to
see at this scale, we also plot the relative difference of the arc radius at each
angle in Figure~\ref{fig:deviation},
also for the smaller polyominoes of Figure~\ref{fig:large-ref-polygons}.
From the small polyominoes, one would not have guessed that
 the curve runs initially \emph{inside} the circle for quite a large range of~$\phi$. One can also recognize some irregularity near
$\arctan \frac12 \approx 26.57^\circ$.


We see in Figure~\ref{fig:deviation} that
  the relative deviation from the circle decreases slowly with
  increasing~$n$; thus, it seems that the shape converges to a circle
  after all.
  
In comparison to the huge number of 15-ominoes, the number of true corners,
which are marked by dots, is
remarkably small.
The reason is that,
for $n$-ominoes, the vectors $\binom xy$ that appear are restricted by
the condition $|x|+|y|\le n$, and practically all of these vectors
appear, because one can easily construct an $n$-omino with a
given
longest reference vector. The only exceptions are vectors $\binom{x}{y}$ that are
shorter than about $2\sqrt{n/\pi}$, i.e.,\ with $x^2+y^2\lesssim 4n/\pi$, because of the isodiametric
inequality.
(The reader may try to guess what the less frequent directions are,
those that appear between close consecutive dots. (Some dots are so
close that they are not distinguishable in the plot.)

Finally, we estimate the
area of the $500$-omino fence. To this end, we combine $8$ copies
of the chain formed by the 13.5 million samples in
Figure~\ref{fig:chain} and calculate the area of the resulting closed
loop. Then we scale it up as if the loop consisted of all
$2.0 \cdot 10^{300}$ free $500$-ominoes instead of just $108$
million vectors. Doing this for our data, we arrive at a reference
polygon area of $2.0 \cdot 10^{603}$.
The estimate for the actual $500$-omino fence is the same, because for
large $n$, the sum of the corrections is negligible. More concretely, since the correction for any $n$-omino can be trivially bounded by $n^2$, the sum of all corrections is at most $n^2 A'(n) \approx 5.0 \cdot 10^{305}$, i.e.,\ around $300$ orders of magnitude smaller than our estimate.
However, due to the enormous uncertainty in the true number of free $500$-ominoes mentioned at the beginning of this subsection, this estimate should only be taken as a rough order of magnitude.

Appendix~\ref{sec:refvec-length-estimate} studies how the expected length of the longest reference vector of a random $n$-omino grows with~$n$.

\paragraph{Interactive lab.} In the \emph{limit-shape lab} (\url{https://erikaroldanroa.github.io/fence-challenge/limit-shape-lab/}) one can see the reference polygons of the optimal fences for $n = 4, 5, 6$ and of the greedy solutions for $n = 7, \dots, 15$ (Figure~\ref{fig:large-ref-polygons}) next to a circle, and experiment with reference vector samples for $n = 500$.

\section{Three-dimensional walls}\label{sec:3d}

In three dimensions, polyominoes become polycubes, fences become \emph{walls}, and the enclosed area becomes an enclosed \emph{volume}.

\subsection{Definitions}\label{sec:3d-def}

The definitions of Section~\ref{sec:def-area} carry over, with the pair $(8,4)$ of connectivities replaced by $(26,6)$ \cite{KongRosenfeld89}. Two cells of $\mathbb{Z}^3$ are \emph{6-adjacent} if they share a face, and \emph{26-adjacent} if they share at least a vertex. As in the plane, free cells are regarded as connected even if they touch only at a vertex or an edge, so the pieces of a wall must be joined face to face.

\begin{definition}\label{def:wall}
  A placement of polycubes is a \emph{wall} if the pieces are pairwise disjoint and their union~$\Phi$ is 6-connected. The cells of~$\Phi$ are the \emph{wall cells}, and all other cells are \emph{free}. The \emph{outside} is the unique infinite 26-connected component of free cells; the other free cells are \emph{enclosed}. The \emph{volume} of the wall is the number of cells of the largest 26-connected component of enclosed cells.
\end{definition}

The \emph{$n$-polycube wall problem} asks for the largest volume of a wall that uses each $n$-polycube exactly once.

\paragraph{Orientation.}
For polycubes we use a different convention than for polyominoes: mirror images count as \emph{distinct} pieces. Thus there are $29$ pentacubes ($23$ if mirror images were identified), $166$ hexacubes, and $1023$ heptacubes. We did not choose this convention deliberately when we started to work on walls, but it can be defended: a polyomino can be flipped over in the three-dimensional world we live in, whereas a polycube cannot be mirrored. In any case, the wall problem is a legitimate question for both interpretations.

\subsection{An isoperimetric upper bound}\label{sec:3d-bound}

Extending the dynamic programming approach, which was very successful in 2D, is not feasible: there are more $n$-polycubes than $n$-ominoes, each polycube has more orientations, and there are drastically more combinations of pairwise adjacencies (instead of the single cyclic sequence of Lemma~\ref{lem:interval}). Instead, we relax the problem to an isoperimetric one: we drop the requirement that the wall is composed of $n$-polycubes and only count its unit cubes.

Since a free cell may escape through 26-connected paths, a wall that encloses a set $S$ of cells must contain every cell that shares at least a vertex with~$S$, that is, the set $N(S)\setminus S$. Here $N(S)$ is the closed neighborhood of $S$ in the graph on $\mathbb{Z}^3$ in which two cells are adjacent if their $\ell_\infty$-distance is~$1$ (i.e.,\ they share a vertex). Hence a wall built from $n$-polycubes with a total of $B$ cubes encloses a volume of at most
\[
  V^*(B) \;=\; \max\{\,|S| : |N(S)\setminus S| \le B\,\}.
\]
This is the \emph{vertex-isoperimetric problem}: for a given number of cells, find the set that can be enclosed by the fewest cubes.

In the plane, the extremal sets are squares and almost-squares; in three dimensions they are cubes and \emph{almost-cubes}~\cite{alonso-cerf}. Formally, an almost-cube $Q$ of dimension $d$ and order $k \in \mathbb{N}$ ($k \ge 2$) is defined inductively: it consists of a primary $d$-dimensional cuboid
\[
  C = \{1, \dots, k_1\} \times \{1, \dots, k_2\} \times \dots \times \{1, \dots, k_d\}
\]
with edge lengths $k_i \in \{k, k+1\}$, and a lower-dimensional almost-cube $Q'$ of dimension $d - 1$ and order $k' \le k$ that is attached to exactly one boundary face of $C$. Veomett and Radcliffe~\cite{VeomettRadcliffe12} solved the vertex-isoperimetric problem exactly: for every number of cells, an initial segment of an explicit ordering of $\mathbb{Z}^3$ minimizes $|N(S)|$ among all sets with this number of cells, and these sets are almost-cubes. Consequently, $V^*(B)$ is attained by an almost-cube~$S$, which we call the \emph{target region}; the set $W=N(S)\setminus S$ is the \emph{target wall}.

The maximum is attained by an almost-cube, but not only by almost-cubes. This already happens in the plane: the $2\times5$ rectangle needs a fence of $18$ cells, the same as the almost-square with $10$ cells (a $3\times3$ square with one extra cell). Taking the budget viewpoint, however, the same $18$ cells enclose a $3\times4$ rectangle with $12$ cells; in the plane, the optimal solutions are squares or rectangles whose side lengths differ by~$1$. In three dimensions, a budget of $146$ cubes (one more than the $29$ pentacubes supply) allows the maximum volume $54$, which is attained by a $3\times3\times6$ box (its wall has $5\cdot5\cdot8-54=146$ cubes) as well as by the almost-cube formed by a $3\times4\times4$ box with a $2\times3$ layer attached to one of its sides.

\paragraph{Surface area.}
For the related problem of minimizing the \emph{surface area} of $n$ unit cubes in $\mathbb{Z}^d$, the optimum is also known. It follows from a theorem of Bollob\'as and Leader~\cite{bl-isoperimetric-91} on the maximum number of edges induced by $n$ vertices of the grid: initial segments of their cube order (which is the order later used by Veomett and Radcliffe, up to relabeling the coordinates) are optimal, and they are special almost-cubes. As above, not every optimal set is an almost-cube. Alonso and Cerf~\cite[Corollaries 2.5 and 3.4]{alonso-cerf}, who were not aware of this work, gave an explicit formula in dimensions two and three. They also described the minimizing shapes, completely in two dimensions and partially in three (their ``canonical shapes'' are exactly the initial segments of Bollob\'as and Leader).\footnote{Theorem 3.5 in~\cite{alonso-cerf}, which gives conditions under which only almost-cubes are optimizers, is subject to exceptions: one family are the $j\times j\times (j+3)$ boxes, another family are the $j\times (j+1)\times (j+1)$ boxes from which $2j$ cubes along two opposite parallel edges of length~$j$ have been removed. Accordingly, as confirmed by Rapha\"el Cerf (personal communication, October 2026), the cases $n=j^2(j+3)=j(j+1)^2+(j-1)j$ and $n=j(j+1)^2-2j=j^2(j+1)+(j-1)j$ must be excluded from the list in Theorem~3.5. The gap in the proof is that, for these values, the projections $p_1,p_2,p_3$ can map a non-principal polyomino to a principal one. This was first noticed by the AI agent Claude from Anthropic and then discussed with Rapha\"el Cerf.}
Empirically, surface area and wall size $|N(S)\setminus S|$ seem to be minimized by the same configurations for a given number $|S|$ of cells, but this has not been proved.

\paragraph{Exact cover problem.}

If the target wall~$W$ of an optimal target region can be covered by our $n$-polycubes, the wall has volume $V^*(B)$, and we have found an optimal solution of the 3D wall problem. Often this is not possible. In that case we try smaller almost-cubes and also allow some polycubes to protrude from the target wall to the outside (but never into~$S$), so that the method provides only an upper and a lower bound.

Thus we have reduced the placement of polycubes to an exact cover problem on the target wall. For small instances such as the pentacubes, exact solvers, including Knuth's Dancing Links algorithm~\cite{kn4b,knuth2000dancing} and integer linear programming, are computationally feasible. For larger polycubes, where the search space explodes, we use a greedy placement heuristic that positions new pieces adjacent to the boundaries of previously placed pieces.

\subsection{Pentacubes}\label{sec:pentacubes}

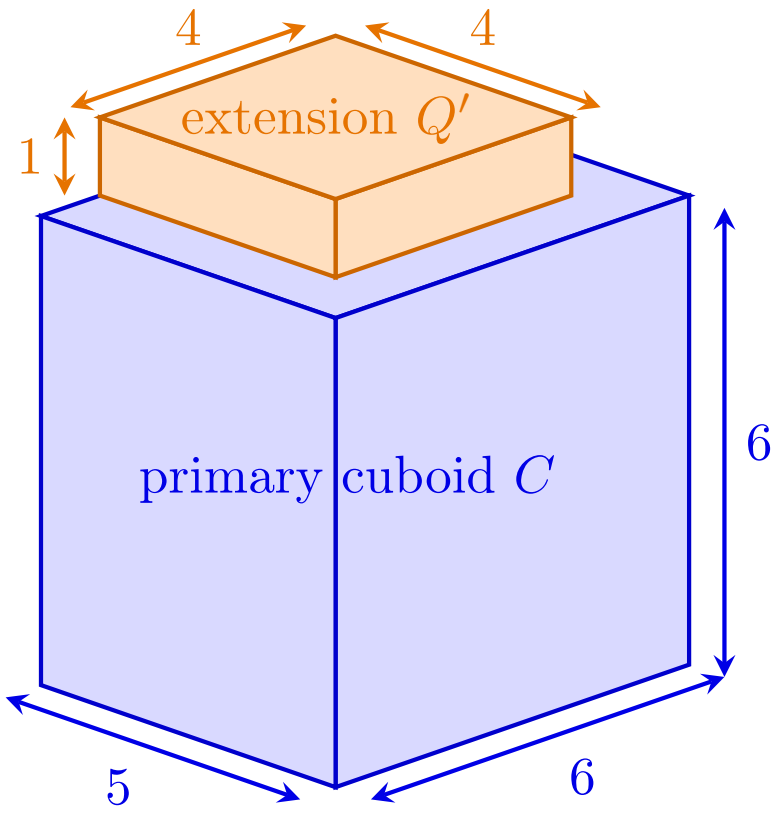
\begin{figure}[htbp]
    \centering
    \begin{tikzpicture}[
        x={(0.866cm, -0.3cm)}, y={(0.866cm, 0.3cm)}, z={(0cm, 1.149cm)},
        scale=0.45,
        >=stealth
    ]
        \tikzset{
            facet main/.style={fill=blue!15, draw=blue!80!black, thick},
            facet sub/.style={fill=orange!25, draw=orange!80!black, thick}
        }

        \filldraw[facet main] (0,0,6) -- (5,0,6) -- (5,6,6) -- (0,6,6) -- cycle;
        \filldraw[facet main] (5,0,0) -- (5,6,0) -- (5,6,6) -- (5,0,6) -- cycle;
        \filldraw[facet main] (0,0,0) -- (5,0,0) -- (5,0,6) -- (0,0,6) -- cycle;

        \draw[<->, thick, blue!90!black] (0,-0.6,0) -- (5,-0.6,0) node[midway, below left, font=\small] {$5$};
        \draw[<->, thick, blue!90!black] (5.6,0,0) -- (5.6,6,0) node[midway, below right, font=\small] {$6$};
        \draw[<->, thick, blue!90!black] (5.6,6,0) -- (5.6,6,6) node[midway, right, font=\small] {$6$};

        \filldraw[facet sub] (4,1,6) -- (4,5,6) -- (4,5,7) -- (4,1,7) -- cycle;
        \filldraw[facet sub] (0,1,6) -- (4,1,6) -- (4,1,7) -- (0,1,7) -- cycle;
        \filldraw[facet sub] (0.,1,7) -- (4,1,7) -- (4,5,7) -- (0,5,7) -- cycle;

        \draw[<->, thick, orange!90!black] (0,5.5,7) -- (4,5.5,7) node[midway, above, font=\small] {$4$};
        \draw[<->, thick, orange!90!black] (-0.5,1,7) -- (-0.5,5,7) node[midway, above, font=\small] {$4$};
        \draw[<->, thick, orange!90!black] (-0.3,0.7,6) -- (-0.3,0.7,7) node[midway, left, font=\small] {$1$};

        \node[blue!90!black, font=\small] at (2.5,3,2.5) {primary cuboid
          $C$\ \ \ };
        \node[orange!90!black, font=\small, above] at (2,3,6.4) {extension $Q'\,\,$};

    \end{tikzpicture}
    \caption{The closed neighborhood $N(S)$ of the optimal target
      region $S$ for the pentacubes: an almost-cube with a $5 \times 6
      \times 6$ primary cuboid $C$ and a $4 \times 4 \times 1$
      extension $Q'$. The target region $S$ (a $3\times4\times4$ box
      with a $2\times2$ layer on top of it) lies inside; the target wall is $W=N(S)\setminus S$.}
    \label{fig:almost-cube-hull}
\end{figure}

The $29$ pentacubes consist of $145$ unit cubes. The minimum wall size is always even; it is $144$ for $|S|=52$ and $146$ for $|S|=53$. Hence $145$ cubes can enclose at most $52$ cells. The optimal target region~$S$ is the almost-cube formed by a $3\times4\times4$ box with a $2\times2$ layer attached to one of its sides, with $|S|=3\cdot4\cdot4+2\cdot2=52$. Its closed neighborhood $N(S)$ is the almost-cube formed by a $5\times6\times6$ box and a $4\times4\times1$ layer (\autoref{fig:almost-cube-hull}), and the target wall $W=N(S)\setminus S$ has $144$ cells.

We encode the covering of $W$ as an exact cover problem. Let $\mathcal{P}$ denote the set of pieces. A placement option $k = (p, \tau, t)$ puts piece $p \in \mathcal{P}$ in orientation~$\tau$ at translation~$t$ such that the set of cells $\text{cells}(k)$ it occupies lies in~$W$. Since the cells of~$W$ must be covered exactly, $|W|$ must be the sum of the piece sizes. Here, however, $|W|=144$ while the pentacubes have $145$ cubes, so one cube has to lie outside $N(S)$. We therefore select one pentacube and the cube of it that protrudes, and use in the model the tetracube formed by its remaining four cubes. Then $|W| = 28\cdot 5 + 4 = 144$.

We construct a binary matrix $M$ whose rows are the placement options. Its columns are partitioned into two groups:
\begin{enumerate}
    \item Cell columns ($c \in W$): primary constraints requiring that each cell of the target wall is covered by exactly one piece.
    \item Piece columns ($p \in \mathcal{P}$): constraints ensuring that each piece is used at most once (and hence exactly once, since the piece sizes add up to $|W|$).
\end{enumerate}
With binary decision variables $x_k \in \{0, 1\}$ indicating whether placement option $k$ is selected, the exact cover formulation is
\begin{align}
    \sum_{k : c \in \text{cells}(k)} x_k &= 1 \quad \forall c \in W, \label{eq:cell-cover}\\
    \sum_{k : \text{piece}(k) = p} x_k &\le 1 \quad \forall p \in \mathcal{P}. \label{eq:piece-usage}
\end{align}
The structure of $M$ is illustrated in Table~\ref{tab:exact-cover-matrix}.

\begin{table}[htbp]
    \centering
    \small
    \begin{tabular}{|l|cccc|cccc|}
        \hline
        & \multicolumn{4}{c|}{Target Wall Cells ($c \in W$)} & \multicolumn{4}{c|}{Pieces ($p \in \mathcal{P}$)} \\
        Placement Option ($k$) & $c_1$ & $c_2$ & $\dots$ & $c_{|W|}$ & $p_1$ & $p_2$ & $\dots$ & $p_{|\mathcal{P}|}$ \\
        \hline
        Option 1: $(p_1, \tau_1, t_1)$ & 1 & 0 & $\dots$ & 0 & 1 & 0 & $\dots$ & 0 \\
        Option 2: $(p_1, \tau_2, t_2)$ & 0 & 1 & $\dots$ & 1 & 1 & 0 & $\dots$ & 0 \\
        Option 3: $(p_2, \tau_1, t_3)$ & 1 & 1 & $\dots$ & 0 & 0 & 1 & $\dots$ & 0 \\
        $\vdots$ & $\vdots$ & $\vdots$ & $\ddots$ & $\vdots$ & $\vdots$ & $\vdots$ & $\ddots$ & $\vdots$ \\
        \hline
    \end{tabular}
    \caption{Schematic representation of the exact cover matrix $M$. Rows represent valid placements of pieces; columns enforce exact coverage of the target wall cells $W$ and the piece usage constraints for~$\mathcal{P}$.}
    \label{tab:exact-cover-matrix}
\end{table}

\begin{figure}[htb]
  \centering
  \tabskip=0pt
  \halign to \textwidth{\hfil#\hfil\tabskip=0pt plus 1fill
    &\hfil#\hfil&\hfil#\hfil\tabskip=0pt\cr
  \includegraphics[scale=0.2]{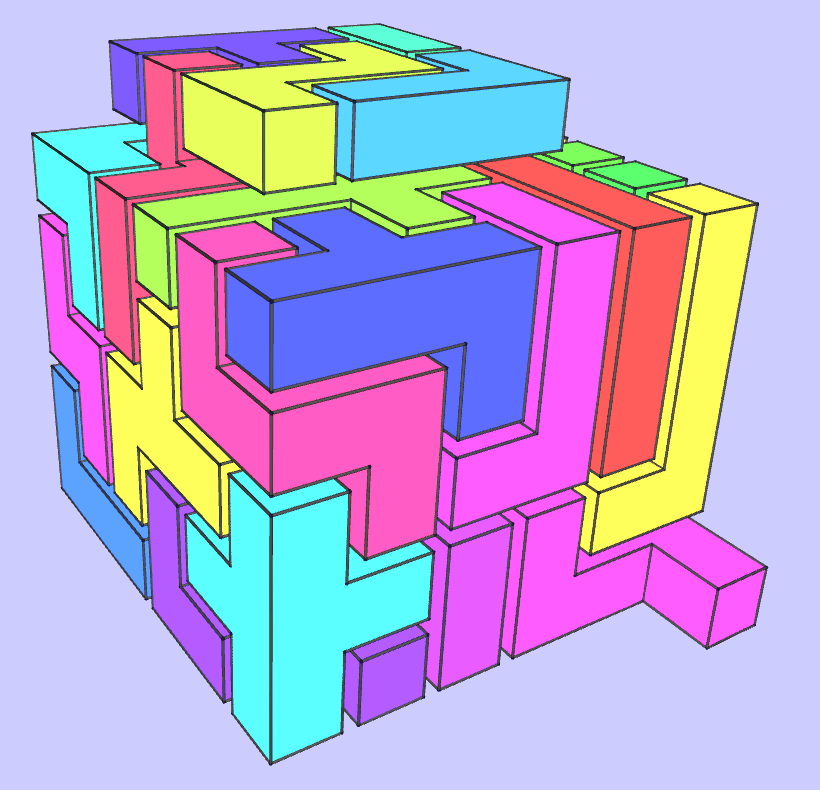}&
  \includegraphics[scale=0.2]{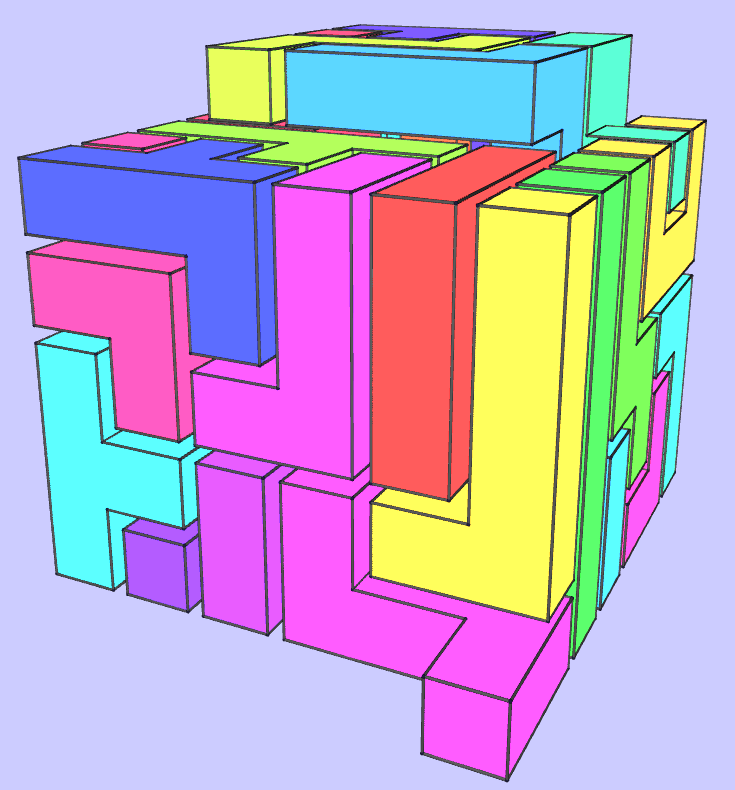}&
  \includegraphics[scale=0.2]{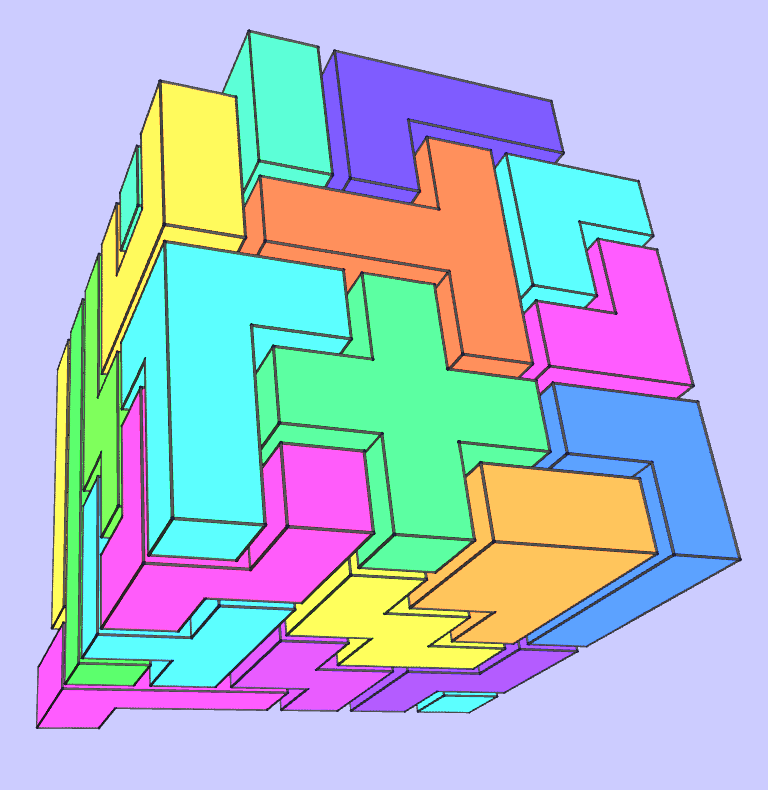}\cr
  (a) view from front left
  &(b) view from front right
  &(c) back and bottom view
  \cr
  }
  \caption{One optimal wall built from 29 pentacubes. It encloses a
    volume of 52.
In contrast to \autoref{fig:almost-cube-hull}, this wall is based on a
$6\times 6\times 5$ primary cuboid (of height 5).
    Each piece is shrunk to make it easier to distinguish.
  These pictures are taken from the viewer
\url{https://erikaroldanroa.github.io/fence-challenge/archive/polycube-viewer/Polycubefences.html}.
  The cube lab \url{https://erikaroldanroa.github.io/fence-challenge/cube-lab/?wall=52}
  offers 
  more interaction: 
  the wall
  can be 
  cut open and spread apart, and the pieces can be manipulated.
  }

  \label{fig:52}
\end{figure}

Using Knuth's \texttt{ssxcc-binary} solver~\cite{knuth-ssxcc}, which is based on his Dancing Links algorithm \cite{kn4b,knuth2000dancing}, we find a valid cover of the target wall (\autoref{fig:52}). Together with the upper bound of Section~\ref{sec:3d-bound}, this proves that the maximum volume for the pentacubes is $52$.

We additionally tried to enumerate all possible constructions, but stopped the backtracking search after several weeks, at which point we had already found more than 35 million solutions.

\subsection{Hexacubes and heptacubes}

For larger polycubes, exact solvers are computationally infeasible. Instead, we find lower bounds by applying a greedy placement heuristic to the target wall of a prescribed box-shaped target region. With this, we achieved the following results:
\begin{itemize}
    \item \textbf{Hexacubes:} We can greedily cover the target wall of the $11 \times 11 \times 11$ cube, thus achieving an enclosed volume of $V = 1\,331$. The optimal isoperimetric solution yields an upper bound of $1\,633$.
    
    \item \textbf{Heptacubes:} Using the $28 \times 33 \times 28$ box as target region (whose closed neighborhood is the $30 \times 35 \times 30$ box), the greedy heuristic achieved an enclosed volume of $V = 25\,872$ using 164 of the 166 heptacubes (the two remaining pieces can be stuck on the outside). The optimal isoperimetric solution yields an upper bound of $37\,446$.
\end{itemize}

Our implementation of the greedy placement algorithm currently only considers box-shaped target regions. In each step, it considers all placements of the remaining polycubes that cover an uncovered target-wall cell adjacent to the already placed pieces, and greedily selects the feasible placement with the highest score. Placements that overlap previously placed cubes or put any cube into the target region are infeasible. The score of a feasible placement option~$k$ of polycube~$p$ is $s(k) = w_H\, H(k) + A(k) - \lambda\, c(p)$, where $H(k)$ is the number of unit cubes of the placed polycube that cover still uncovered cells of the target wall, and $A(k)$ is the number of face contacts between the placed polycube and previously placed cubes. The weight was chosen as $w_H = 2$ for hexacubes and $w_H = 4$ for heptacubes. The term $c(p)$ is the largest number of unit cubes of~$p$ in a single axis-parallel plane, with $\lambda = 0$ for hexacubes and $\lambda = 4$ for heptacubes; it is independent of the placement and thus only influences which polycube is chosen.

In the current solutions, many cubes lie outside the target wall, creating bumps rather than contributing to the enclosed volume. Local rearrangements could flatten these bumps and thereby improve the volume; we leave this computational refinement for future work.

\paragraph{Interactive lab.}
Readers can try building a polycube wall and explore our
solutions
at the \emph{cube lab}:
\url{https://erikaroldanroa.github.io/fence-challenge/cube-lab/}.
Additional visualizations of the solutions are available at \url{https://zopfbit.github.io/polycube-solutions/}.

\section{Conclusion}

We hope that this paper conveys how much mathematical richness is
hidden in this recreational geometric problem from the 1960s. Using
the new dynamic programming approach we have proven optimality for the
pentomino and hexomino fence, and the abstraction of the reference
polygon allowed us to characterize all optimal solutions. Furthermore,
we now have good bounds even for very large $n$-omino fences and many
results on the generalization to three dimensions.

We would be delighted to see further research on these fence problems and hope that the bounds on unsolved fences can be further improved with new ideas. Another interesting avenue is the application of the dynamic programming approach to non-square grids. In our previous paper~\cite{previous-german} we used the ILP encoding to find optimal fences for the 7 tetrahexes and the 12 hexiamonds respectively. We believe that the multi-step method laid out in Section~\ref{sec:hexomino} can determine the currently unknown optimal fence area for the 22 pentahexes and 24 heptiamonds respectively.

\paragraph{AI Usage Statement.}
We have developed all fence algorithms, in particular the
dynamic programming approach including all its optimizations and the
various approaches for the 3D walls, entirely without AI
assistance. This includes the programs that create good fences or walls and that establish rigorous upper bounds.

In many other ways, the creation of this paper was supported by
the AI assistant
Claude of Anthropic.
In particular, the AI assistant translated parts of
Section~\ref{sec:hexomino} from an earlier German draft,
created and improved
figures, and was instrumental in keeping a consistent notation.
It provided
optimized programs for running
and evaluating the computational experiments in
Section~\ref{sec:limit}, and it identifed literature on isoperimetric
problems.
The AI assistant helped to formulate the precise problem specification
in Sections~\ref{sec:precise-def} and~\ref{sec:3d}, and it
provided drafts for the proofs of
Lemmas~\ref{lem:area-decomposition}, \ref{lem:integer-vector},
\ref{non-overlap} and~\ref{lem:simple-walk}.

Crucially, the AI assistant pointed out that we had
forgotten to consider the
permutations of parallel reference vectors in
\autoref{sec:all-6}, and it provided the program 
to count the 500 permuted copies of the 244 reference polygons.
It also discovered a mistake in treating the two families 
with six half-integer vectors that can be used with
rotation $\pm 90^\circ$:
Originally we treated them as 4 ``convex'' corner vectors
with $+90^\circ$ (plus two
reflex ones
with $-90^\circ$).
As discussed in \autoref{sec:all-6}
the correct partition is $5+1$.
All these errors have been fixed.
The authors take full responsibility for this work.

\paragraph{Funding.}
ALR acknowledges funding from a postdoctoral fellowship of the Fonds de Recherche du Québec -- Nature et Technologies (ref.~326641) during the course of this project.

\paragraph{Acknowledgments.}

The initial inspiration for this project came from Dr.\ Erich Fuchs, who sparked Mia Müßig's enthusiasm for this problem many years ago when she was still a school student. We thank Corbin Masur, Frederik Safner, and Lennard Scharein, who helped with the maximum wall search  in Section~\ref{sec:3d}. We would also like to thank Jonas Handwerker for proofreading the initial draft of Section~\ref{sec:pentomino}. Finally, we thank the participants of the Fence Challenge workshops.

\addcontentsline {toc}{section}{\numberline {--}References}%
\printbibliography
\appendix

\section{Optimizing the total enclosed area}
\label{sec:total-area}

For fences whose enclosed cells form several regions,
Definition~\ref{def:fence} counts only the \emph{largest} region. This
was convenient, because it gave us a cyclic order of the pieces around
the boundary of this region (Section~\ref{sec:pieces-around}), and
thereby a grip on optimizing the fence along this order by dynamic
programming. An alternative, equally natural definition counts all
enclosed cells: the \emph{total enclosed area} $A_{\mathrm{tot}}(\Phi)$ of a fence~$\Phi$
is the number of its enclosed cells. For the arrangement of
\autoref{fig:weird}, this is the hatched plus the cross-hatched area,
$51+21+21+5+2+1=101$ cells.
It seems unlikely that a fence whose enclosed cells fall apart into
several regions could be optimal, at least when all pieces have the
same size. For pentominoes, we will confirm this: the largest total
enclosed area is again~$128$, and it is attained by exactly the same
fences (Corollary~\ref{cor:total-area}).

Holes inside pieces, like those of \textit{C} and~\textit{H} in
\autoref{fig:weird}, count towards $A_{\mathrm{tot}}(\Phi)$. Free cells that are 8-adjacent
only to cells of a single piece, like the holes of \textit{C}
and~\textit{H}, are enclosed in every placement. They add a constant to
the total area, and they do not affect the optimization.
This is not true for every hole, however: For example, the heptomino with a hole
leaves a diagonal gap next to its hole; this hole can either be connected
to the outside region or be part of the enclosed area, depending on how
the piece is built into the fence.
Pentominoes and hexominoes have no holes.

\paragraph{The walk around the outside.}
We construct a walk along the boundary of the fence, as in
Section~\ref{sec:pieces-around}, except that we walk along the outer boundary.

Like a region, the outside of a fence is 8-connected. Its complement,
which consists of the fence cells and the enclosed cells, is finite
and 4-connected: the fence~$\Phi$ is 4-connected, and from every
enclosed cell, a horizontal row of enclosed cells leads to a fence
cell. Hence, as in Section~\ref{sec:pieces-around}, the boundary of
the outside forms a single closed walk. We orient it clockwise around
the fence, so that the outside lies on the left side and the fence on
the right side. At a vertex where two cells of the outside touch
diagonally, like the vertex where \textit{G} and~\textit{R} touch in
\autoref{fig:weird}, the walk turns right and goes through the gap.
We record the fence cells along this
walk; 
at a left turn, we also record the 
diagonally opposite cell that meets the walk only at a vertex.
The recorded
cells are the fence cells that are 8-adjacent to the outside, and
consecutive recorded cells are 4-adjacent. Replacing every recorded
cell by the piece that contains it and merging consecutive repetitions,
we obtain the cyclic sequence of the \emph{pieces around the outside}.

In \autoref{fig:weird}, the pieces around the outside, in clockwise
order starting at the top left, are
\textit{ABCDEFQRQFGHIEJKLMNOCP}. The pieces \textit{C}, \textit{E},
\textit{F}, and~\textit{Q} occur twice, and \textit{S} and~\textit{T}
do not occur at all.

\begin{definition}\label{def:simple-fence}
  A fence~$\Phi$ is \emph{simple} if it consists of a single piece, or if
  its pieces can be arranged in a cyclic sequence
  $P_1,P_2,\ldots,P_k$ with $k\ge2$ (indices taken modulo~$k$) such
  that there is a Jordan curve~$Z$ with the following
  properties:
  \begin{itemize}
  \item $Z$ runs in the interior of the fence, that is, in the
    interior of the union~$\Phi$ of the pieces (regarded as closed
    subsets of the plane);
  \item $Z$ visits every piece exactly once, in the given order: for
    each~$i$, the intersection $Z\cap P_i$ is an arc that lies in the
    interior of~$P_i$, except for its endpoints. One endpoint is an
    interior point of an edge that $P_i$ shares with~$P_{i-1}$, and
    the other one is an interior point of an edge that $P_i$ shares
    with~$P_{i+1}$;
  \item all enclosed cells of~$\Phi$ lie inside~$Z$.
  \end{itemize}
\end{definition}
  The edges through which $Z$ passes from one piece to the next are
  the 
  {attachment edges} from Section~\ref{sec:pentomino}.

\begin{figure}[htb]
  \centering
  \begin{tikzpicture}[scale=0.45, y=-1cm]
    \foreach \tx/\ty in {
      5/4, 6/4, 7/4, 4/5, 5/5, 6/5, 7/5, 3/6, 4/6, 5/6, 6/6, 7/6,
      3/7, 4/7, 5/7, 6/7, 7/7, 3/8, 4/8, 5/8, 6/8, 7/8, 4/9, 5/9,
      6/9, 7/9, 5/10, 6/10, 7/10
    }
      \fill[pattern=north east lines, pattern color=black!60] (\tx,\ty) rectangle ++(1,1);
    \foreach \tx/\ty in {
      9/4, 10/4, 11/4, 9/5, 10/5, 11/5, 12/5, 9/6, 10/6, 11/6, 12/6, 13/6,
      9/7, 10/7, 11/7, 12/7, 13/7, 9/8, 10/8, 11/8, 12/8, 13/8, 9/9, 10/9,
      11/9, 12/9, 9/10, 10/10, 11/10
    }
      \fill[pattern=north west lines, pattern color=black!60] (\tx,\ty) rectangle ++(1,1);
    \foreach \x/\y in {13/5, 14/5, 14/6, 15/6, 14/7}
      \draw[fill=colorF!45] (\x,\y) rectangle ++(1,1);
    \foreach \x/\y in {8/3, 8/4, 8/5, 8/6, 8/7}
      \draw[fill=colorI!45] (\x,\y) rectangle ++(1,1);
    \foreach \x/\y in {8/8, 8/9, 8/10, 8/11, 9/11}
      \draw[fill=colorL!45] (\x,\y) rectangle ++(1,1);
    \foreach \x/\y in {12/10, 13/10, 10/11, 11/11, 12/11}
      \draw[fill=colorN!45] (\x,\y) rectangle ++(1,1);
    \foreach \x/\y in {5/11, 6/11, 7/11, 5/12, 6/12}
      \draw[fill=colorP!45] (\x,\y) rectangle ++(1,1);
    \foreach \x/\y in {10/1, 10/2, 9/3, 10/3, 11/3}
      \draw[fill=colorT!45] (\x,\y) rectangle ++(1,1);
    \foreach \x/\y in {12/2, 13/2, 12/3, 12/4, 13/4}
      \draw[fill=colorU!45] (\x,\y) rectangle ++(1,1);
    \foreach \x/\y in {2/10, 3/10, 4/10, 4/11, 4/12}
      \draw[fill=colorV!45] (\x,\y) rectangle ++(1,1);
    \foreach \x/\y in {3/4, 4/4, 2/5, 3/5, 2/6}
      \draw[fill=colorW!45] (\x,\y) rectangle ++(1,1);
    \foreach \x/\y in {14/8, 13/9, 14/9, 15/9, 14/10}
      \draw[fill=colorX!45] (\x,\y) rectangle ++(1,1);
    \foreach \x/\y in {5/2, 4/3, 5/3, 6/3, 7/3}
      \draw[fill=colorY!45] (\x,\y) rectangle ++(1,1);
    \foreach \x/\y in {1/7, 2/7, 2/8, 2/9, 3/9}
      \draw[fill=colorZ!45] (\x,\y) rectangle ++(1,1);
    \draw[line width=1.5pt, red!75!black, rounded corners=3pt] (12.5,3.5) -- (12.5,4.5) --
      (13.5,4.5) -- (13.5,5.5) -- (14.5,5.5) -- (14.5,10.5) -- (12.5,10.5) -- (12.5,11.5) --
      (6.5,11.5) -- (6.5,12.5) -- (4.5,12.5) -- (4.5,10.5) -- (2.5,10.5) -- (2.5,5.5) --
      (3.5,5.5) -- (3.5,4.5) -- (4.5,4.5) -- (4.5,3.5) -- cycle;
    \node[red!75!black, fill=white, inner sep=1pt] at (3.75,8.5) {$Z$};
  \end{tikzpicture}
  \caption{A simple fence of the $12$ pentominoes.
The I-pentomino and the L-pentomino form a wall that separates the enclosed
area into two regions
  of $29$~cells each (hatched in different directions). 
  The total
    enclosed area is~$58$. The Jordan curve~$Z$ visits the pieces in the cyclic
    order in which they occur around the outside, and it passes from one piece
    to the next through an attachment edge.}
  \label{fig:several-parts}
\end{figure}
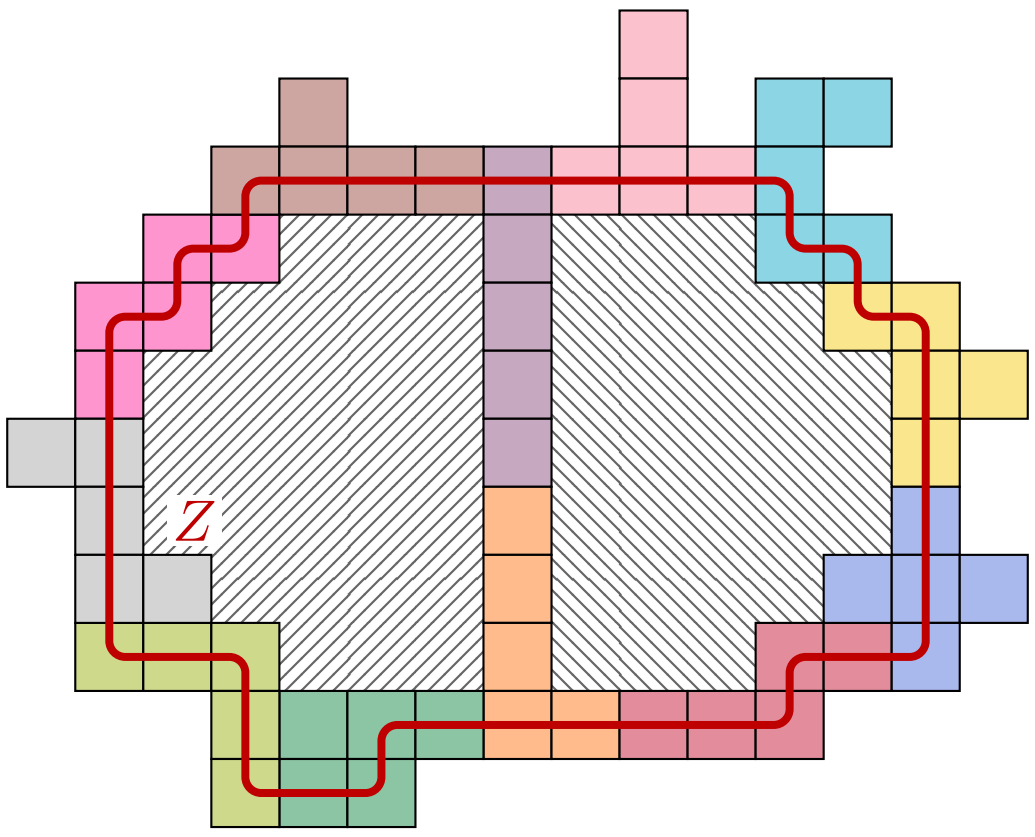

Since $Z$ runs in the interior of the fence, a free cell inside~$Z$
cannot be 8-connected to a free cell outside~$Z$. Hence the enclosed
cells of a simple fence are precisely the free cells inside~$Z$.
However, they can form several regions, as in
\autoref{fig:several-parts}.

\begin{lemma}\label{lem:decomposition}
  Let $\Phi$ be a fence.
  \begin{enumerate}
  \item[(a)] If every piece of~$\Phi$ occurs exactly once in the cyclic
    sequence of pieces around the outside, then $\Phi$ is simple.
  \item[(b)] If some pieces of~$\Phi$ do not occur in this sequence, then
    the pieces that do occur form a fence~$\Phi'$ with $A_{\mathrm{tot}}(\Phi')>A_{\mathrm{tot}}(\Phi)$.
  \item[(c)] If some piece~$P$ occurs more than once in this sequence,
    then $\Phi$ can be decomposed into two fences $\Phi_1$ and~$\Phi_2$, each
    consisting of at least two pieces, that share the piece~$P$ and no
    other piece, such that every enclosed cell of~$\Phi$ is an enclosed
    cell of $\Phi_1$ or of~$\Phi_2$. In particular,
    $A_{\mathrm{tot}}(\Phi)\le A_{\mathrm{tot}}(\Phi_1)+A_{\mathrm{tot}}(\Phi_2)$.
  \end{enumerate}
  Consequently, if $\Phi$ is not simple, then (b) or (c) applies.
\end{lemma}

\begin{figure}[htb]
  \centering
  \begin{tikzpicture}[scale=0.45, y=-1cm]
    \foreach \tx/\ty in {
      3/3, 4/3, 5/3, 2/4, 3/4, 4/4, 5/4, 2/5, 3/5, 4/5, 5/5, 2/6,
      3/6, 4/6, 5/6, 3/7, 4/7, 5/7
    }
      \fill[pattern=north east lines, pattern color=black!60] (\tx,\ty) rectangle ++(1,1);
    \foreach \tx/\ty in {
      9/3, 10/3, 11/3, 12/3, 9/4, 10/4, 11/4, 12/4, 13/4, 9/5, 10/5, 11/5,
      12/5, 13/5, 9/6, 10/6, 11/6, 12/6, 13/6, 9/7, 10/7, 11/7, 12/7
    }
      \fill[pattern=north east lines, pattern color=black!60] (\tx,\ty) rectangle ++(1,1);
    \foreach \x/\y in {3/1, 2/2, 3/2, 4/2, 2/3}
      \draw[fill=colorF!45] (\x,\y) rectangle ++(1,1);
    \foreach \x/\y in {1/3, 1/4, 1/5, 1/6, 1/7}
      \draw[fill=colorI!45] (\x,\y) rectangle ++(1,1);
    \foreach \x/\y in {8/6, 8/7, 8/8, 8/9, 9/9}
      \draw[fill=colorL!45] (\x,\y) rectangle ++(1,1);
    \foreach \x/\y in {6/6, 6/7, 5/8, 6/8, 5/9}
      \draw[fill=colorN!45] (\x,\y) rectangle ++(1,1);
    \foreach \x/\y in {9/8, 10/8, 11/8, 10/9, 11/9}
      \draw[fill=colorP!45] (\x,\y) rectangle ++(1,1);
    \foreach \x/\y in {2/7, 2/8, 3/8, 4/8, 2/9}
      \draw[fill=colorT!45] (\x,\y) rectangle ++(1,1);
    \foreach \x/\y in {14/3, 15/3, 14/4, 14/5, 15/5}
      \draw[fill=colorU!45] (\x,\y) rectangle ++(1,1);
    \foreach \x/\y in {8/2, 9/2, 10/2, 8/3, 8/4}
      \draw[fill=colorV!45] (\x,\y) rectangle ++(1,1);
    \foreach \x/\y in {14/6, 13/7, 14/7, 12/8, 13/8}
      \draw[fill=colorW!45] (\x,\y) rectangle ++(1,1);
    \foreach \x/\y in {7/4, 6/5, 7/5, 8/5, 7/6}
      \draw[fill=colorX!75] (\x,\y) rectangle ++(1,1);
    \foreach \x/\y in {6/1, 5/2, 6/2, 6/3, 6/4}
      \draw[fill=colorY!45] (\x,\y) rectangle ++(1,1);
    \foreach \x/\y in {11/1, 11/2, 12/2, 13/2, 13/3}
      \draw[fill=colorZ!45] (\x,\y) rectangle ++(1,1);
    \draw[line width=1.5pt, blue!70!black, dashed, rounded corners=8pt]
      (7.5,5.5) -- (7.5,0.5) -- (0.5,0.5) -- (0.5,10.5) -- (7.5,10.5) -- cycle;
    \node[blue!70!black, fill=white, inner sep=1pt] at (-0.1,5.5) {$\gamma$};
    \node[fill=white, inner sep=1pt] at (4,5.5) {$\Phi_1$};
    \node[fill=white, inner sep=1pt] at (11.5,5.5) {$\Phi_2$};
  \end{tikzpicture}
  \caption{A fence of the $12$ pentominoes in which the X-pentomino
    (dark) occurs twice around the outside: once from above, and once from below.
    The dashed closed curve~$\gamma$ runs through the outside and through~X. It cuts
    the fence into two simple fences that share the X-pentomino: $\Phi_1$, with six
    pieces, encloses $18$~cells, and $\Phi_2$, with seven pieces, encloses $23$~cells.}
  \label{fig:cut-into-subfences}
\end{figure}
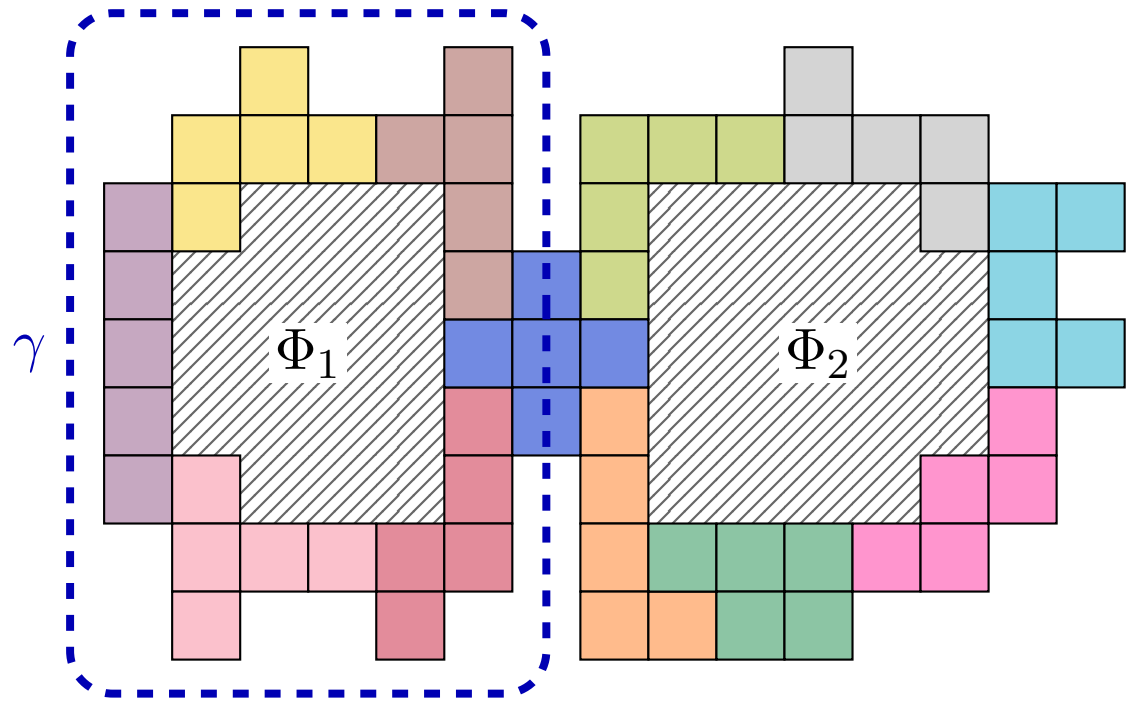

\begin{proof}
  Let $\omega$ be the walk around the outside.

  (a) Let $Z$ be the closed curve that runs parallel to~$\omega$ at
  offset
  ~$\frac14$ on its right side. Along an edge of~$\omega$, $Z$
  runs through the fence cell on the right side of this edge. At a
  right turn of~$\omega$, it stays in this cell, and at a left turn, it
  goes around the vertex through the diagonally opposite cell. Thus,
  $Z$ passes through the recorded cells in the order in which they are
  recorded, and it passes from one recorded cell to the next through an
  interior point of their common edge. In particular, $Z$ runs in the
  interior of the fence. The walk~$\omega$ uses every edge only once,
  and at a vertex that it passes twice, it turns right both times.
  Therefore, $Z$ is a Jordan curve. If every piece occurs exactly once
  around the outside, then $Z$ visits every piece exactly once, as
  required in Definition~\ref{def:simple-fence}. Finally, $Z$ surrounds
  all points of the complement of the outside that have distance more
  than~$\frac14$ from~$\omega$, because $\omega$ is the boundary of
  this complement. An enclosed cell is not 8-adjacent to
  the outside, and so it has distance at least~$1$ from~$\omega$.
  Hence, all enclosed cells lie inside~$Z$, and $\Phi$ is simple.

  (b) Since consecutive recorded cells are 4-adjacent, the pieces
  around the outside form a 4-connected set, and hence a fence~$\Phi'$.
  Removal of the pieces that are not used in $\Phi'$ increases
  the enclosed area, and therefore
  $A_{\mathrm{tot}}(\Phi')>A_{\mathrm{tot}}(\Phi)$.

  (c)
  We cut the fence into two parts along
  a Jordan curve $\gamma$ through $P$ and the outside region, see
  \autoref{fig:cut-into-subfences}.
More precisely,
  we pick two occurrences of~$P$ in the cyclic
  sequence. Then the sequence has the form $P,\sigma_1,P,\sigma_2$,
  where $\sigma_1$ and~$\sigma_2$ are nonempty.
  Let $a$ be a cell of~$P$ that is recorded at the
  first occurrence, and let $o_1$ be a cell of the outside that is
  8-adjacent to~$a$ (the cell on the left side of the corresponding
  edge of~$\omega$, or the corner cell of the outside at a left turn).
  Similarly, choose $c$ and~$o_2$ for the second occurrence.
  We form the Jordan curve~$\gamma$ from four parts: a straight segment from the center
  of~$a$ to the center of~$o_1$, a path from $o_1$ to~$o_2$ through
  the centers of cells of the outside (which is 8-connected), a
  straight segment from $o_2$ to~$c$, and a path from $c$ back to~$a$
  through the centers of cells of~$P$ (which is 4-connected).
  The fences $\Phi_1$ and $\Phi_2$ will then consist of $P$ together
  with the pieces inside $\gamma$ and outside $\gamma$, respectively.
  The remainder of the proof will establish formally that
  this construction achieves the claimed goals.

First, 
 $\gamma$ is indeed a Jordan curve,
  if we
  choose shortest paths in the above construction:
  In a shortest
  path, two cells that are not consecutive are not 8-adjacent, and so
  the path cannot cross itself.

  The curve~$\gamma$ crosses $\omega$
  exactly twice, between $a$ and~$o_1$ and between $o_2$ and~$c$.
  Elsewhere, it can touch~$\omega$ only at vertices where it passes
  diagonally between two cells of the outside. So, the part
  of~$\omega$ that runs along~$\sigma_1$ lies on one side of~$\gamma$,
  and the part that runs along~$\sigma_2$ lies on the other side. Every
  piece other than~$P$ meets $\gamma$ at most in isolated vertices, and
  hence it lies entirely inside or entirely outside~$\gamma$. In
  particular, the pieces of~$\sigma_1$ other than~$P$ lie on one side of~$\gamma$, say
  inside, and the pieces of~$\sigma_2$ other than~$P$
  lie outside; thus
  $\sigma_1$ and $\sigma_2$ have no piece other than~$P$ in common.

  Let $\Phi_1$ consist of~$P$ and the pieces inside~$\gamma$, and let
  $\Phi_2$ consist of~$P$ and the pieces outside~$\gamma$. Each of them
  contains at least two pieces. A 4-connected path of fence cells from
  a piece inside~$\gamma$ to~$P$, drawn through the centers of its
  cells, can leave the inside of~$\gamma$ only through~$P$, since
  $\gamma$ runs only through~$P$ and through the outside. Hence $\Phi_1$ is 4-connected, and similarly,
  $\Phi_2$ is 4-connected; thus, both are fences.

  Finally, let $e$ be an enclosed cell of~$\Phi$. It does not
  meet~$\gamma$, and we assume without loss of generality that it lies inside~$\gamma$.
  Suppose that $e$ is not enclosed in~$\Phi_1$. Then there is an
  8-connected path of free cells of~$\Phi_1$ from $e$ to a cell far away,
  outside the bounding box of~$\Phi$. Let $g$ be the first cell on this
  path that is not an enclosed cell of~$\Phi$; it exists because the last
  cell of the path is not enclosed in~$\Phi$. It is 8-adjacent
  to an enclosed cell of~$\Phi$, so it does not belong to the outside
  of~$\Phi$. Since $g$ is free in~$\Phi_1$, it belongs to a piece of~$\Phi_2$
  other than~$P$, and thus it lies outside~$\gamma$. But the path from
  $e$ to~$g$, drawn through the centers of its cells, does not
  meet~$\gamma$: it
  never enters~$P$, and up to~$g$, it runs through enclosed cells
  of~$\Phi$, which are not even 8-adjacent to the outside of~$\Phi$. This
  contradiction shows that $e$ is enclosed in~$\Phi_1$.
\end{proof}

In the
example of \autoref{fig:cut-into-subfences}, the total enclosed area is the sum of the total areas of the
two subfences. This is not true in general, because a hole inside the
shared piece~$P$ is counted in both subfences.

For a more elaborate example, we look again at \autoref{fig:weird}.
Although this arrangement is not a fence, the walk around its outside
is defined in the same way, and we have seen that the pieces
\textit{S} and~\textit{T} do not occur. The other $18$ pieces form a
fence of total area~$110$, in which the cells of \textit{S} and
\textit{T} become free and join the regions that surround them.
Cutting this fence through the pieces \textit{C}, \textit{E},
\textit{F}, and~\textit{Q}, which occur twice, in any order, we obtain
five simple fences, consisting of the pieces \textit{ABCP},
\textit{CDEJKLMNO}, \textit{EFGHI}, \textit{FQ}, and~\textit{QR},
respectively. The last two, in the upper right corner, enclose
nothing. The other three have total areas $22$, $60$, and~$29$. Their
sum, $111$, exceeds $110$ by one, because the hole of~\textit{C} is
counted twice.

Lemma~\ref{lem:decomposition} reduces every fence to simple fences.
Since the area that can be enclosed by $i$ pieces grows roughly
quadratically with $i$,
splitting a fence into two parts that share a
piece is expensive. The following lemma makes this precise. A fence
\emph{built from pieces of~$\mathcal{S}$} uses each piece of~$\mathcal{S}$ at most once.

\begin{lemma}\label{lem:subadditive}
  Let $\mathcal{S}$ be a family of $N$ polyomino pieces, and let
  $u_1,u_2,\ldots,u_N$ be numbers with the following properties.
  \begin{enumerate}
  \item[(i)] For $1\le i\le N$, $u_i$ is an upper bound on the total
    enclosed area of every simple fence that is built from $i$ pieces
    of~$\mathcal{S}$.
  \item[(ii)] $u_i+u_j\le u_{i+j-1}$ for all $i,j\ge 2$ with
    $i+j\le N+1$.
  \item[(iii)] $u_1\le u_2\le\cdots\le u_N$.
  \end{enumerate}
  Then:
  \begin{enumerate}
  \item[(a)] Every fence~$\Phi$ that is built from $i$ pieces of~$\mathcal{S}$,
    $1\le i\le N$, satisfies $A_{\mathrm{tot}}(\Phi)\le u_i$. In particular,
    $u_N$ is an upper bound on the total area enclosed by any fence
    built from pieces of~$\mathcal{S}$.
  \item[(b)] If, in addition, the inequalities in~(ii) are strict,
    then every fence that is built from all $N$ pieces of~$\mathcal{S}$ and is
    not simple satisfies $A_{\mathrm{tot}}(\Phi)<u_N$.
  \end{enumerate}
\end{lemma}

The inequalities $u_2\le u_3\le\cdots\le u_N$ in~(iii) follow from~(ii)
when $u_2\ge0$, by taking $j=2$; only $u_1\le u_2$ must be checked
separately. For larger polyominoes, which may enclose holes inside the
pieces, it would likely be advantageous to take those holes out of the game
and adjust the lemma. Then, $u_1$ becomes zero, and double-counting of
those holes is avoided.

\begin{proof}
  (a) We use induction on~$i$. If $\Phi$ is simple, in particular if
  $i=1$, then $A_{\mathrm{tot}}(\Phi)\le u_i$ by~(i). Otherwise, by
  Lemma~\ref{lem:decomposition}, case (b) or~(c) of this lemma
  applies. In case~(b), the fence~$\Phi'$ consists of $i'<i$ pieces, and
  by induction and~(iii),
  $A_{\mathrm{tot}}(\Phi)<A_{\mathrm{tot}}(\Phi')\le u_{i'}\le u_i$.
  In case~(c), the fences $\Phi_1$ and~$\Phi_2$ consist of $i_1\ge2$ and
  $i_2\ge2$ pieces with $i_1+i_2=i+1$, and hence $i_1,i_2<i$. By
  induction and~(ii),
  \begin{equation}\label{eq:subadditive}
    A_{\mathrm{tot}}(\Phi)\le A_{\mathrm{tot}}(\Phi_1)+A_{\mathrm{tot}}(\Phi_2)\le u_{i_1}+u_{i_2}\le u_{i_1+i_2-1}=u_i .
  \end{equation}
  The second statement holds because $u_i\le u_N$ by~(iii).

  (b) Let $\Phi$ be a fence that is built from all $N$ pieces and is not
  simple. In case~(b) of Lemma~\ref{lem:decomposition}, we get
  $A_{\mathrm{tot}}(\Phi)<A_{\mathrm{tot}}(\Phi')\le u_{i'}\le u_N$ as before, by part~(a). In case~(c),
  the last inequality in~\eqref{eq:subadditive} is strict, and
  $A_{\mathrm{tot}}(\Phi)<u_N$.
\end{proof}

To apply Lemma~\ref{lem:subadditive}, we need upper bounds for simple
fences. For pieces without holes, these are provided by the dynamic
program of Section~\ref{sec:pentomino}, since the total area of a
simple fence is the area of a closed walk of the kind that the
dynamic program considers:

\begin{lemma}\label{lem:simple-walk}
  Let $\Phi$ be a simple fence with $k\ge2$ pieces whose boundaries are
  simple closed curves, and let $Z$ be a curve as in
  Definition~\ref{def:simple-fence}. We orient $Z$ counterclockwise
  and number the pieces $P_1,\ldots,P_k$ in the order in which $Z$
  visits them.
  \begin{enumerate}
  \item[(a)] If $k=2$ and $Z$ crosses the same edge twice, then
    $A_{\mathrm{tot}}(\Phi)=0$.
  \item[(b)] Otherwise, for each~$i$, let $\beta_i$ be the part of
    the boundary of~$P_i$ that lies inside~$Z$ between the attachment
    edge through which $Z$ enters~$P_i$ and the attachment edge
    through which $Z$ leaves~$P_i$, excluding these two edges. Then
    $\beta_1,\beta_2,\ldots,\beta_k$ form a closed walk~$\beta$, and
    the area enclosed by~$\beta$, as computed by the Shoelace
    formula~\eqref{eq:shoelace}, is~$A_{\mathrm{tot}}(\Phi)$.
  \end{enumerate}
\end{lemma}

\begin{proof}
  Every piece~$P_i$ is a topological disk, and the arc $Z\cap P_i$
  cuts it into two parts. Let $P_i^+$ be the part inside~$Z$.

  (a) In this case, the arc $Z\cap P_i$, for $i=1,2$, and the segment
  between the two crossing points on the common edge bound a part
  of~$P_i$. These two parts lie on different sides of the common edge,
  and together they form a disk bounded by~$Z$. Hence they are the
  parts $P_1^+$ and~$P_2^+$, they cover the inside of~$Z$, there is no
  free cell inside~$Z$, and $A_{\mathrm{tot}}(\Phi)=0$.

  (b) The curve~$Z$ meets each attachment edge only at one point, where
  it crosses the edge. Hence exactly one endpoint of the edge lies
  inside~$Z$, and we call the segment from the crossing point to this
  endpoint the \emph{inner part} of the edge. The boundary of~$P_i^+$
  consists of the arc $Z\cap P_i$, the inner parts of the two
  attachment edges of~$P_i$, and~$\beta_i$. Thus, $\beta_i$ leads from
  the endpoint inside~$Z$ of the attachment edge where $Z$
  enters~$P_i$ to the endpoint inside~$Z$ of the attachment edge where
  $Z$ leaves~$P_i$. The latter edge is the one where $Z$
  enters~$P_{i+1}$, so consecutive sections fit together to a closed
  walk~$\beta$.

  When we replace the arc $Z\cap P_i$ by the path along the boundary
  of~$P_i^+$ that consists of the inner part of an attachment edge,
  $\beta_i$, and the inner part of the other attachment edge, the
  winding number around a point decreases by~$1$ if the point lies in
  the interior of~$P_i^+$, and otherwise it does not change. After
  doing this for all~$i$, the inner parts of the attachment edges
  cancel, and we obtain~$\beta$. Hence $\beta$ winds once around every
  point inside~$Z$ that is not covered by a piece, and zero times
  around every other point that is not on~$\beta$.
  Since the Shoelace formula integrates the winding number, it yields
  the number of free cells inside~$Z$, which are precisely the
  enclosed cells of~$\Phi$.
\end{proof}

The walk~$\beta$ runs counterclockwise around the enclosed cells, with
the pieces on its right side, and $\beta_i$ is a section of the boundary
of~$P_i$ between an incoming and an outgoing attachment edge, as
described in Section~\ref{sec:boundary}; it may degenerate to a single point, like the
section of~\textit{K} in \autoref{fig:weird}. Consecutive pieces share
the attachment edge between them. So, $\beta$ is one of the closed
boundaries whose Shoelace area is maximized by the dynamic program of
Section~\ref{sec:pentomino}, and for pentominoes, $A_{\mathrm{tot}}(\Phi)$ is
bounded by the entry $A_{\max}(k)$ of Table~\ref{tab:max-subset}.

\begin{corollary}\label{cor:total-area}
  The largest total enclosed area of a fence built from pentominoes is
  $128$. Every such fence with total enclosed area~$128$ uses all $12$
  pentominoes, is simple, and has a single region of $128$~cells.
  Thus, the optimal fences with respect to the total enclosed area are
  the same as the optimal fences with respect to the area of
  Definition~\ref{def:fence}, namely the $5760$ solutions found in
  Section~\ref{sec:dp-unrestricted}.
\end{corollary}

\begin{proof}
  We apply Lemma~\ref{lem:subadditive} to the family~$\mathcal{S}$ of the
  $N=12$ pentominoes, with $u_i=A_{\max}(i)$ from
  Table~\ref{tab:max-subset}. By the construction in
  Section~\ref{sec:dp-unrestricted}, $A_{\max}(i)$ is the largest area,
  computed by the Shoelace formula, of a closed boundary that the
  dynamic program can form from $i$ distinct pentominoes, maximized
  over all choices of the pentominoes, of their cyclic order, and of
  their incoming and outgoing attachment edges, regardless of whether
  the pieces overlap.
  By Lemma~\ref{lem:simple-walk} and the remark following it,
  $A_{\max}(i)$ is an upper bound on the total enclosed area of every
  simple fence of $i$ pentominoes, for $i\ge2$ (in the case of
  Lemma~\ref{lem:simple-walk}(a), trivially). A single pentomino
  encloses nothing, and $A_{\max}(1)=0$. Thus, condition~(i) holds.
  Condition~(iii) is obvious from the table. Condition~(ii) holds even
  with strict inequality, as can be checked in the table for all pairs
  $2\le i\le j$ with $i+j-1\le 12$. For example,
  $A_{\max}(4)+A_{\max}(6)=13+34=47<A_{\max}(9)=78$. The smallest gap
  occurs for $i=j=2$, where $1+1<6=A_{\max}(3)$.

  By Lemma~\ref{lem:subadditive}(a), every fence~$\Phi$ built from
  pentominoes satisfies
  $A_{\mathrm{tot}}(\Phi)\le A_{\max}(12)=128$, and the optimal fences of
  Section~\ref{sec:dp-unrestricted} attain this bound. Now let $\Phi$ be a
  fence with $A_{\mathrm{tot}}(\Phi)=128$.
  Lemma~\ref{lem:subadditive}(a) excludes the case that
 $\Phi$ does not use all  $12$ pentominoes, because
 $
 A_{\max}(11)=111$.
  By Lemma~\ref{lem:subadditive}(b), $\Phi$ is simple. By
  Lemma~\ref{lem:simple-walk}, the walk~$\beta$ of~$\Phi$ is a closed
  boundary of the dynamic program with area~$128$, the optimum.
  Therefore, $\Phi$ is one of the optimal placements that were enumerated
  in Section~\ref{sec:dp-unrestricted}, and all of them consist of a
  single region of $128$ cells.
\end{proof}

The proof gives more: a fence of pentominoes that is not simple
encloses at most $112$ cells in total. Indeed, in case~(b) of
Lemma~\ref{lem:decomposition},
$A_{\mathrm{tot}}(\Phi)<A_{\mathrm{tot}}(\Phi')\le A_{\max}(i')\le A_{\max}(11)=111$,
and in case~(c),
$A_{\mathrm{tot}}(\Phi)\le A_{\max}(i_1)+A_{\max}(i_2)$ with
$i_1+i_2=i+1\le13$, which is at most $A_{\max}(2)+A_{\max}(11)=112$.

For hexominoes, whose boundaries are also simple closed curves,
Lemma~\ref{lem:simple-walk} holds as well. The bounds of
Section~\ref{sec:hexomino} should also apply to the walk~$\beta$ of a
simple fence that uses all $35$ hexominoes. However, $\beta$ may touch
itself, and so the rotation condition
$D_1+\dots+D_{35}=360^\circ$, which was derived for the boundary of a
region, would have to be justified with the help of the Jordan
curve~$Z$ (for example, by applying the theorem on turning tangents to
the boundary of each part~$P_i^+$). To show that no hexomino fence encloses a
total area of $1598$ or more, it would suffice to derive upper bounds~$u_i$ for
simple fences of $i<35$ hexominoes that satisfy the conditions of
Lemma~\ref{lem:subadditive} together with $u_{35}=1597$, for example
with the simpler techniques of Section~\ref{sec:ub-hex}. Since
nonsimple fences pay a substantial penalty, rather weak bounds should
suffice. However, we have not carried out this proof.

\section{Hexomino area distribution of random orders}
\label{sec:hex-dist}

\begin{figure}[h!]
    \centering
    \begin{tikzpicture}
        \begin{axis}[
            width=15cm,
            height=8cm,
            ybar,
            bar width=0.7,
            ymin=0,
            xmin=1455,
            xmax=1565,
            enlarge x limits=false,
            xtick={
            1460,
            1470,1480,1490,
            1500,1510,1520,1530,1540,1550,
            1560
            },
            xticklabel style={
                /pgf/number format/1000 sep={},
                font=\tiny,
                rotate=45,
                anchor=east
            },
            ylabel style={yshift=0.5em},
            yticklabel style={text width=3em, align=right},
            ytick scale label code/.code={},
            ymajorgrids=true,
            grid style=dashed,
        ]
          \addplot+[fill=blue, draw=none] coordinates {
            (1458, 3) (1459, 2) (1460, 1) (1461, 6) (1462, 8) (1463, 20) (1464, 37) (1465, 64)
            (1466, 125) (1467, 161) (1468, 299) (1469, 412) (1470, 628) (1471, 869) (1472, 1081) (1473, 1243)
            (1474, 1505) (1475, 1671) (1476, 1878) (1477, 1994) (1478, 2110) (1479, 2310) (1480, 2297) (1481, 2488)
            (1482, 2521) (1483, 2652) (1484, 2674) (1485, 2586) (1486, 2747) (1487, 2710) (1488, 2834) (1489, 2697)
            (1490, 2745) (1491, 2668) (1492, 2698) (1493, 2623) (1494, 2528) (1495, 2457) (1496, 2471) (1497, 2328)
            (1498, 2366) (1499, 2246) (1500, 2188) (1501, 2073) (1502, 1983) (1503, 1908) (1504, 1860) (1505, 1766)
            (1506, 1643) (1507, 1591) (1508, 1353) (1509, 1360) (1510, 1281) (1511, 1180) (1512, 1102) (1513, 999)
            (1514, 930) (1515, 887) (1516, 775) (1517, 772) (1518, 701) (1519, 602) (1520, 582) (1521, 543)
            (1522, 487) (1523, 453) (1524, 382) (1525, 332) (1526, 335) (1527, 287) (1528, 225) (1529, 228)
            (1530, 183) (1531, 180) (1532, 155) (1533, 112) (1534, 133) (1535, 97) (1536, 98) (1537, 68)
            (1538, 77) (1539, 51) (1540, 50) (1541, 29) (1542, 33) (1543, 14) (1544, 23) (1545, 30)
            (1546, 17) (1547, 8) (1548, 14) (1549, 11) (1550, 13) (1551, 10) (1552, 4) (1553, 3)
            (1554, 2) (1555, 4) (1556, 3) (1557, 2) (1559, 3) (1560, 1) (1561, 1)
          };
        \end{axis}
    \end{tikzpicture}
    \caption{Distribution of the maximum achievable fence area across $10^5$ random hexomino orders.}
    \label{fig:hex-ord}
\end{figure}
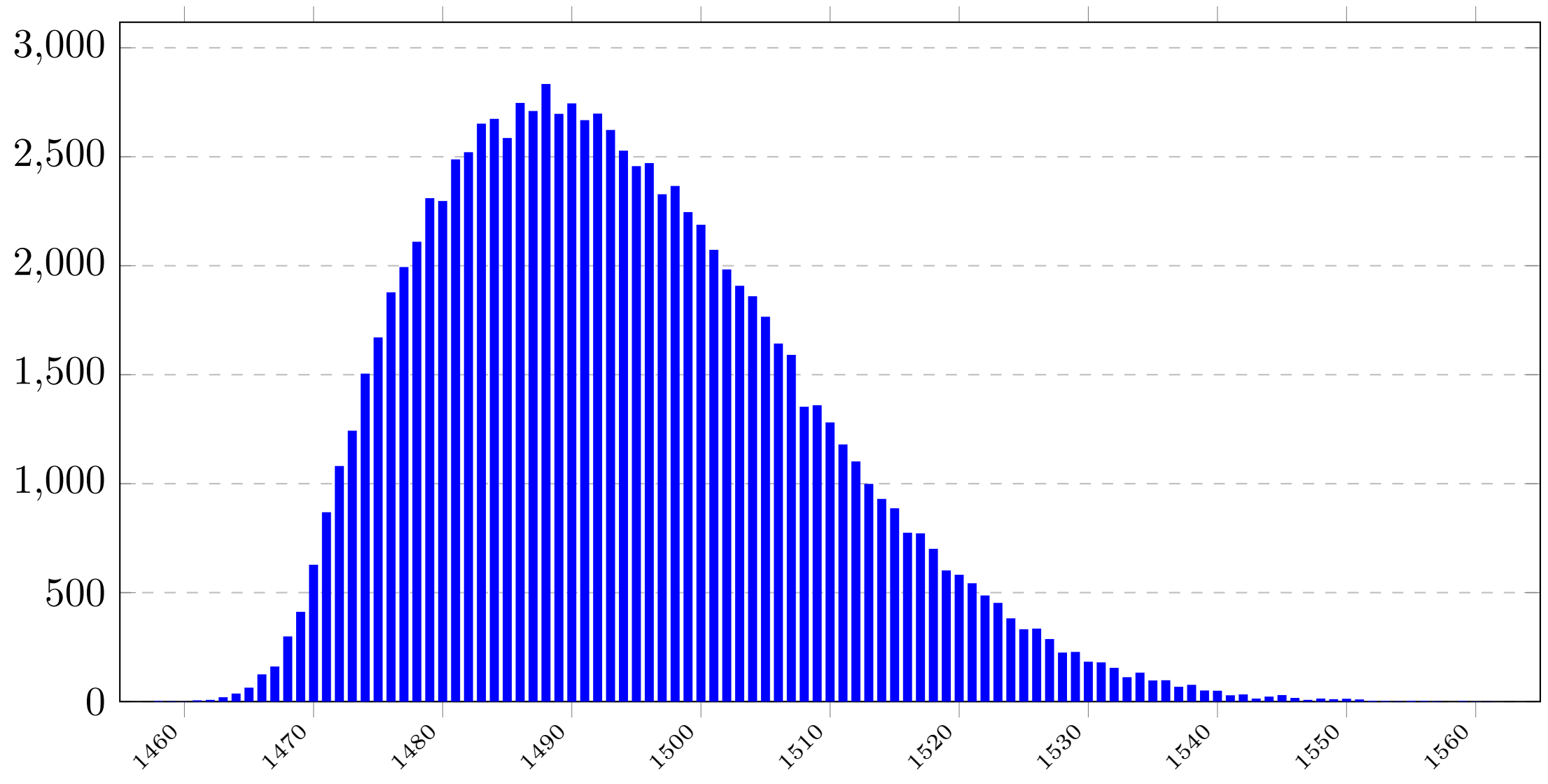

We also ran \hyperref[alg:optimal-ordered-fence]{\textsc{OptimalOrderedFence}} on various orders for the 35 hexominoes. Since here the total number of orders is infeasibly large ($\frac{34!}{2} \approx 1.5 \cdot 10^{38}$), we sampled $10^5$ orders uniformly at random and show the resulting distribution in Figure~\ref{fig:hex-ord}.

The shape closely mirrors the pentomino distribution in Figure~\ref{fig:pent-ord}: a roughly bell-shaped curve with a sharp left tail and a gradual right tail. The smallest maximum hexomino fence area was $1458$, which we observed three times among the samples.

\begin{figure}[h!]
    \centering 
    \includegraphics[width=0.48\textwidth,angle=270]{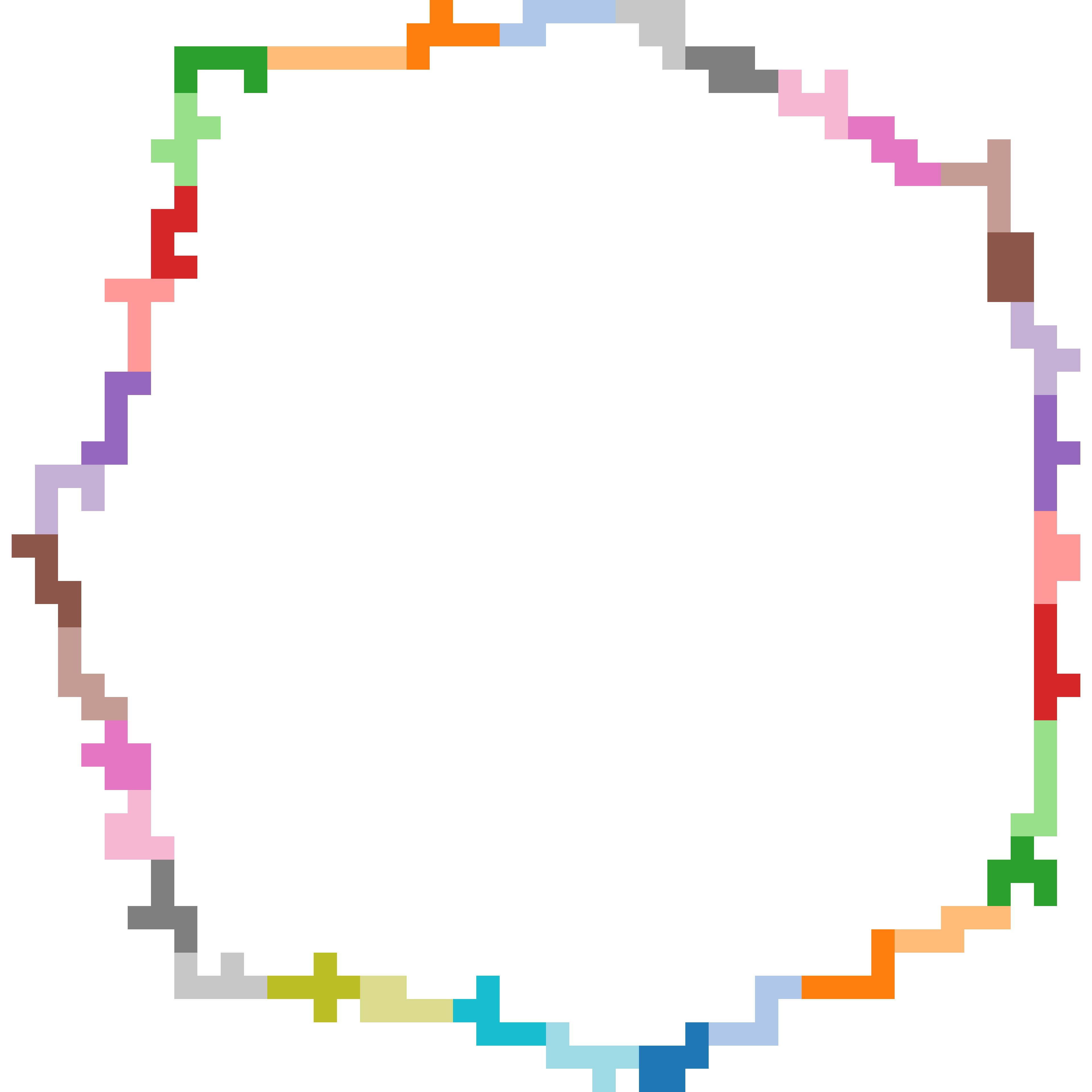}
    \caption{A hexomino order whose best fence has area $1448$, the smallest maximum among all orders that we found.}
    \label{fig:hex-worst}
\end{figure}

We also experimented with a simple genetic algorithm that tries to minimize the maximum hexomino fence area by swapping hexominoes in the order. Figure~\ref{fig:hex-worst} shows the best order found after 6 hours of running time, with a maximum area of $1448$.

\section{Expected length of reference vectors}
\label{sec:refvec-length-estimate}

\begin{figure}[h!]
\centering
\begin{tikzpicture}
\begin{axis}[
    width  = 14cm,
    height = 9cm,
    ymin = 0,
    xmin = 0,
    xmax = 500,
    xlabel = {$n$},
    ylabel = {$\mathbb{E}[|v^{\max}|]$},
    legend pos  = north west,
    legend style = {font=\small},
    grid = both,
    grid style = {line width=0.2pt, draw=gray!30},
    minor tick num = 1,
]

\addplot[
    only marks,
    mark      = *,
    mark size = 1.2pt,
    color     = blue!70,
]
table[
    x=n,
    y=mean_length,
]{\datatable};
\addlegendentry{sample average}

\addplot[
    domain  = 10:500,
    samples = 200,
    thick,
    color   = red,
    dashed,
]
{1.464 * pow(x,0.643) - 0.479};
\addlegendentry{$1.464 \cdot n^{0.643} - 0.479$}
\end{axis}
\end{tikzpicture}
\caption{Monte Carlo estimate of $\mathbb{E}[|v^{\max}|]$ using
  $100\,000$ random $n$-ominoes each for $n = 10, 20, \ldots,\allowbreak 490, 500$.}
\label{fig:refec-fit}
\end{figure}
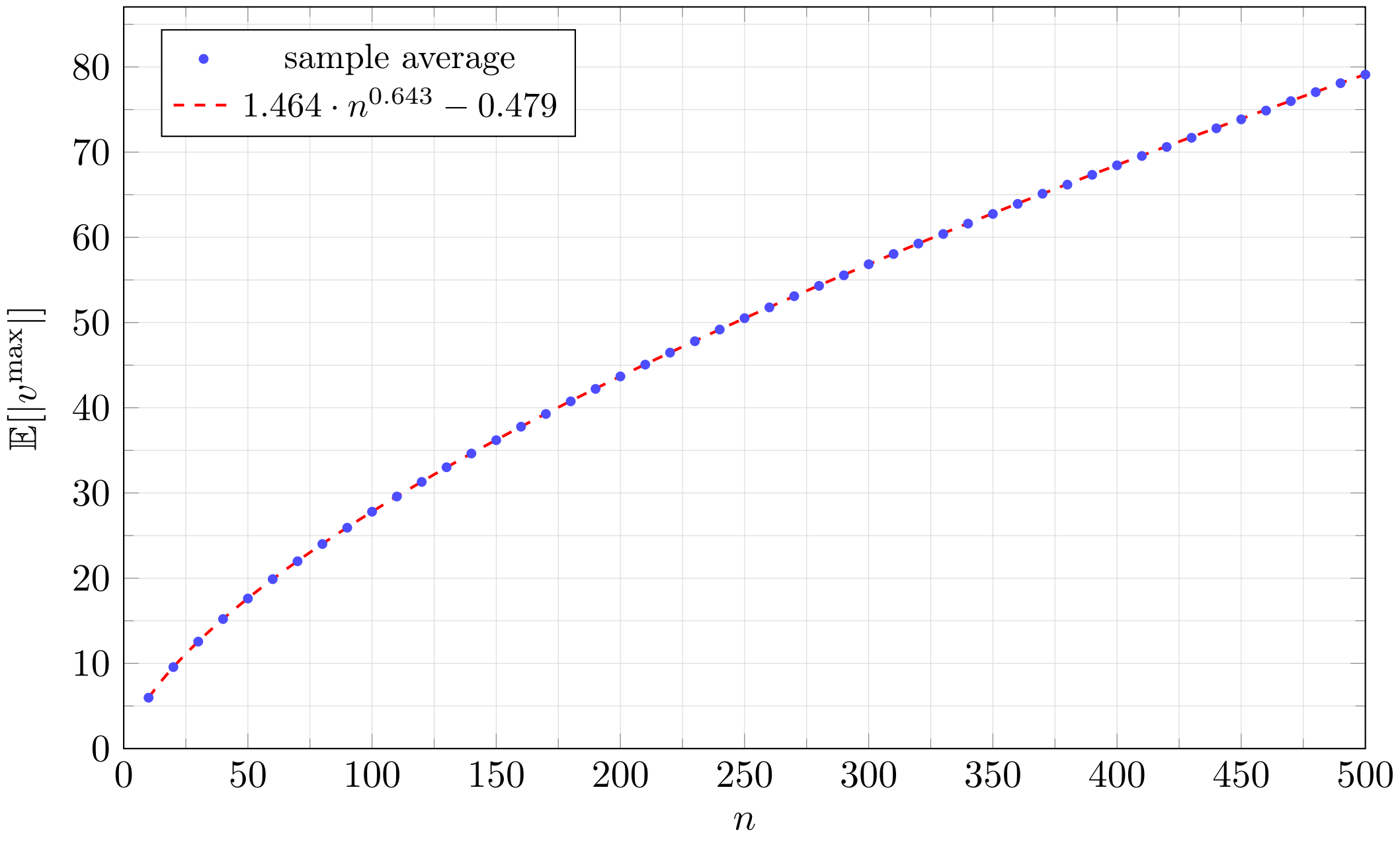

As an extension of Section~\ref{sec:limit}, we now determine the expected length of the longest reference vector $v^{\max}$ of a random $n$-omino, which we denote by $\mathbb{E}[|v^{\max}|]$. We reuse the Metropolis algorithm to sample $100\,000$ random $n$-ominoes for each $n \in \{10, 20, \ldots, 490, 500\}$ and show their average reference vector length in Figure~\ref{fig:refec-fit}.

Fitting a curve of the form $a \cdot n^{b} -c$ to these data, we obtain $1.464 \cdot n^{0.643} - 0.479$. We want to contrast the exponent with a related quantity: Jensen and Guttmann~\cite{Iwan_Jensen_2000} concluded from numerical experiments that the radius of gyration (the root-mean-square distance of the cells from the centroid) of an $n$-omino scales with the exponent $0.64115$. It is remarkable that the exponents are quite close. This means that the square of the diameter, i.e., of the \emph{maximum} distance, does not grow much faster than the \emph{average} squared distance. Indeed, as can be seen in Figure~\ref{fig:500}, typical polyominoes do not fill the area densely, but they grow mainly along linear structures.

\end{document}